\documentclass[11pt,a4paper]{amsart}
\usepackage{amsmath,amssymb,times, amscd, graphicx, xspace,mathrsfs,a4wide,hyperref,url,chngcntr,mathtools,tikz-cd,comment,cleveref}
\usepackage[title]{appendix}
\hypersetup{
  colorlinks   = true, 
  urlcolor     = blue, 
  linkcolor    = blue, 
  citecolor   = blue 
}

\usepackage[OT2,T1]{fontenc}
\DeclareSymbolFont{cyrletters}{OT2}{wncyr}{m}{n}
\DeclareMathSymbol{\Sha}{\mathalpha}{cyrletters}{"58}

\newcommand{\id}{\operatorname{id}}
\newcommand{\C}{\mathbb{C}}
\newcommand{\Z}{\mathbb{Z}}
\newcommand{\Q}{\mathbb{Q}}
\newcommand{\R}{\mathbb{R}}
\newcommand{\N}{\mathbb{N}}
\newcommand{\A}{\mathbb{A}}
\newcommand{\p}{\mathfrak{p}}
\newcommand{\calRSN}{\calR_{\bfS_N}}
\newcommand{\Gal}{\operatorname{Gal}}
\newcommand{\GL}{\operatorname{GL}}
\renewcommand{\O}{\mathscr{O}}
\renewcommand{\d}{\mathrm{d}}
\renewcommand{\mod}{\operatorname{mod}}
\newcommand{\Mod}{\operatorname{Mod}}
\renewcommand{\phi}{\varphi}

\newcommand{\Ind}{\operatorname{Ind}}
\newcommand{\End}{\operatorname{End}}
\newcommand{\coker}{\operatorname{coker}}
\newcommand{\rank}{\operatorname{rank}}
\newcommand{\Hom}{\operatorname{Hom}}
\newcommand{\Tor}{\operatorname{Tor}}
\newcommand{\Frob}{\operatorname{Frob}}
\newcommand{\Spec}{\operatorname{Spec}}
\newcommand{\ad}{\operatorname{ad}}

\newcommand{\surj}{\twoheadrightarrow}

\newcommand{\bt}{\bullet}
\newcommand{\colim}{\operatorname{colim}}

\newcommand{\calDvlambdatau}{\mathscr{D}_{\rhobar_v}^{\lambda_v,\preceq \tau_v}}
\newcommand{\frakm}{\mathfrak{m}}
\newcommand{\opL}{\operatorname{L}}
\newcommand{\AniMod}{\mathscr{Mod}}
\newcommand{\Map}{\operatorname{Map}}
\newcommand{\rhobar}{{\bar{\rho}}}
\newcommand{\Xboxrhobar}{\mathscr{X}_{\rhobar}^\square}
\newcommand{\Xrhobar}{\mathscr{X}_{\rhobar}}
\newcommand{\Xboxrhobarxi}{(\mathscr{X}^\square_{\rhobar})_{(\xi)}}
\newcommand{\ArtE}{\mathscr{Art}_{E}}
\newcommand{\Artk}{\mathscr{Art}_{k}}
\newcommand{\CNLk}{\mathscr{CNL}_{k}}
\newcommand{\CNLE}{\mathscr{CNL}_{E}}
\newcommand{\Int}{\mathscr{Int}}
\newcommand{\otimesL}{\otimes^{\mathbf{L}}}
\newcommand{\Fun}{\mathscr{Fun}}

\newcommand{\RxExihat}{(\calRS\otimes_\O E)_\xi^\wedge}
\newcommand{\IFv}{I_{F_v}}

\newcommand{\fraka}{\mathfrak{a}}
\newcommand{\Sinfahat}{\widehat{S}_{\infty,\fraka}^\circ}

\newcommand{\slashL}{/^{\mathbf{L}}}
\newcommand{\Csc}{C_{\operatorname{sc}}}
\newcommand{\Set}{\mathbf{Set}}
\newcommand{\simto}{\stackrel{\sim}{\to}}
\newcommand{\Sinf}{S_\infty}
\newcommand{\Soinf}{\Sinf^\circ}

\newcommand{\calC}{\mathscr{C}}
\newcommand{\calX}{\mathscr{X}}
\newcommand{\Ani}{\mathscr{Ani}}

\newcommand{\calD}{\mathscr{D}}
\newcommand{\Dv}{\calD_{v}}

\newcommand{\RN}{R_N}
\newcommand{\Rinf}{R_\infty}
\newcommand{\Rinfphat}{\widehat{R}_{\infty,\p}}

\newcommand{\RNbar}{\overline{R}_N}
\newcommand{\SNbar}{\overline{S}_N}
\newcommand{\SNurbar}{\overline{S}^{\operatorname{ur}}_N}
\newcommand{\calS}{\mathscr{S}}

\newcommand{\RMbar}{\overline{R}_M}
\newcommand{\bbT}{\mathbb{T}}
\newcommand{\TS}{\bbT^S}

\newcommand{\Rvlv}{R_{\rho_v}^{\lambda_v}}
\newcommand{\Vlto}{\Vlt^{\circ}}

\newcommand{\calT}{\mathscr{T}}

\newcommand{\bfG}{\mathbf{G}}
\newcommand{\WD}{\operatorname{WD}}
\newcommand{\calR}{\mathscr{R}}
\newcommand{\opT}{\operatorname{T}}
\newcommand{\RG}{R\Gamma}
\newcommand{\Vlam}{\mathscr{V}_\lambda}
\newcommand{\Minfphat}{\widehat{M}_{\infty,\p}}
\newcommand{\bfZ}{\mathbf{Z}}
\newcommand{\calXbox}{\calX^\square}
\newcommand{\Pro}{\operatorname{Pro}}
\newcommand{\Setcp}{\operatorname{Set}^{\operatorname{cp}}}
\newcommand{\dDn}{\partial \Delta^n}

\newcommand{\bfS}{\mathbf{S}}

\newcommand{\calSNur}{\calS_N^{\operatorname{ur}}}
\newcommand{\SN}{S_N}

\newcommand{\calRN}{\calR_N}
\newcommand{\calSN}{\calS_N}
\newcommand{\RSSloc}{R_{\bfS}^{S,\operatorname{loc}}}
\newcommand{\RSprSloc}{R_{\bfS'}^{S,\operatorname{loc}}}
\newcommand{\RSSQN}{R_{\bfS_{Q_N}}^{S}}
\newcommand{\SNo}{S^\circ_N}
\newcommand{\SNobar}{\overline{S}_N^\circ}
\newcommand{\SMobar}{\overline{S}^{\circ}_M}
\newcommand{\calRNpr}{\calRN^\prime}

\newcommand{\GFSQN}{\Gamma_{F,SQ_N}}

\newcommand{\calF}{\mathcal{F}}
\newcommand{\FN}{\calF_N}
\newcommand{\Finfty}{\calF_{\infty}}

\newcommand{\calRS}{\calR_{\bfS}}
\newcommand{\calRSpr}{\calR_{\bfS'}}
\newcommand{\XGammaG}{\mathscr{X}_{\Gamma,\bfG}}
\newcommand{\XGammaGbox}{\XGammaG^\square}
\newcommand{\GFS}{\Gamma_{F,S}}
\newcommand{\GFv}{\Gamma_{F_v}}

\newcommand{\XGammaFSGbox}{\calX_{\GFS,\bfG}^{\square}}

\newcommand{\bfT}{\mathbf{T}}
\newcommand{\uprightAlgO}{\operatorname{Alg}_\O}
\newcommand{\CAlgO}{\mathscr{Alg}_{\O/}}
\newcommand{\Grp}{\operatorname{Grp}}
\newcommand{\bfPG}{\mathbf{PG}}

\newcommand{\calRrho}{\calR_{\rho}}
\newcommand{\Rrho}{R_{\rho}}
\newcommand{\Rrhoss}{R_\rho^{\operatorname{ss}}}
\newcommand{\Drhoss}{D_{\rho}^{\operatorname{ss}}}
\newcommand{\Rlocrho}{R^{\operatorname{loc}}_{\rho}}
\newcommand{\Rlocssrho}{R^{\operatorname{loc},\operatorname{ss}}_{\rho}}
\newcommand{\calRrhoss}{\mathscr{R}_{\rho}^{\operatorname{ss}}}
\newcommand{\TSQN}{\bbT^{S\cup Q_N}}
\newcommand{\TSQNQN}{\TSQN_{Q_N}}
\newcommand{\ON}{\O_N}
\newcommand{\OM}{\O_M}
\newcommand{\Vlt}{\mathscr{V}_{\lambda,\tau}}

\newcommand{\DrhoSss}{\calD_{\rho}^{\operatorname{ss}}}
\newcommand{\Drhobar}{\mathscr{D}_{\rhobar}}
\newcommand{\Dboxrhobar}{\mathscr{D}_{\rhobar}^\square}
\newcommand{\Dboxrhobarv}{\mathscr{D}_{\rhobar_v}^\square}
\newcommand{\rhobarv}{\rhobar_v}
\newcommand{\Dboxrhov}{\calD^{\square}_{\rho_v}}
\newcommand{\Dboxrho}{\calD^{\square}_{\rho}}
\newcommand{\DboxS}{\calD_{\bfS}^\square}
\newcommand{\DS}{\calD_{\bfS}}
\newcommand{\Rvbox}{R_{\rhobarv}^{\square}}
\newcommand{\calDTS}{\calD^{T}_{\bfS}}
\newcommand{\DSxi}{(\DS)_{(\xi)}}
\newcommand{\DSbox}{\DS^{\square}}
\newcommand{\DSboxxi}{(\calD^{\square}_{\bfS})_{(\xi)}}
\newcommand{\Dboxrhobarxi}{(\calD^\square_{\rhobar})_{(\xi)}}
\newcommand{\Dxi}{\calD_{(\xi)}}

\newcommand{\length}{\operatorname{length}}
\newcommand{\Rvlambdatau}{R_{\rhobarv}^{\lambda_v,\preceq\tau_v}}
\newcommand{\Drho}{\calD_\rho}

\newcommand{\GF}{\Gamma_F}
\newcommand{\cL}{\mathcal{L}}
\newcommand{\cLp}{\cL^{\perp}}
\newcommand{\cLv}{\cL_v}
\newcommand{\cLpv}{\cLp_v}
\newcommand{\calPro}{\mathscr{Pro}}

\newcommand{\Ghat}{\widehat{\bfG}}
\newcommand{\PGhat}{\widehat{\mathbf{PG}}}
\newcommand{\PGhatk}{\PGhat_k}
\newcommand{\PGhatE}{\PGhat_E}
\newcommand{\Ghatk}{\Ghat_k}
\newcommand{\GhatE}{\Ghat_E}
\newcommand{\vbar}{\bar{v}}
\newcommand{\wbar}{\bar{w}}
\newcommand{\RHom}{\operatorname{RHom}}

\newcommand{\colimAIntB}{\underset{A\in \Int B}{\operatorname{\colim}}}
\newcommand{\G}{\mathbf{G}}

\newcommand{\AniModcn}{\AniMod^{\operatorname{cn}}}
\newcommand{\RNSSQNpr}{R_{\bfS_{Q_N}'}^S}
\newcommand{\RNSSQN}{R_{\bfS_{Q_N}}^S}
\newcommand{\calA}{\mathcal{A}}

\DeclareFontFamily{U}{mathc}{}
\DeclareFontShape{U}{mathc}{m}{it}%
{<->s*[1.03] mathc10}{}

\DeclareMathAlphabet{\mathscr}{U}{mathc}{m}{it}

\newtheorem{theorem}{Theorem}[subsection]
\newtheorem{proposition}[theorem]{Proposition}
\newtheorem{corollary}[theorem]{Corollary}
\newtheorem{lemma}[theorem]{Lemma}
\newtheorem*{schematic}{Schematic}
\newtheorem*{slogan}{Slogan}

\newtheorem*{conjecture}{Conjecture}

\theoremstyle{definition}
\newtheorem{definition}[theorem]{Definition}

\newtheorem{remark}[theorem]{Remark}

\newtheorem{example}[theorem]{Example}

\crefname{proposition}{Proposition}{Propositions}
\crefname{lemma}{Lemma}{Lemmas}
\crefname{theorem}{Theorem}{Theorems}
\crefname{corollary}{Corollary}{Corollaries}
\crefname{remark}{Remark}{Remarks}
\crefname{definition}{Definition}{Definitions}
\crefname{section}{Section}{Sections}
\crefname{conjecture}{Conjecture}{Conjectures}

\mathchardef\mhyphen="2D
\title{On derived Galois deformation rings in characteristic 0 and a conjecture of Venkatesh}
\author{Douglas Molin}

\address{Department of Mathematical Sciences, Chalmers University of Technology and the University of Gothenburg, 412 96 Gothenburg, Sweden}
\email{molind@chalmers.se}

\input xy
\xyoption{all}
\usepackage{amscd}

\begin{document}

\begin{abstract}
	Let $\Pi$ be a regular algebraic cuspidal automorphic representation of $\GL_n$ over a CM field $F$, and $\rho$ its associated Galois representation in characteristic $0$. With mild assumptions, we use the Taylor--Wiles method to define a free action of the potentially semistable derived deformation ring of $\rho$ on the $\Pi$-part of the $p$-adic cohomology of $\GL_n/F$, confirming a prediction of Venkatesh. As a key step, we establish an equivalence between the derived deformation ring of $\rho$ with the completion at the point corresponding to $\rho$ of the derived deformation ring of the mod~$p$ reduction $\rhobar$, generalising well-known results by Kisin.
\end{abstract}
\maketitle

\tableofcontents

\counterwithin{equation}{subsection}

\section{Introduction}
The Langlands programme predicts correspondences between automorphic forms and Galois representations. Classically, these take the form of bijections between sets of isomorphism classes. In its modern formulation, the programme predicts these correspondences to admit a geometric interpretation, as in the following schematic, inspired by \cite[§4.7]{zhu2025} (see also \cite[Conj.~1.5]{scholze2026}).
\begin{schematic}
    Let $\calA$ be a space of automorphic forms with an action of a Hecke algebra $\mathbb{T}$. Then there exists a moduli stack $\calX$ of Galois representations with a morphism
    $
      \calX \to \Spec(\mathbb{T})
    $
    such that
    \[
      \calA \cong \Gamma(\calX,\omega_{\calX}),
    \]
    the global sections of the dualising complex $\omega_\calX$.
\end{schematic}
The objects are left deliberately unspecified, but the relationships between them are made precise. The morphism maps Galois representations to systems of Hecke eigenvalues, and $\calA$ is pulled back along this morphism to a sheaf on the Galois stack, where it is identified with the dualising complex. Since the dualising complex is generically a trivial line bundle, we may extract the following slogan.
\begin{slogan}
    Locally near a sufficiently generic point, a space of automorphic forms is naturally a free module of rank $1$ over a local ring of functions of a moduli stack of Galois representations.
\end{slogan}
This slogan has appeared in different guises in a series of works by Venkatesh and collaborators on which the present work builds (\cite{galatius_venkatesh2018,prasanna_venkatesh2021,venkatesh2019,harris_venkatesh2019}). To explain this in detail in the $p$-adic setting, let us now fix the meaning of the objects appearing in the slogan. 

Fix a prime $p$, a CM field $F$ with maximal totally real subfield $F^+$, an isomorphism $\iota\colon\C \cong \overline{\Q}_p$ and a finite extension $E/\Q_p$ assumed large enough. Let $\G = \GL_n/F$ and suppose $\Pi$ is a regular algebraic cuspidal automorphic representation of $\G(\A_F)$ of weight $\lambda\in(\Z^n_+)^{\Hom(F,E)}$. Fix an open compact subgroup $K \subset \G(\A^\infty_F)$ and consider the locally symmetric space $X_K$ of level $K$, which is a smooth manifold of dimension $\dim X_K = [F^+\colon \Q]n^2-1$ equipped with a locally constant sheaf of $E$-vector spaces $\Vlam$ determined by the weight $\lambda$. The singular homology $H_\ast(X_K,\Vlam)$ is a finite-dimensional $E$-vector space equipped with an action of spherical Hecke operators, and the system of eigenvalues associated to $\Pi$ defines an eigenspace
\[
    H_\ast(X_K,\mathscr{V}_{\lambda})_\Pi \subset H_\ast(X_K,\mathscr{V}_{\lambda}),
\]
which is non-zero in degrees $[q_0,q_0+l_0]$, where $l_0 = [F^+\colon\Q]n-1$ and $2q_0+l_0 = \dim X_K$. Up to a fixed factor depending on the prime-to-$p$ level, the dimension of $H_{q_0+i}(X_K,\Vlam)_\Pi$ for $i=0,\dots,l_0$ is $\binom{l_0}{i}$. This suggestive numerology and the Bloch--Kato conjecture underlies the following conjecture.
\begin{conjecture}\label[conjecture]{conj_ven}\cite[Conj.~1.2.1]{prasanna_venkatesh2021}
    $H_\ast(X_K,\mathscr{V}_{\lambda})_\Pi$ is a free graded module over $\wedge^\ast V_\Pi$, where $V_\Pi$ is the dual adjoint Bloch--Kato Selmer group attached to $\Pi$ (defined below), a vector space of dimension $l_0$.
\end{conjecture}
Properly speaking, this is the $p$-adic realisation of a conjecture for homology with $\Q$-coefficients where the role of $V_\Pi$ is played by a conjectural motivic cohomology group. Thus, one can think of this action as evincing a `hidden' motivic action. The statement that $\dim V_\Pi = l_0$ essentially follows from the Bloch--Kato conjecture. In order for the conjecture to fit with the slogan above, we ought to explain how the dual adjoint Bloch--Kato Selmer group can be thought of as a `local ring of functions' on a stack of Galois representations.

By work of Harris--Lan--Taylor--Thorne and Scholze \cite{scholze_2015,hltt_2016}, one may attach to $\Pi$ a continuous Galois representation
\[
    \rho \colon \GF = \Gal(\overline{F}/F) \to \GL_n(E),
\]
uniquely characterised up to isomorphism by its values at Frobenius elements matching the Hecke eigenvalues of $\Pi$ in a precise way. The representation $\rho$ is unramified outside a finite set of places $S$ and thus factors through a quotient denoted $\GFS$; to simplify the exposition in the introduction, we ignore the places in $S$ not above $p$.

The representation $\rho$ and its Hodge--Tate weights $\operatorname{HT}(\rho)$ define a (classical) deformation functor ${\Drhoss\colon \operatorname{Art}_E \to \operatorname{Set}}$ where $\operatorname{Art}_E$ is the category of local artinian $E$-algebras, given by
\[
    B \mapsto \{r_B\colon \GFS \to \GL_n(B) \mid r_B \equiv \rho \mod \frakm_B,r_B\text{ potentially semistable, }\operatorname{HT}(r_B)=\operatorname{HT}(\rho)\}/\sim,
\]
where $\sim$ denotes strict equivalence. When $\rho$ is absolutely irreducible, $\Drhoss$ is pro-represented by a complete noetherian local $E$-algebra $\Rrho^{\operatorname{ss}}$, whose spectrum $\Spec(\Rrhoss)$ we view as an intersection
\[
    (\text{deformations of }\rho) \cap (\text{deformations of }\prod_{v\in S} \rho_v\text{ w. local conditions}) \subset (\text{ deformations of }\prod_{v\in S} \rho_v).
\]
The construction of $\rho$ from $\Pi$ furnishes an action
\[
    \Rrhoss \circlearrowright H_\ast(X_K,\Vlam)_\Pi,
\]
and thus, our space of automorphic forms is lifted to a module over a ring describing infinitesimal deformations of $\rho$. As we shall see, $\Rrhoss$ can indeed be thought of as a local ring of functions on a stack of Galois representations.

In terms of rings, the intersection above is represented by a tensor product
\[
    \Rrhoss \cong \Rrho \otimes_{\Rlocrho} \Rlocssrho.
\]
By \cite[Prop.~1.3.12]{allen2016}, the tangent space of $\Rrhoss$ is given by the adjoint Bloch--Kato Selmer group $H^1_g(\GFS, \ad\rho)$, which is known to vanish -- confirming the Bloch--Kato conjecture's prediction -- in many cases by the main result of \cite{a'campo2024}. In other words, $\Spec \Rrhoss$ is a point. However, tangent-obstruction theory predicts the dimension of $\Rrhoss$ to be given by an Euler characteristic which comes out to $-l_0<0$. This suggests that the intersection above is not transverse and that $\Rrhoss$ underlies a `derived object' which encodes this statement. 

A natural candidate is given by replacing $\otimes$ with $\otimesL$, but the resulting object only has good formal properties when the rings involved are complete intersections of the expected dimension, which is expected to hold but hard to verify for $\Rrho$. Supposing for the moment that this was known, we would have an isomorphism
\[
    \Tor_\ast^{\Rlocrho}(\Rrho,\Rlocssrho) \cong \bigwedge H^1_g(\GFS,\ad\rho),
\]
and one might hope to define an action of this $\Tor$-algebra on $H_{\ast}(X_K,\Vlam)_\Pi$. This is the approach taken in \cite{hansen_thorne2017} (using eigenvarieties) and \cite{molin2024} (using completed homology), but one can circumvent the issue entirely by replacing the $\Tor$-algebra with a more natural object.

This is the motivation behind the theory of derived Galois deformation rings developed by Galatius--Venkatesh \cite{galatius_venkatesh2018}, which yields derived enhancements of classical deformation rings for mod~$p$ representations. Later, Zhu \cite{zhu2025} gave a more general construction which also works in characteristic $0$, In fact, Zhu constructs a derived moduli stack of Galois representation as featured in the beginning of this introduction, whose local rings of functions are derived deformation rings. In both setups, a \emph{derived deformation ring} is an animated (or simplicial) ring whose underlying static ring is canonically isomorphic to a classical deformation ring; for example, there exists a derived deformation ring denoted $\calRrhoss$ which satisfies
\[
   \pi_0(\calRrhoss) \cong \Rrhoss 
\] 
The benefit of the derived framework can be explained in terms of tangent complexes. To $\calRrhoss$ and $\Rrhoss$ we can associate a cochain complexes of $E$-vector spaces denoted $\opT^\bullet_E(\calRrhoss)$ and $\opT^\bullet_E(\Rrhoss)$, respectively. In both cases, the $0$th cohomology is the adjoint Bloch--Kato Selmer group:
\[
    \opT^0_E(\calRrhoss) \cong \opT^0_E(\Rrhoss) \cong H^1_g(\GFS, \ad\rho).
\]
For the derived deformation ring, we have
\[
    \opT^1_E(\calRrhoss) \cong H^{2}_g(\GFS, \ad\rho) \cong H^1_g(\GFS, \ad\rho(1))^\vee,
\]
and $\opT^i_E(\calRrhoss)=0$ for $i\geq 2$. In contrast, supposing that $\Rrhoss\cong E$ we have $\opT^i_E(\Rrhoss)=0$ for all $i$. As another example, the unrestricted derived deformation ring of $\rho$ satisfies
\[
    \opT^i_E(\calRrho) \cong H^{i+1}(\GFS, \ad\rho),\ i\geq 0,
\]
whereas there is no clear-cut description of $\opT^i_E(\Rrho)$ in terms of Galois cohomology for $i\geq 1$. Arguably, this closer relation to Galois cohomology means $\calRrho$ is a more natural object to consider; more importantly, it means the derived framework is more amenable to computations, since many questions about derived deformations reduce to the analysis of tangent complexes.

Without further ado, let us state a simplified version of our main result (cf. \Cref{conj_ven} above).
\begin{theorem}[\Cref{theorem_mainresult}]\label{theorem_mainresult_intro}
With notation as above and under technical assumptions on $F$ and $\rho$, there exists a free action of $\pi_\ast(\calRrhoss)$ on $H_\ast(X_K,\mathscr{V}_\lambda)_\Pi$ and a canonical isomorphism $\pi_\ast(\calRrhoss)\cong \wedge V_\Pi$.
\end{theorem}
In \cite{galatius_venkatesh2018}, a free action of the crystalline version of $\pi_\ast(\calRS)$ on the \emph{integral} cohomology is defined, in the Fontaine--Laffaille setting and assuming the Calegari--Geraghty vanishing conjecture for mod~$p$ cohomology \cite[Conj. B]{calegari_geraghty_2018}, together with some other simplifying assumptions (some of which have since been relaxed, see \cite{cai2021,cai_tilouine2022}). Inverting $p$ then yields the desired action. In a general setting, this strategy is not applicable, and the naive generalisation of the integral result is not expected to always hold. Our proof takes a modified approach with the aim of avoiding restrictive local assumptions, at the cost of obtaining only a result in characteristic 0.

The proof of \Cref{theorem_mainresult_intro} builds on the Taylor--Wiles method as applied to derived deformation rings by Galatius--Venkatesh. Since this method concerns mod~$p$ deformations and we are interested in a result in characteristic 0, we relate the derived deformation ring of $\rho$ to that of $\rhobar$ by proving derived analogues of Kisin's widely referenced results \cite[Lem.~2.3.3, Prop.~2.3.5]{kisin2009}, which relate the two notions in the classical setting. For example, we obtain the following result as a consequence of \Cref{prop_FrhobarxiFrho,theorem_charptochar0rings}, which are applicable in broad generality. 
\begin{corollary}\label[corollary]{cor_introcompletiondemo}
    Let $\calR_{\rhobar} \in \CNLk$ and $\calRrho\in \CNLE$ denote the derived deformation rings of $\rhobar$ and $\rho$, respectively, and let $\xi\colon\calR_{\rhobar} \to \O$ be a point corresponding to an $\O$-lattice $\xi \subset \rho$. Then
    \[
        \calRrho \simeq \lim_s\ (\calR_{\rhobar}\otimes_\O E)\slashL I_\xi^s
    \]
    the derived completion of $\calR_{\rhobar} \otimes_\O E$ at $I_\xi = \ker(\pi_0(\calR_{\rhobar})\otimes_\O E \to E)$, is equivalent to $\calRrho$.
\end{corollary}
Using the same results, and under our technical assumptions, we prove an equivalence
\[
    \calRrhoss \simeq (\calRS \otimes_\O E)_\xi^\wedge,
\]
where $\calRS$ represents a mod~$p$ (derived) deformation problem $\bfS$ defined using Kisin's potentially semistable deformation rings (see \Cref{cor_RSphat_equals_calRrhoss}). This equivalence allows us to carry out patching arguments in the mod~$p$ setting and then pass to characteristic 0.

Let us describe in simplified terms the Taylor--Wiles argument. The goal, following \cite{galatius_venkatesh2018}, is to approximate the mod~$p$ deformation ring $\calR \coloneqq \calRS$ by a derived tensor product of static rings. For every $N\geq 1$, one defines an augmented deformation problem $\bfS_N$ of Taylor--Wiles level $N$ and a corresponding derived deformation ring $\calRN \coloneqq \calRSN$ with underlying static ring $\RN \coloneqq \pi_0(\calRN)$. Descending from Taylor--Wiles level yields an equivalence
\[
    \calR \simeq \calRN \otimesL_{\SNo} \O,
\]
where $\SNo$ is the $\O$-group ring of a finite $p$-group. Via this equivalence, one defines a morphism
\[
    \calR \simeq \calRN \otimesL_{\SNo} \O \to \RN \otimesL_{\SNo} \O \to \RNbar \otimesL_{\SNobar} \O_N \eqqcolon \calC_N
\]
where $\RN \surj \RNbar$ and $\SNo \surj \SNobar$ are artinian quotients and $\O_N \coloneqq \O/\varpi^N$. We think of $\calC_N$ as a `mod~$\varpi^N$ approximation' of $\calR$. A compactness argument implies the existence of a compatible system of approximations, yielding a map
\[
    g\colon \calR \to \calC_\infty = \lim_N\ \calC_N
\]
for some sequence of $N$. Localising and completing at the point corresponding to $\xi \subset \rho$ as in \Cref{cor_introcompletiondemo}, we obtain a morphism
\[
    (\calRS\otimes_\O E)_\xi^\wedge \simeq \calRrhoss \to (\calC_\infty \otimes_\O E)_\xi^\wedge.
\]
We prove that this morphism is an equivalence by comparing tangent complexes. Carrying out this non-constructive argument while retaining sufficient control over the induced morphism of tangent complexes is the main technical difficulty in the proof. Applying an `$R=T$' theorem from \cite{a'campo2024}, we identify $(\calC_\infty \otimes_\O E)_\xi^\wedge$ with a derived tensor product of power series rings:
\[
    (\calC_\infty \otimes_\O E)_\xi^\wedge \simeq \Rinfphat \otimesL_{\Sinfahat} E \simeq E[[X_1,\dots,X_s]] \otimesL_{E[[X_1,\dots,X_{s+l_0}]]} E.
\]
A'Campo proves that $\Rinfphat = E[[X_1,\dots,X_s]]$ acts freely on a patched homology module $\Minfphat$, and descending from Taylor--Wiles level, we obtain a free action
\[
    \pi_\ast(\calRrhoss) \circlearrowleft \pi_\ast\big(\Minfphat \otimesL_{\Sinfahat} E\big) \cong H_\ast(X_K,\Vlam)_{\Pi}.
\]
A priori, the action thus obtained depends on the non-canonical choice made in defining the morphism $g$ above, though we expect it to be canonical (see \Cref{rmk_noncanonical}).

\subsection*{Outline}
Let us outline the content of the sections below. In \Cref{section_preliminaries}, we recall generalities concerning the cohomology of locally symmetric spaces and animated rings. \Cref{sec_ddt} discusses the derived deformation theory of Galois representations: we define various derived deformation functors and prove all their necessary properties. In particular, \Cref{subsec_localisationatachar0pt} concerns the link between deformations of a characteristic 0 representation and its mod~$p$ reduction. In \Cref{section_mainresults}, we use the patching and the results of the preceding sections to prove the main result, \Cref{theorem_mainresult}.

\subsection*{Notation and conventions}
The following notation will be used throughout the article. $F$ is a CM field with ring of integers $\O_F$. When $S$ is a set of places of $F$, we use the notation $S_p$ for $\{v\in S\colon v\mid p\}$, and the completion of $\O_F$ at a place $v$ is $\O_{F_v}$. Our field of coefficients is denoted $E$ and is a finite extension of $\Q_p$ with ring of integers $\O$, uniformiser $\varpi$ and residue field $k \coloneqq \O/\varpi$. The extension $E$ is tacitly assumed to be large enough, so that all of our Galois representations are defined over $E$ and $E$ contains all embeddings $F \hookrightarrow \overline{\Q}_p$. For any $r\geq 1$, we let $\O_r \coloneqq \O/\varpi^r$. Galois groups are denoted by $\Gamma$, e.g. $\GFS = \Gal(\overline{F}_S/F)$ and $\GFv = \Gal(\overline{F}_v/F_v)$. The inertia subgroup at a place $v$ of $F$ is denoted $I_{F_v}$. When $\rhobar$ or $\rho$ is a representation of a global Galois group such as $\GFS$, we denote restrictions by subscripts, e.g. $\rhobar_v = \rhobar|_{\GFv}$.   

This article uses the language of quasicategories (from now on, $\infty$-categories). Generally, we use calligraphic letters to denote $\infty$-categories, so that e.g. $\CAlgO$ is the animation of the category $\operatorname{Alg}_\O$ of $\O$-algebras (see \Cref{subsec_animatedalgebras}). When we are interested in both an animated object and its underlying static object, we sometimes distinguish them by fonts; this is always explicitly declared at some point in the text, and if so then e.g. $R=\pi_0(\calR)$. Unless otherwise specified, limits and colimits are taken in the $\infty$-categorical sense. We denote by $\Map(X,Y)$ the mapping space between objects. If $Y$ is an object of an $\infty$-category $\mathscr{C}$, we denote by $\mathscr{C}_{/Y}$ the over-$\infty$-category with objects $(X \to Y)$ \cite[§1.2.9]{lurie_htt}. Given an $\infty$-category $\calC$ and a full subcategory $h'\calC \hookrightarrow h \calC$ of its homotopy category, the pullback of
$
    N(h'\calC) \to N(h \calC) \leftarrow \calC,
$
where $N(\cdot)$ denotes the nerve, is the $\infty$-subcategory of $\calC$ spanned by the vertices lying in $h'\calC$ \cite[Ex.~3.9.5]{cisinski_highercat}.

\subsection*{Acknowledgements}
I am grateful to Christian Johansson for suggesting this problem, generously sharing his knowledge, and providing guidance throughout the project. I would also like to thank Patrick Allen, Dustin Clausen, Wushi Goldring, and Vincent Pilloni for their helpful comments and constructive criticism during my thesis defence. Finally, I thank Håvard Damm-Johnsen for providing feedback on an earlier version of this text.

\subsection*{Statement of AI use} ChatGPT was used to improve the presentation and to assist with literature searches. The results and arguments are the author's own.

\section{Preliminaries}\label{section_preliminaries}
The goal of this section is to describe our general setup and collect some useful results from the literature for easy reference.

\subsection{Locally symmetric spaces}
Let $\bfG = \GL_n/F$. We begin by recalling how automorphic representations of $\bfG(\mathbb{A}_F)$ contribute to the cohomology of locally symmetric spaces, following the exposition in \cite[§2.1]{hevesi2024}. We let $X^{\bfG}$ be the symmetric space of $\operatorname{Res}_{F/\Q}\bfG$ as in \cite[§2]{borel_serre}, which is a homogeneous $\bfG(F\otimes_\Q \R)$-space.

An open compact subgroup $K \subset \bfG(\A^\infty_F)$ is \textit{good} if $K = \prod_{v} K_v$ with $K_v\subseteq \bfG(\O_{F_v})$ and $K$ is neat (see \cite[§2.1.1]{10authorpaper}). If $S$ is a finite set of places of $F$, we write $K = K^S \times K_S \coloneqq \prod_{v\in S} K_v \times \prod_{v\notin S} K_v$. For a good subgroup $K$, we define the \textit{locally symmetric space of level $K$} as
\[
    X_K \coloneqq \bfG(F)\backslash (X^{\bfG} \times \bfG(\A^\infty_F)/K),
\]
a smooth manifold of dimension $[F^+\colon\Q](n^2-1)$. The \textit{defect} of $\bfG$ is the positive integer
\[
    l_0 \coloneqq [F^+\colon\Q]n-1.
\]

Let $\lambda = (\lambda_\iota)_{\iota\in \Hom(F,E)} \in (\Z^n_+)^{\Hom(F,E)}$ be a dominant weight and suppose we are given inertial representations $\tau_v\colon I_{F_v} \to \GL_n(E)$ with open kernel for $v\mid p$. Then \cite{schneider_zink1999} associate to each $\tau_v$ an irreducible smooth $\O$-representation $\sigma^\circ(\tau_v)$ of $\GL_n(\O_{F_v})$, and we let $\sigma^\circ(\tau) = \bigotimes_{v\mid p} \sigma^\circ(\tau_v)$. From the pair $(\lambda,\tau)$ we then obtain a locally constant sheaf of $\O$-modules $\Vlto = \Vlam^\circ \otimes \sigma^\circ(\tau) \in \operatorname{Sh}(X_K)$ where $\Vlam^\circ$ is an $\O$-lattice in the highest weight representation of $(\GL_n)_E^{\Hom(F,E)}$ of highest weight $\lambda$. We define $\sigma(\tau) \coloneqq \sigma^\circ(\tau) \otimes_\O E$ and $\Vlt \coloneqq \Vlto \otimes_\O E$. The complex
\[
    \RG(X_K,\Vlto) \in \calD(\O)
\]
carries an action of the \emph{abstract Hecke algebra}
$
    \bbT^S \coloneqq \mathscr{H}(\bfG(\mathbf{A}^{\infty,S}_F),K^S)
$
of $K^S$-biinvariant functions $\mathbf{A}_F^{\infty,S}\to\O$, where $S$ is a finite set of places such that $K_v = \GL_n(\O_{F_v})$ for $v\notin S$. 

Given a regular algebraic cuspidal representation $\Pi$ of $\bfG(\mathbf{A}^\infty_F)$, we define
\[
    H_\ast(X_K,\Vlt)_\Pi \subseteq H_{\ast}(X_K,\Vlt)
\]
as the eigenspace (over $E$) corresponding to the system of Hecke eigenvalues of $\Pi$, i.e. the subspace on which $\bbT^S$ acts via the character determined by $\Pi$. We write
\[
    \TS\big(\RG(X_K,\Vlto)\big) \subset \End_{\calD(\O)}\big(\RG(X_K,\Vlto)\big)
\]
for the corresponding faithful Hecke algebra. There exists a good subgroup $K_\tau \subseteq K$ such that 
\[
    \RG(X_K,\Vlt)_\Pi = \RG\big(K/K_\tau,\RG(X_{K_\tau},\Vlam)\otimes_E \sigma(\tau)\big),
\]
and since taking invariants is an exact functor in characteristic 0, we have
\[
    H_\ast(X_K,\Vlt)_\Pi \cong \big(H_\ast(X_{K_\tau},\Vlam)_\Pi\big)_{K/K_\tau}.
\]
\begin{proposition}\label[proposition]{prop_spreadingout}\cite[Thm.~2.4.10]{10authorpaper} Let $2q_0 = \dim X_K - l_0$, fix a regular algebraic cuspidal representation $\Pi$ of weight $\lambda$ and inertial type $\tau = (\tau_v)_{v\mid p} = (\operatorname{LL}(\Pi_v)|_{I_{F_v}})_{v\mid p}$. Let $m_\Pi$ denote $\frac{1}{\dim \sigma(\tau)}\dim_\C(\Pi^\infty)^{K^p}$. Then $H_i(X_K,\Vlt)_\Pi = 0$ for $i\notin [q_0,q_0+l_0]$ and for $i=0,\dots,l_0$ we have
    \[
        \dim H_{q_0+i}(X_K,\Vlt)_\Pi = m_\Pi \binom{l_0}{i}.
    \]
\end{proposition}
\subsection{Animated algebras}\label{subsec_animatedalgebras}
In this seciton, we recall the definition of animated algebras, following the presentation in \cite[§5.1]{cesnavicius_scholze2024}. Let $\uprightAlgO$ denote the 1-category of commutative $\O$-algebras, and let
\[
    \CAlgO \coloneqq \Ani(\uprightAlgO)
\]
denote the $\infty$-category of \emph{animated $\O$-algebras}, i.e. the $\infty$-category generated under sifted colimits by the subcategory $\uprightAlgO^{\operatorname{ffg}}\subset \uprightAlgO$ of finite free $\O$-algebras; it is the $\infty$-category underlying the simplicial model category
$
    \operatorname{Fun}(\Delta^{\operatorname{op}},\uprightAlgO)    
$
 of simplicial $\O$-algebras \cite[Cor.~5.5.9.3]{lurie_htt}. 
 
 If $A$ is an animated $\O$-algebra, an \emph{$A$-module} is a module over the underlying $\mathbb{E}_1$-ring of $A$ (see \cite[§2.4]{lurie_dag}). We let $\AniMod_A$ denote the category of $A$-modules. When $A = \pi_0(A)$ is static, we denote by $\Mod_A$ the 1-category of $A$-modules. In this case, the sub-$\infty$-category $\AniMod_A^{\operatorname{cn}}$ of connective modules coincides with $\Ani(\Mod_A)$, and is the $\infty$-category underlying the derived category $D^{\geq 0}(A)$ of connective complexes of $A$-modules, equivalently the simplicial model category $\operatorname{Fun}(\Delta^{\operatorname{op}},\Mod_A)$ of simplicial $A$-modules.

\begin{definition}\cite[Def.~2.5.9]{lurie_dag}
    An animated $\O$-algebra $A\in \CAlgO$ is:
    \begin{enumerate}
        \item[(1)] \emph{Noetherian} if $\pi_0(A)$ is noetherian and $\pi_i(A)$ is a finitely generated $\pi_0(A)$-module for every $i\geq 0$.
        \item[(2)] \emph{Local} if $\pi_0(A)$ is local.
    \end{enumerate}
\end{definition}
The content of the following proposition is found immediately below \cite[Rem.~3.1.4]{lurie_dag}. 
\begin{proposition}\cite{lurie_dag}\label[proposition]{prop_finitelypresented}
    Let $A \in \CAlgO$ and suppose $B$ lies in the smallest subcategory of $\mathscr{Alg}_{A/}$ containing $A[X]$ which is closed under finite colimits. Then $\Map_A(B,-)$ commutes with filtered colimits. 
\end{proposition}
If the assumption of \Cref{prop_finitelypresented} holds, we say that $B$ is a \emph{finitely presented $A$-algebra.}
\begin{proposition}\label[proposition]{prop_almostfinitelypresented}\cite[Prop.~3.1.5]{lurie_dag} Let $A\in \CAlgO$ be noetherian, $j\geq 0$ and suppose $B$ is an $A$-algebra. The following are equivalent: 
    \begin{enumerate}
        \item[(1)] The functor $\Map_A(B,-)$ commutes with filtered colimits when restricted to the subcategory $\mathscr{Alg}_A^{\tau_{\leq j}}$ of $j$-truncated $A$-algebras.
        \item[(2)] The $A$-algebra $\tau_{\leq j} (B)$ is noetherian and $\pi_0(B)$ is a finitely presented $\pi_0(A)$-algebra. 
    \end{enumerate}
\end{proposition}
If the equivalent conditions in \Cref{prop_almostfinitelypresented} hold, we say that $B$ is \emph{of finite presentation to order $j$ as an $A$-algebra}. If the conditions hold for every $j\geq 0$, we say that $B$ is \emph{almost of finite presentation over $A$}. Note that if $A$ is a static noetherian ring and $A \surj B$, then $B$ is almost of finite presentation over $A$ by (2) above.

Given an animated $\O$-algebra $A\in \CAlgO$, an animated $A$-module $M \in \AniMod_\O$ and an element $a \colon \{\ast\} \to A$, we obtain a map $\O[t] \to A$ and define
\[
    M \slashL a \coloneqq M \otimesL_{\O[t]} \O,
\]
where the map $\O[t] \to \O$ corresponds to $0\in \O$. More generally, given $a_1,\dots,a_m\in A$, we define
\[
    M \slashL (a_1,\dots,a_m) \coloneqq M \otimesL_{\O[t_1,\dots,t_m]} \O \simeq M \slashL a_1\slashL \dots \slashL a_m.
\]
There are natural maps 
\[
    M \slashL (a_1^{r+1},\dots,a_m^{r+1}) \to M \slashL (a_1^r,\dots,a_m^r) 
\]
forming an inverse system of animated $\O$-modules. If $M\in \CAlgO$, this is an inverse system of animated $\O$-algebras and the limit depends only on the closed subspace of $\Spec \pi_0(A)$ defined by $(a_1,\dots,a_m)$, as the following result shows. 
\begin{proposition}\label[proposition]{prop_lurie_prohomspaces}\cite[Prop.~6.1.1]{lurie_dag} Let $A\in \CAlgO$ and let $I = (a_1,\dots,a_m) \subset \pi_0(A)$ be a finitely generated ideal. Then, for any $B\in \CAlgO$, the natural map
        \[
            \underset{r}{\colim} \Map_{\CAlgO}(A \slashL (a_1^r,\dots a_m^r),B) \to \Map_{\CAlgO}(A,B) 
        \]
        identifies the left-hand side with the subspace spanned by morphisms $f\colon A \to B$ such that the induced morphism $\pi_0(A) \to \pi_0(B)$ factors through $\pi_0(A)/I^s$ for some $s\geq 1$.
\end{proposition}
\begin{definition}\cite[5.6.1]{cesnavicius_scholze2024}
    Let $A\in \CAlgO$ be noetherian and $I=(a_1,\dots,a_m)\subset \pi_0(A)$. Then $A$ is \emph{derived $I$-adically complete} if the natural map $A \to \lim_r A\slashL (a_1^r,\dots,a_m^r)$ is an equivalence.
\end{definition}
\begin{definition}\label[definition]{definition_CNLk_CNLE}
    A local noetherian animated $\O$-algebra $A$ is \emph{complete} if $A$ is derived $\mathfrak{m}_{\pi_0(A)}$-adically complete. We let $\CNLk \subset (\CAlgO)_{/k}$ and $\CNLE\subset(\CAlgO)_{/E}$ denote the $\infty$-categories of complete noetherian local animated $\O$-algebras with residue field $\pi_0(A)/\frakm_{\pi_0(A)}$ equal to $k$ and $E$, respectively.
\end{definition}
\begin{lemma}\label[lemma]{lemma_derived_complete_pi_complete}\cite[Lem.~3.4.13]{bhatt-scholze2015}
    Let $A\in \AniMod_\O$ and suppose $\pi_i(A)$ is $\varpi$-adically complete in the usual sense for every $i\geq 0$. Then $A$ is derived $\varpi$-adically complete.
\end{lemma}
\subsection{Bloch--Kato Selmer groups}
In this section, we recall the definition of the Bloch--Kato Selmer group mentioned in the introduction, and prove a dimension formula for it. We also define a complex whose cohomology is given by Bloch--Kato Selmer groups, which will be related to derived deformation rings in \Cref{cor_RSphat_equals_calRrhoss}. See \cite{milne} for a complete reference on Galois cohomology.

Let $S$ be a finite set of places of $F$ containing $S_p$ and suppose $\rho\colon\GFS\to \GL_n(E)$ is a continuous representation which is de Rham at all places above $p$. Let $v\in S$. If $v \mid p$, then
\[
	H^1_g(\GFv, \ad\rho) \coloneqq \ker\big(H^1(\GFv, \ad\rho) \to H^1(\GFv, \ad\rho\otimes_{\Q_p}B_{\operatorname{dR}})\big),
\]
where the map on the right is induced by the tensor product and $B_{\operatorname{dR}}$ is Fontaine's ring of de Rham periods. If $v\nmid p$, let $\IFv\subset \GFv$ denote the inertia subgroup and set
\[
	H^1_g(\GFv, \ad\rho) \coloneqq \ker\big(H^1(\GFv, \ad\rho) \to H^1(\IFv, \ad\rho)\big).
\]
where the map is given by restriction. The (geometric) \textit{adjoint Bloch--Kato Selmer group} of $\rho$ is the space of global classes with restrictions lying in the groups defined above, i.e.
\[
	H^1_g(\GFS, \ad\rho) \coloneqq \ker\left( H^1(\GFS, \ad\rho) \to \prod_{v\in S} \frac{H^1(\GFv, \ad\rho)}{H^1_g(\GFv, \ad\rho)}  \right).
\] 
Similarly, we define the dual adjoint Bloch--Kato Selmer group as
\[
	H^1_g(\GFS, \ad\rho(1)) \coloneqq \ker\left( H^1(\GFS, \ad\rho) \to \prod_{v\in S} \frac{H^1(\GFv, \ad\rho(1))}{H^1_g(\GFv, \ad\rho(1))}  \right).
\] 
We use lowercase $h$ to denote the dimension of cohomology groups or subspaces thereof, e.g.
\[
    h^1_g(\GFS, \ad\rho) \coloneqq \dim_E H^1_g(\GFS, \ad\rho).
\]
\begin{proposition}\label[proposition]{prop_h1gdimension}
    Let $\rho\colon\GFS \to \GL_n(E)$ be a continuous representation, such that:
    \begin{enumerate}
        \item[(1)] $\rho$ is absolutely irreducible.
        \item[(2)] For every $v\in S$, the Weil--Deligne representation $\WD(\rho_v)$ is generic (\cite[Def.~1.1.2]{allen2016}).
        \item[(3)] For every $v\in S_p$, $\rho_v$ is potentially semistable with distinct Hodge-Tate weights.
    \end{enumerate}
    Then, with $l_0 = [F^+\colon \Q]n-1$,
    \[
        h^1_g(\GFS, \ad \rho) = h^1_g(\GFS, \ad\rho(1))-l_0. 
    \]
\end{proposition}
\begin{proof}
    By the Greenberg--Wiles formula,
    \begin{multline*}
            h^1_g(\GFS, \ad \rho) - h^1_g(\GFS, \ad\rho(1)) = \\ \sum_{v\in S} \big(h^1_g(\GFv, \ad \rho) - h^0(\GFv, \ad \rho)\big)
            +h^0(\GFS, \ad\rho) - h^0(\GFS, \ad\rho(1))
            -\sum_{v\mid \infty} h^0(\GFv, \ad \rho).
    \end{multline*}
    Assumption (1) implies
    $
        h^0(\GFS, \ad\rho) = 1
    $
    and $h^0(\GFS, \ad\rho(1)) = 0$, and assumption (2) implies that the local terms coming from $v\in S\setminus S_p$ vanish. Thus, we have
    \[
        h^1_g(\GFS, \ad \rho) - h^1_g(\GFS, \ad\rho(1)) = \sum_{v\in S_p} \big(h^1_g(\GFv, \ad \rho) - h^0(\GFv, \ad \rho)\big)
            + 1 - [F^+\colon \Q]n^2.
    \]
    Now, when $\WD(\rho_v)$ is generic, the sum over $v\in S_p$ above equals the sum of total multiplicities of the negative Hodge-Tate weights of $\ad\rho$ \cite[Cor.~3.8.4]{bloch_kato1990}, and by assumption (3) this number is $[F\colon \Q]n(n-1)/2$. Since $l_0 = [F^+\colon \Q]n-1$, the result follows.
\end{proof}

\subsection{Poitou--Tate duality and local conditions}\label[section]{ssec_localconditions}
In this section, we construct complexes with local conditions whose cohomology are given by Selmer groups. Later on, these complexes will be related to tangent complexes of derived deformation rings. In the setting of mod~$p$ derived deformations, the required complexes are constructed in \cite[§B]{galatius_venkatesh2018} (the idea behind the construction was known before, see e.g. \cite[§6]{nekovar2006}). In fact, the authors give a more general construction than we shall describe here; we restrict to the case where the local conditions come from a Selmer system in degree 1. The arguments are essentially the same as those of \cite{galatius_venkatesh2018}.

Throughout this section, we assume that $\rho\colon \GFS\to \GL_n(E)$ is potentially semistable at all places above $p$, and that $S$ contains all places above $p$. We let $C^\bullet(\GFS, \ad\rho)$ denote the usual complex of continuous cochains with cohomology groups given by
\[
    H^i(\GFS, \ad\rho)=Z^i(\GFS, \ad\rho)/B^i(\GFS,\ad\rho).
\]
We make similar definitions for $\GFv$ and $\ad\rho(1)$. Given a \textit{Selmer system} $\cL=(\cLv)_{v\in S}$, i.e. a collection of subspaces
\[
    \cLv \subseteq H^1(\GFv, \ad\rho),
\]
we let $\cLpv \subseteq H^1(\GFv, \ad\rho(1))$ be the orthogonal complement of $\cLv$ with respect to the cup product, and $\cLp = (\cLpv)_{v\in S}$. As in \cite[§B.3]{galatius_venkatesh2018}, $\mathcal{L}$ defines subcomplexes
    $
        C^\bullet_{\cLv}(\GFv, \ad\rho) \subseteq C^\bullet(\GFv, \ad\rho)
    $
    for $v\in S$, given by
    \[
C^j_{\mathcal{L}_v}(\GFv, \ad\rho) =
\begin{cases}
C^0(\GFv, \ad\rho), & j = 0,\\[4pt]
\ker\big(Z^1(\GFv, \ad\rho) \to H^1(\GFv, \ad\rho)/\mathcal{L}_v\big), & j = 1,\\[4pt]
0, & j \geq 2.
\end{cases}
\]
We define
    \[
        C^\bullet_{\mathcal{L}}(\GFS, \ad\rho) \coloneqq \operatorname{fib}\left( C^\bullet(\GFS, \ad\rho)\to \bigoplus_{v\in S}\frac{C^\bullet(\GFv, \ad\rho)}{C^\bullet_{\mathcal{L}_v}(\GFv, \ad\rho)}\right),
    \]
    i.e. the cone given in degree $i$ by
    \[
        C^i_{\cL}(\GFS, \ad\rho) = C^i(\GFS, \ad\rho) \oplus \bigoplus_{v\in S} \frac{C^{i-1}(\GFv, \ad\rho)}{C^{i-1}_{\cLv}(\GFv, \ad\rho)}
    \]
    whose cohomology we denote by
    \[
        H^\ast_{\cL}(\GFS, \ad\rho) \coloneqq H^\ast\big(C^\bullet_{\mathcal{L}}(\GFS, \ad\rho)\big).
    \]
    Note that the degree-wise quotient $C^\bullet(\GFv, \ad\rho)/C^\bullet_{\mathcal{L}_v}(\GFv, \ad\rho)$ is the complex
    \[
        0 \to C^1(\GFv, \ad\rho)/C^1_{\cLv}(\GFv, \ad\rho) \to C^2(\GFv, \ad\rho) \to \dots,
    \]
    whose cohomology groups are given by
    \[
        0, \quad H^1(\GFv, \ad\rho)/\cLv, \quad H^2(\GFv, \ad \rho).
    \]
    It follows that we have a long exact sequence
\[\begin{tikzcd}[row sep=small]
	0 & {H^0_{\cL}(\GFS, \ad\rho)} & {H^0(\GFS, \ad\rho)} & 0 \\
	{} & {H^1_{\cL}(\GFS, \ad\rho)} & {H^1(\GFS, \ad\rho)} & {\bigoplus_{v\in S} {H^1(\GFv, \ad\rho)/\cLv}} \\
	{} & {H^2_{\cL}(\GFS, \ad\rho)} & {H^2(\GFS, \ad\rho)} & {\bigoplus_{v\in S} H^2(\GFv, \ad\rho)} \\
	{} & {H^3_{\cL}(\GFS, \ad\rho).}
	\arrow[from=1-1, to=1-2]
	\arrow[from=1-2, to=1-3]
	\arrow[from=1-3, to=1-4]
	\arrow[from=2-1, to=2-2]
	\arrow[from=2-2, to=2-3]
	\arrow[from=2-3, to=2-4]
	\arrow[from=3-1, to=3-2]
	\arrow[from=3-2, to=3-3]
	\arrow[from=3-3, to=3-4]
	\arrow[from=4-1, to=4-2]
\end{tikzcd}\]
Similarly, we define $H^\ast_{\cLp}(\GFS, \ad\rho(1))$ as the cohomology of the complex
\[
    C^\bullet_{\cLp}(\GFS, \ad\rho(1)) = \operatorname{fib}\left(  C^\bullet(\GFS, \ad\rho(1))\to \bigoplus_{v\in S}\frac{C^\bullet(\GFv, \ad\rho(1))}{C^\bullet_{\cLpv}(\GFv, \ad\rho(1))}\right),
\]
where $C^\bullet_{\cLpv}(\GFv, \ad\rho(1))$ denotes the complex
\[
    0 \to C^0(\GFv, \ad\rho(1)) \to \ker\big(Z^1(\GFv, \ad\rho(1)) \to H^1(\GFv, \ad\rho(1))/\cLpv \big) \to 0 \to \dots.
\]
We now prove a duality theorem for these cohomology groups. The following theorem is a characteristic 0 version of \cite[Thm.~B.1]{galatius_venkatesh2018} in the special case where the local conditions come from a Selmer system in degree 1. 
    \begin{theorem}\label{theorem_HicL_duality}
        With notation as above, there is a natural perfect pairing
        \[
            H^i_{\cL}(\GFS, \ad\rho) \times H^{3-i}_{\cLp}(\GFS, \ad\rho(1)) \to E.
        \]
    \end{theorem}
    \begin{proof}
        In the case $i=0$, we have
        $
            H^0_{\cL}(\GFS, \ad\rho) \cong H^0(\GFS, \ad\rho),
        $
        and, by Poitou--Tate duality,
        \[
            H^3_{\cLp}(\GFS, \ad\rho(1)) \cong \coker\bigg( H^2(\GFS, \ad\rho(1)) \to \bigoplus_{v\in S} H^2(\GFv, \ad\rho(1)) \bigg) \cong H^0(\GFS, \ad\rho)^\vee,
        \]
        see e.g. \cite[§6, (8.6.10)]{neukirch2000} for a statement with torsion coefficients, from which the characteristic 0 statement follows from passage to the limit. It remains to verify the statement for $i=1$, i.e. to construct a perfect pairing
        \[
            H^1_{\cL}(\GFS, \ad\rho) \times H^2_{\cLp}(\GFS, \ad\rho(1)) \to E.
        \]
        We follow the argument of \cite[§B.5]{galatius_venkatesh2018}. By the long exact sequence above,
        \[
            H^1_{\cL}(\GFS, \ad\rho) \cong \ker\big(Z^1(\GFS, \ad\rho) \to \prod_{v\in S} H^1(\GFv, \ad\rho)/\cLv\big),
        \]
        whereas an element of $H^2_{\cLp}(\GFS, \ad\rho(1))$ is represented by a tuple
        \[
            (x',(y'_v)_{v\in S}) \in Z^2(\GFS, \ad\rho(1)) \times \prod_{v\in S} C^1(\GFv, \ad\rho(1))
        \]
        such that $x'_v + \d y'_v \in \cLpv$ for every $v\in S$, where  $(\cdot)_v$ denotes the restriction map at $v$, and $\d$ is the differential. Given $x$ and $(x',(y'_v)_{v\in S})$ as above, we choose an arbitrary $z \in C^2(\GFS, E(1))$ such that $\d z = x \smallsmile x'$ and obtain elements
        \[
            -x_v\smallsmile y'_v + z_v \in C^2(\GFv, E(1)).
        \]
        One verifies that this defines a cocycle, and hence we obtain elements of $H^2(\GFv, E(1)) \cong E$. Taking the sum over $v$ produces an element of $E$, which we define to be $\langle x, (x',(y'_v)_{v\in S}) \rangle$.
        
        This pairing is independent of choice of $z$ since $\sum [z_v] = 0$. Moreover, replacing $y'_v$ by $y'_v+w'_v$ for some $w'_v\in C^1_{\cLp}(\GFv, \ad\rho(1))$ changes the local contribution at $v$ by
        $
            x_v \smallsmile w'_v,
        $
        which is a coboundary since $x_v$ and $w'_v$ define orthogonal classes. Lastly, if we replace $x$ by $x-\d a$, one may choose $z-a\smallsmile x'$ in place of $z$, and similarly for $x'$. Thus, $\langle\cdot,\cdot\rangle$ is a well-defined pairing.

        It remains to prove non-degeneracy. Consider the following two-part filtrations
        \begin{align*}
            \Sha^1(\GFS, \ad\rho) \coloneqq \{[x] \in H^1(\GFS, \ad\rho)\mid \forall v\in S\colon [x_v]=0\ \} &\subseteq H^1_{\cL}(\GFS, \ad\rho), \\
            Q^1_{\cLp} \coloneqq \coker\big(H^1(\GFS, \ad\rho(1)) \to \prod_{v\in S} H^1(\GFS, \ad\rho(1))/\cLpv\big) &\subseteq H^1_{\cLp}(\GFS, \ad\rho(1)).
        \end{align*}
        The subspace $Q^1_{\cLp}$ is identified with the classes represented by tuples $(x,(y'_v)_{v\in S})$ such that
        \[
            x'_v = x'_v+\d y'_v = 0
        \]
        for all $v \in S$. It follows that the pairing respects this filtration; i.e. the restriction
        \[
            \Sha^1(\GFS, \ad\rho) \times Q^1 \to E
        \]
        is identically 0. Moreover, we have an isomorphism
        \[
            H^1_{\cLp}(\GFS, \ad\rho(1))/Q^1 \cong \Sha^2(\GFS, \ad\rho(1)).
        \]
        The problem thus reduces to verifying the non-degeneracy for the induced pairings
        \begin{align*}
            \Sha^1(\GFS, \ad\rho) \times \Sha^2(\GFS, \ad\rho(1)) &\to E,\\
            H^1_{\cL}(\GFS, \ad\rho)/\Sha^1(\GFS, \ad\rho) \times Q^1_{\cLp} &\to E.
        \end{align*}
        In the first case, the obtained pairing is the Tate pairing (\cite[Chap.~VIII, Thm.~8.6.7]{neukirch2000}). In the second case, we choose a representative $(x',(y'_v)_{v\in S}) = (0,(y'_v)_{v\in S})$ for the right-hand class such that $\d y'_v = 0$, and $z=0$. The local contribution at $v$ of the pairing reduces to $-x_v\smallsmile y'_v$, which is the usual local pairing between in degree $1$ (up to a minus sign). 
    \end{proof}
    We record the following special case which will be referenced later in \Cref{cor_RSphat_equals_calRrhoss}.
    \begin{corollary}\label[corollary]{cor_localconditions}
        Let $
        \cL = \big(H^1_g(\GFv, \ad\rho)\big)_{v\in S_p} \cup \big(H^1(\GFv, \ad \rho)\big)_{v\in S\setminus S_p}  
        $, and suppose $\WD(\rho_v)$ is generic for $v\in S_p$. Then, with notation as above, there is a natural isomorphism
        \[
            H^i_{\cL}(\GFS, \ad\rho) \cong H^{3-i}_{\cLp}(\GFS, \ad\rho(1))^\vee.
        \]
        In particular, 
        $
            H^2_{\cL}(\GFS, \ad\rho) \cong H^1_g(\GFS, \ad\rho(1))^\vee.
        $
    \end{corollary}
    \begin{proof}
        The first statement follows from \Cref{theorem_HicL_duality}. Now, if $v\in S_p$ and $\WD(\rho_v)$ is generic, the inclusion
        \[
            H^1_f(\GFS, \ad\rho) \subseteq H^1_g(\GFS, \ad\rho)
        \]
        is an equality by \cite[Lem.~1.1.7]{allen2016}, where $H^1_f$ denotes the Bloch--Kato Selmer group defined using the crystalline period ring $B_{\operatorname{cris}}$ in place of $B_{\operatorname{dR}}$. Since the orthogonal of $H^1_f(\GFS, \ad\rho)$ is $H^1_f(\GFS, \ad\rho(1))$ (see \cite[Prop.~3.8]{bloch_kato1990}), the second statement follows.
    \end{proof}

    \section{Derived deformation theory}\label[section]{sec_ddt}
\subsection{Artinian algebras}
In this section, we recall the notion of artinian animated algebras from \cite{lurie_dag}.
\begin{definition}\label[definition]{def_Artk}
    The category $\Artk$ (resp. $\ArtE$) of \emph{artinian animated $\O$-algebras over $k$} (resp. over $E$) is defined as the sub-$\infty$-category of $(\CAlgO)_{/k}$ (resp. $(\CAlgO)_{/E}$) spanned by $(B\to k)$ (resp. $(B\to E)$) such that
    \begin{enumerate}
        \item[(1)] $\pi_0(B)$ is an artinian ring with residue field $k$ (resp. $E$).
        \item[(2)] $\pi_\ast(B)$ is a finitely generated $\pi_0(B)$-module.
    \end{enumerate}
\end{definition}
Any artinian animated $\O$-algebra is truncated. Note that if $B\in \ArtE$ then $\pi_0(B)$ is an $E$-algebra and $\pi_\ast(B)$ is a finite-dimensional $E$-vector space.
\begin{proposition}\label[proposition]{prop_artinianhomspaces}
    Let $A\in \CNLk$ and $B\in \Artk$. Then, for any proper ideal $I\subset \pi_0(A)$, the natural morphism
    \[
        \underset{s}{\colim} \Map_{/k}(A\slashL I^s,B) \to \Map_{/k}(A,B) 
    \]
    is an equivalence.
\end{proposition}
\begin{proof}
    The maximal ideal of $\pi_0(B)$ is nilpotent, hence given any map $A\to B$, the induced morphism $\pi_0(A)\to \pi_0(B)$ factors through a power of $I\subset \pi_0(A)$. We conclude using \Cref{prop_lurie_prohomspaces}. 
\end{proof}
A useful description of artinian algebras allowing for inductive proofs is stated in terms of small extensions.
\begin{definition}\cite[Def.~3.3.1]{lurie_dag}
    Let $B \in \CAlgO$ and $M \in \AniMod_B$ a connective $B$-module. A \emph{small extension} of $B$ by $M$ is a pullback
\[\begin{tikzcd}
	{B'} & B \\
	B & {B\oplus M[1]}
	\arrow[from=1-1, to=1-2]
	\arrow[from=1-1, to=2-1]
	\arrow["s", from=1-2, to=2-2]
	\arrow["{(\id,0)}"', from=2-1, to=2-2]
\end{tikzcd}\]
where $B \oplus M[1]$ is the trivial square-zero extension (see below \cite[Prop.~3.1.5]{lurie_dag}) and $s$ is a section of the projection $B\oplus M[1] \to B.$
\end{definition}
\begin{lemma}\cite[Lem.~6.2.6]{lurie_dag}\label[lemma]{lemma_ArtEsmallextensions}
    Let $(B\to E)\in \ArtE$. Then there exists a factorisation
    \[
        B \simeq B_m \to B_{m-1} \to \dots \to B_0 \simeq E
    \]
    in $\ArtE$ such that $B_j$ is a small extension of $B_{j-1}$ by $E[n_j]$ for an increasing sequence $n_j\geq 0$ and $\dim_E(\pi_\ast(B_j))=j+1$. The analogous statement is true for $(B\to k)\in \Artk$.
\end{lemma}

\begin{lemma}\label[lemma]{lemma_artinian_mapping_spaces_finite}
    If $\calR\in\CNLk$ and $B\in\Artk$, then $\Map_{/k}(\calR,B)$ has finitely many connected components and finite homotopy groups.
\end{lemma}
\begin{proof}
    We argue by induction on $m$, the length of a presentation of $B$ in terms of small extensions as in \Cref{lemma_ArtEsmallextensions}. For $m=0$, we have $B\simeq k$, and $\Map_{/k}(\calR,k)$ is contractible. For the induction step, suppose the claim holds for every $B'\in\Artk$ given by series of length $m'<m$. By \Cref{lemma_ArtEsmallextensions}, there is an algebra $B'\in\Artk$ and an equivalence
    \[
        B\simeq B'\times_{B'\oplus k[n+1]}B'
    \]
    for some $n\geq 0$, where $B'$ is given by a series of length $m'< m$. Applying $\Map_{/k}(\calR,-)$ gives a pullback
    \[
        \Map_{/k}(\calR,B)\simeq\Map_{/k}(\calR,B')\times_{\Map_{/k}(\calR,B'\oplus k[n+1])}\Map_{/k}(\calR,B')
    \]
    and an associated long exact sequence of homotopy groups
    \[
        \cdots\to\pi_{i+1}\big(\Map_{/k}(\calR,B'\oplus k[n+1])\big)\to\pi_i\big(\Map_{/k}(\calR,B)\big)\to \big(\pi_i\big(\Map_{/k}(\calR,B')\big)\big)^2\to\cdots.
    \]
    By the induction hypothesis, $\pi_i\big(\Map_{/k}(\calR,B')\big)$ is finite. To complete the proof, it therefore suffices to prove the same statement for
    $\pi_{i+1}\big(\Map_{/k}(\calR,B'\oplus k[n+1])\big)$. Now,
    we have a map
    \[
        \Map_{/k}(\calR,B'\oplus k[n+1])\longrightarrow\Map_{/k}(\calR,B')
    \]
    induced by the projection $B'\oplus k[n+1]\to B'$, whose fibre over a map $f$ equals
    \[
        \Map_{/B'}(\calR,B'\oplus k[n+1]) \simeq \Map_{/k}(\calR,k\oplus k[n+1]).
    \]
    In the notation of \Cref{section_tangentcomplexes} below, the homotopy groups of the fibre are given by 
    \[
        \pi_i\big(\Map_{/k}(\calR,k\oplus k[n+1])\big)\cong\opT_k^{n+1-i}(\calR),
    \]
    which is finite by \Cref{prop_Tifinitelygenerated} below and the noetherianity assumption on $\calR$. The long exact sequence of homotopy groups
    \[
       \dots \to \pi_i\big(\Map_{/k}(\calR,k\oplus k[n+1])\big)
       \to \pi_i\big(\Map_{/k}(\calR,B'\oplus k[n+1])\big)
       \to \pi_i\big(\Map_{/k}(\calR,B')\big)
       \to \dots
    \]
    now implies the sought statement, which completes the induction step.
\end{proof}
\subsection{Integral models of artinian \texorpdfstring{$E$}{E}-algebras}
In this section, we introduce a notion of integral models of artinian $E$-algebras. The definition is made in analogy with \cite[2.3.5]{kisin2009} and will be used in \Cref{subsec_localisationatachar0pt}. Note that we have an equivalence of $\infty$-categories 
\[
    ((\CAlgO)_{/E})_{/(B\to E)} \simeq (\CAlgO)_{/B}.
\]
\begin{definition}\label[definition]{def_IntB}
    Let $(B\to E)\in \ArtE$. The $\infty$-category $\Int B$ of \emph{integral models of $B$} is defined as the sub-$\infty$-category of $(\CAlgO)_{/B}$ spanned by all $(A\to B)$ such that:
    \begin{enumerate}
        \item[(1)] $\pi_\ast(A)$ is a finitely generated $\O$-module.
        \item[(2)] The induced map $\pi_\ast(A)\to\pi_\ast(B)$ is injective.
        \item[(3)] The induced map $A\otimesL_\O E \to B \otimesL_\O E$ is an equivalence.  
    \end{enumerate}
\end{definition}
By \Cref{prop_almostfinitelypresented}(2) and \Cref{def_IntB}(1), any $A\in \Int B$ is derived $\varpi$-adically complete and almost of finite presentation over $\O$. When $B$ is clear from context, we sometimes denote objects of $\Int B$ simply by $A$. The key property we shall need of the category of integral models is the following. 
\begin{proposition}\label[proposition]{prop_IntBfiltered}
    For any $B\in \ArtE$, the $\infty$-category $\Int B$ is filtered.
\end{proposition}
Before embarking on the proof of \Cref{prop_IntBfiltered}, we establish some preliminary results about integral models.
\begin{lemma}\label[lemma]{lemma_AmapstoBtimesO}
    Let $(B\to E)\in \ArtE$ and $(A\to B)\in \Int B$. The natural map
    \[
        \Map_{\CAlgO}(A,B\times_E \O) \to \Map_{\CAlgO}(A,B)
    \]
    is an equivalence.
\end{lemma}
\begin{proof}
    Since $\O$ and $E$ are static and $\pi_0(A)$ is a finitely generated $\O$-module, we have equivalences of mapping spaces
    \[
    \Map_{\CAlgO}(A,E) \simeq \Map_{\CAlgO}(\pi_0(A),E) \simeq \Map_{\CAlgO}(\pi_0(A),\O)\simeq \Map_{\CAlgO}(A,\O).    
    \]
    Consequently,
\[
\Map_{\CAlgO}(A,B\times_E \O) \simeq \Map_{\CAlgO}(A,B) \times_{\Map_{\CAlgO}(A,E)} \Map_{\CAlgO}(A,\O) \simeq
    \Map_{\CAlgO}(A,B).
\]
\end{proof}
To prove \Cref{prop_IntBfiltered}, we will show that integral models can be constructed via small extensions, in the following sense.
\begin{lemma}\label[lemma]{lemma_IntBsmallextensions}
    Let $B\in \ArtE$ be given by a sequence of small extensions
    \[
        B \simeq B_m \to B_{m-1} \to \dots \to B_0 \simeq E
    \]
    as in \Cref{lemma_ArtEsmallextensions}. Then, for any $A\in \Int B$, the morphism $A\to \O$ factors as 
    \[
        A \simeq A_m \to A_{m-1} \to \dots \to A_0 \simeq \O,
    \]
    where, for every $j\geq 1$, $A_j\in \Int B_j$ is a small extension of $A_{j-1}$ by $\varpi^{-r_j}\O[n_j]$ for a non-decreasing sequence $0\leq n_1 \leq \dots \leq n_m$ and some $r_j\geq 0$.
\end{lemma}
Before embarking on the proof of the Lemma, let us describe the homotopy groups of a small extension. Let $A \in \CAlgO$ be $n$-truncated and suppose $M\in \AniMod_A$ has homotopy concentrated in degree $n-1$, so that $M[1]$ has homotopy concentrated in degree $n$. Given a small extension $A' = A \times_{A\oplus M[1]} A$, we obtain a long exact sequence
\[
    \dots \to \pi_{i+1}(A) \to \pi_{i-1}(M) \to \pi_i(A') \to \pi_i(A) \to \pi_{i-2}(M)\to\dots
\]
which implies that $\pi_i(A')\cong \pi_i(A)$ for every $i\neq n$ and that
\[
    \ker(\pi_n(A')\to \pi_n(A)) \cong \pi_{n-1}(M).
\]
Thus, a single small extension of the form considered in \Cref{lemma_IntBsmallextensions} only changes the homotopy groups in the top degree by adding a module which is free of rank $1$ over $\O$.
\begin{proof}[Proof of \Cref{lemma_IntBsmallextensions}]
    Suppose we have a factorisation $B\simeq B_m \to \dots \to B_0\simeq E$. We prove the lemma by induction on $m$. The case $m=0$ is trivial since $\Int E \simeq \{\O\}$. For the induction step, suppose the statement holds for any length $m-1$ sequence $B'_{m-1}\to\dots \to B'_0\simeq E$. Let $n=n_m$ denote the largest integer such that $\pi_{n}(B)\neq 0$, and let $A\in \Int B_m$. The kernel
    \[
        I = \ker\big(\pi_\ast(A) \hookrightarrow \pi_\ast(B_m) \to \pi_\ast(B_{m-1})\big) \subset \pi_{n}(A)
    \]
    is a $\pi_0(A)$-submodule of $\pi_\ast(A)$, which is free of rank $1$ as an $\O$-module. Using \cite[Prop.~3.3.3, Prop.~3.3.5]{lurie_dag}, we construct a morphism
    $
        A \to A_{m-1}
    $ in $\CAlgO$, where $A_{m-1}$ (`$A/I$' in the notation of \emph{loc. cit.}) satisfies the following properties:
    \begin{enumerate}
        \item[(1)] For any $C\in \CAlgO$, the map $A\to A_{m-1}$ induces an identification
    \[
        \Map_{\CAlgO}(A_{m-1},C) \hookrightarrow \Map_{\CAlgO}(A,C)
    \]
    with the subspace spanned by maps $A\to C$ for which $I\subset \ker(\pi_\ast(A)\to\pi_\ast(C))$.
    \item[(2)] The induced map $\pi_\ast(A)\to \pi_\ast(A_{m-1})$ is an isomorphism in degrees $i\neq n$ and, in degree $i=n$, is surjective with kernel $I$.
    \item[(3)] $A$ is a small extension of $A_{m-1}$ by $I[n]$, i.e. there is an equivalence
    \[
        A \simeq A_m = A_{m-1} \times_{A_{m-1}\oplus I[n+1]} A_{m-1}.
    \]
    \end{enumerate}
    We deduce that the map $A_{m-1}\to B_{m-1}$ defines an object of $\Int B_{m-1}$ and that the square
\[\begin{tikzcd}
	A & {A_{m-1}} \\
	B & {B_{m-1}}
	\arrow[from=1-1, to=1-2]
	\arrow[from=1-1, to=2-1]
	\arrow[from=1-2, to=2-2]
	\arrow[from=2-1, to=2-2]
\end{tikzcd}\]
commutes up to homotopy. By the induction hypothesis, $A_{m-1}\to \O$ factors as a series of small extensions of the required form. This completes the proof.
\end{proof}

\begin{lemma}\label[lemma]{lemma_intBBtimesO}
    Let $B\in \ArtE$. Then $B\times_E \O$ is equivalent in $(\CAlgO)_{/\O}$ to a filtered colimit of algebras of the form of $A_m$ in \Cref{lemma_IntBsmallextensions}.
\end{lemma}
\begin{proof}    
To prove the second claim, we use induction on $m$ again. The case $m=0$ is trivial since $E\times_E \O \simeq E$. For the induction step, suppose the statement is true for any $B'$ given by a length $m-1$ sequence of small extensions. Let $n=n_m+1$ (notation as in \Cref{lemma_ArtEsmallextensions}) and note that
\[
    E[n] = \underset{r\geq 0}{\colim}\ \varpi^{-r}\O[n].
\]
Since filtered colimits commute with finite limits, under the induction hypothesis we have
\begin{equation*}
    \begin{split}
        B_m \times_E \O 
        &\simeq (B_{m-1} \times_{B_{m-1}\oplus E[n]} B_{m-1})\times_E \O \\
        &\simeq (B_{m-1}\times_E \O) \times_{(B_{m-1}\times_E \O)\oplus E[n]}(B_{m-1}\times_E \O)\\
        &\simeq \underset{A_{m-1}}{\colim}\ (A_{m-1} \times_{A_{m-1}\oplus E[n]} A_{m-1})
    \end{split}
\end{equation*}
where the colimit runs over a filtered diagram of $A_{m-1}\in \Int B_{m-1}$ of the sought form. We have equivalences
\begin{equation*}
    \begin{split}
        A_{m-1} \oplus E[n] &\simeq A_{m-1}\oplus \underset{r\geq 0}{\colim}\ \varpi^{-r}\O[n] \\
        &\simeq \underset{r\geq 0}{\colim}\ A_{m-1} \oplus \varpi^{-r}\O[n].
    \end{split}
\end{equation*}
Indeed, $A_{m-1}$ preserves $\O[n]\subset E[n]$, since, if $\mathscr{End}$ denotes the derived endomorphism ring, we have
\[
    A_{m-1} \to \mathscr{End}_{\O}(E[n]) \simeq \pi_0\big(\mathscr{End}_{\O}(E)\big) \cong E,
\]
and hence the action map $A_{m-1} \to E$ factors through $\pi_0(A_{m-1})$ and lands in $\O$. It follows that
\[
    \underset{A_{m-1}}{\colim}\ (A_{m-1} \times_{A_{m-1}\oplus E[n]} A_{m-1}) \cong \underset{A_{m-1}}{\colim}\ \underset{r\geq 0}{\colim}\ (A_{m-1}\times_{A_{m-1}\oplus \varpi^{-r}\O[n]} A_{m-1}),
\]
and thus $B_m$ is a filtered colimit of the sought form.
\end{proof}
We are now in a position to prove the main result of this section, \Cref{prop_IntBfiltered}.
\begin{proof}[Proof of \Cref{prop_IntBfiltered}]
    By the local characterisation of filtered $\infty$-categories \cite[Tag 02PS]{kerodon}, it suffices to prove the following statement: Given $A,A'\in \Int B$, and a morphism of anima
    \[
        \gamma\colon \partial \Delta^n \to \Map_{/B}(A,A')
    \]
    for some $n\geq 0$, there exists an $A''\in \Int B$ and a morphism $g\colon A'\to A''$ such that the composition
    \[
        g \circ \gamma\colon \partial\Delta^n \to \Map_{/B}(A,A') \to \Map_{/B} (A,A'')
    \]
    is null-homotopic. In other words, we would like to define the dotted arrow in the diagram
\[\begin{tikzcd}
	{\partial\Delta^n} & {\Map_{/B} (A,A')} \\
	{\Delta^n} & {\Map_{/B} (A,A'')}
	\arrow["\gamma", from=1-1, to=1-2]
	\arrow[from=1-1, to=2-1]
	\arrow["{g_\ast}", from=1-2, to=2-2]
	\arrow[dotted, from=2-1, to=2-2]
\end{tikzcd}\]
Since $B$ is final in $(\CAlgO)_{/B}$ and by \Cref{lemma_intBBtimesO}, we have
\[
    \{\ast\} \simeq \Map_{/B}(A,B) \simeq \Map_{/B}(A,B\times_E \O) \simeq  \Map_{/B}(A,\underset{A''}{\colim}\ A''),
\]
where $A''$ runs over a filtered subcategory of $\Int B$. Thus, the composition
\[
    \dDn \to \Map_{/B}(A,A') \to \Map_{/B}(A,\underset{A''}{\colim}\ A'')
\]
is null-homotopic. Now, since $A'\in \Int B$ is a finitely presented $\O$-algebra by \Cref{prop_almostfinitelypresented}, the map
$
    A' \to \colim A''
$
factors through some $A''$. It follows that there exists a $g\colon A'\to A''$ such that the composition $g\circ \gamma$ is null-homotopic, which completes the proof.
\end{proof}

\begin{lemma}\label[lemma]{lemma_IntB}
    Let $(B\to E)\in \ArtE$. We have the following:
    \begin{enumerate}
        \item[(1)] There is an equivalence
        \[
            \underset{A\in \Int B}{\colim} A \simeq B \times_E \O.
        \]
        \item[(2)] There is a canonical isomorphism of graded rings
        $
            \pi_\ast(B\times_E \O) \cong \pi_\ast(B) \times_E \O.
        $
        \item[(3)] Any $A\in \Int B$ is derived $\varpi$-adically complete.
        \item[(4)] For any $A\in \Int B$ and $r\geq 1$, $(A\slashL \varpi^r \to k) \in \Artk$. 
    \end{enumerate}
\end{lemma}
\begin{proof}
    (1) By \Cref{lemma_AmapstoBtimesO}, any $(A\to B)\in \Int B$ defines a map $A\to B\times_E \O$. Taking the colimit over $\Int B$, we produce
    \[
        \underset{A\in \Int B}{\colim} A \to B \times_E \O,
    \]
    and we claim that this map is an equivalence. Since $\Int B$ is filtered (\Cref{prop_IntBfiltered}) and filtered colimits commute with $\pi_\ast$, it is enough to prove that the map
    \[
        \underset{A\in \Int B}{\colim} \pi_\ast(A) \to \pi_\ast(B \times_E \O)
    \]
    is an isomorphism. Since $\pi_\ast(A)\to\pi_\ast(B)$ is injective for every $A\in \Int B$ (\Cref{def_IntB}(3)), injectivity follows. Surjectivity follows from the existence of an isomorphism
    \[
        \pi_\ast(B \times_E \O) \cong \underset{A''}{\colim}\  \pi_\ast(A''),
    \]
    provided by \Cref{lemma_intBBtimesO} ($A''$ runs over a filtered subcategory of $\Int B$).

    (2) There is a Mayer-Vietoris sequence
    \[ 
    \dots \to \pi_n(B\times_E \O) \to \pi_n(B) \times \pi_n(\O)\to \pi_{n}(E) \to \dots,
    \]
    and hence for every $n\geq 1$, we have
    $
        \pi_n(B\times_E \O) \cong \pi_n(B)
    $
    since $\pi_n(\O)=\pi_n(E)=0$. For $n=0$, the sequence implies
    $\pi_0(B \times_E \O) \cong \ker(\pi_0(B)\times \O \to E) \cong \pi_0(B) \times_E \O$.

    (3) Since $A\in \Int B$, $\pi_\ast(A)$ is a finitely generated $\O$-module and hence $\varpi$-adically complete. Thus, $A$ is derived $\varpi$-adically complete by \Cref{lemma_derived_complete_pi_complete}.

    (4) Let $A\in \Int B$ and $r\geq 1$. We have an isomorphism $\pi_0(A\slashL \varpi^r)\cong \pi_0(A)/\varpi^r$, which is artinian since $\pi_0(A)$ is a finitely generated $\O$-module by definition. Moreover, there exists a spectral sequence \cite[Thm.~2.6]{quillen_ha}
    \[
        E^2_{i,j}\colon \Tor_{i}^{\O[t]}\big(\pi_j(A),\O\big) \implies \pi_{i+j}(A\slashL \varpi^r).
    \]
    Thus, for any $n\geq 1$, the spectral sequence induces a short exact sequence
    \[
        0 \to \pi_n(A)/\varpi^r \to \pi_n(A\slashL \varpi^r) \to \pi_{n-1}(A)[\varpi^r]\to 0.
    \]
    Since $A$ is truncated, the short exact sequence implies that $\pi_\ast(A\slashL \varpi^r)$ is truncated, in particular finitely generated over $\pi_0(A)/\varpi^r \cong \pi_0(A\slashL \varpi^r)$.
\end{proof}
\begin{lemma}\label[lemma]{lemma_IntB_to_C/O}
    Let $(A\to B)\in \Int B$. Then the map $A\to E$ factors through $\O\to E$, and we have a colimit-reflecting lifting
\[\begin{tikzcd}
	{\Int B} & {(\CAlgO)_{/O}} \\
	{(\CAlgO)_{/B}} & {(\CAlgO)_{/E}}.
	\arrow[dashed, from=1-1, to=1-2]
	\arrow[from=1-1, to=2-1]
	\arrow[from=1-2, to=2-2]
	\arrow[from=2-1, to=2-2]
\end{tikzcd}\]
\end{lemma}
\begin{proof}
    The functor $\Int B \to (\CAlgO)_{/B}$ factors through
    \[
        (\CAlgO)_{/\colimAIntB A} \simeq (\CAlgO)_{/B\times_E \O} \simeq (\CAlgO)_{/B} \times_{(\CAlgO)_{/E}} (\CAlgO)_{/\O},
    \]
    where we have used \Cref{lemma_IntB}(1) in the first equivalence. Thus, we obtain the sought lifting. Reflection of colimits follows from \cite[Prop.~1.2.13.8]{lurie_htt}, which states that colimits in an over-$\infty$-category $\calC_{/X}$ are computed in $\calC$.
\end{proof}
\subsection{Topologies on animated algebras}
In this section, we discuss the \emph{ind-$\varpi$-adic topology} in the derived setting which we, following \cite[§2.4]{zhu2025}, will use to define derived moduli stacks of continuous representations. Let $\Setcp$ denote the category of finite sets. The starting point is the observation that any static $\O$-module $M$ defines an ind-profinite set, also denoted $M$:
\[
    M = \big(( M'/\varpi^r )_{r\geq 1} \big)_{M'\subset M} \in \Ind\big(\Pro(\Setcp)\big),
\]
where $M'$ runs over the filtered category of finitely generated submodules of $M$. If $M$ is finitely generated, its associated ind-profinite set lies in the essential image of $\Pro(\Setcp)\hookrightarrow \Ind(\Pro(\Setcp))$. The following lemma provides a notion of ind-$\varpi$-adic topology in the derived setting.
\begin{lemma}\label[lemma]{lemma_Csc}\cite[Lem.~2.53]{zhu2025}
    Let $\Gamma$ be a profinite set and define
    \begin{align*}
        \Csc(\Gamma,-) \colon \Mod_\O &\to \Mod_\O \\
        M &\mapsto \Map_{\Ind(\Pro(\Setcp))}(\Gamma, M).
    \end{align*}
    Then $\Csc(\Gamma,-)$ is exact, lax symmetric monoidal, and extends to a $t$-exact functor $\AniMod_\O \to \AniMod_\O$ which preserves colimits. Moreover, $\Csc(\Gamma,-)$ lifts to a 
    functor $\CAlgO \to \CAlgO$ which preserves finite limits and sifted colimits.
\end{lemma}
\begin{lemma}\label[lemma]{lemma_csc_fgA}
    Let $A \in \AniModcn_\O$ and suppose $\pi_\ast(A)$ is a finitely generated $\O$-module. Then
    \[
    \Csc(\Gamma,A)\simeq \lim_{r\geq 1} \Csc(\Gamma,A\slashL\varpi^r).
    \]
    Moreover, if $A \in \CAlgO$, then the equivalence holds in $\CAlgO$.
\end{lemma}
\begin{proof}
    For any $A\in \AniModcn_\O\simeq D^{\geq 0}(\O)$, we have natural morphisms (of $\O$-algebras, if $A\in \CAlgO$)
    \[
        \Csc(\Gamma,A) \to \Csc(\Gamma,\lim A\slashL \varpi^r) \to \lim \Csc(\Gamma,A\slashL \varpi^r).
    \]
    If $A$ is derived $\varpi$-adically complete -- in particular, if $\pi_\ast(A)$ is finitely generated over $\O$ (\Cref{lemma_derived_complete_pi_complete}) -- the first map is a weak equivalence. Let $P(A)$ be the property that the second map is a weak equivalence. Then $P(\O[0])$ holds, $P(A)$ implies $P(A[1])$ and for any fibre sequence $A \to B \to C \to A[1]$, if $P$ holds for two of $A,B,C$ then it holds for the third by the five lemma. Moreover, for any $A$ and $B$, $P(A\oplus B)$ is implied by $P(A)$ and $P(B)$ since $\Csc(\Gamma,-)$ commutes with colimits (\Cref{lemma_Csc}). Finally, any retract $A \to B \to A$ induces a commutative diagram
\[\begin{tikzcd}
	{\Csc(\Gamma,\lim_r A\slashL \varpi^r)} & {\Csc(\Gamma,\lim_r B\slashL \varpi^r)} & {\Csc(\Gamma,\lim_r A\slashL \varpi^r)} \\
	{\lim_r\Csc(\Gamma,A\slashL \varpi^r)} & {\lim_r\Csc(\Gamma,B\slashL \varpi^r)} & {\lim_r\Csc(\Gamma,A\slashL \varpi^r)}
	\arrow[from=1-1, to=1-2]
	\arrow[from=1-1, to=2-1]
	\arrow[from=1-2, to=1-3]
	\arrow[from=1-2, to=2-2]
	\arrow[from=1-3, to=2-3]
	\arrow[from=2-1, to=2-2]
	\arrow[from=2-2, to=2-3]
\end{tikzcd}\]
where the composition of horizontal arrows is the identity,
so that, if $P(B)$ holds, the induced map
\[
  \pi_\ast\big(\Csc(\Gamma, \lim_r A\slashL \varpi^r)\big) \to \pi_\ast\big(\lim_r\Csc(\Gamma, A\slashL \varpi^r)\big)  
\]
is both injective and surjective, i.e. $P(A)$ holds. The lemma now follows from the fact that the category $D_{\operatorname{perf}}(\O)$ of perfect complexes (i.e. complexes with finitely generated homology) is the smallest strictly full, saturated, triangulated subcategory of $D(\O)$ containing $\O[0]$ (\cite[Tag 0ATI]{stacks-project}).  
\end{proof}
\subsection{Moduli stacks of Galois representations}
In this section, we discuss moduli stacks of representations of a profinite group $\Gamma$ valued in a smooth affine group scheme $\bfG/\O$ of finite type, following \cite[§2.4]{zhu2025}. We only consider what is referred to as \emph{strongly continuous representations} in \emph{loc. cit.} 

For $m \in \N$, we let $\O[\bfG^m]$ denote the coordinate ring of $\bfG^m$. With boundary maps defined using the group law on $\bfG$ and degeneracies defined by inclusions $\bfG^{m_1} \times \{\id\} \times \bfG^{m_2} \to \bfG^{m_1+m_2+1}$, we obtain a cosimplicial animated $\O$-algebra $\O[\bfG^\bullet] \in \Fun(\Delta,\CAlgO)$.

Given $A\in \CAlgO$, there is an equivalence (see \cite[Rem.~2.9]{zhu2025}):
\[
    \Map_{\Ani(\Grp)}\big(\Gamma,\bfG(A)\big) \simeq \Map_{\Fun(\Delta,\CAlgO)}\big(\O[\bfG^\bullet], C_{\operatorname{disc}}(\Gamma^\bullet,A)\big),
\]
where $\Ani(\Grp)$ is the $\infty$-category of animated groups, and, for any $[m] \in \Delta$, $C_{\operatorname{disc}}(\Gamma^m,-)$ is the animation of the functor which maps a static $\O$-module $M$ to the set of all functions $\Gamma^m\to M$, cf. the definition of $\Csc(\Gamma^\bullet,-)$ in \Cref{lemma_Csc}. Replacing $C_{\operatorname{disc}}(\Gamma^\bullet,-)$ by $\Csc(\Gamma^\bullet,-)$, we obtain the following definition.
\begin{definition}\label[definition]{def_XGammaMbox}\cite[Def.~2.56]{zhu2025} Let $\Gamma$ be a profinite group and $\bfG/\O$ a smooth affine group scheme of finite type. The \emph{derived moduli stack of continuous framed $\bfG$-valued $\Gamma$-representations} is the functor
    \begin{align*}
        \XGammaGbox\colon \CAlgO &\to \Ani \\
        A &\mapsto  \Map_{\Fun(\Delta,\CAlgO)}\big(\O[\bfG^\bullet], \Csc(\Gamma^\bullet,A)\big).
    \end{align*}
\end{definition}
\begin{lemma}\label[lemma]{lemma_classicalpoints}
    \cite[Lem.~2.57]{zhu2025} If $A\in \CAlgO$ is static and $\bfG = \GL_n$, then $\XGammaGbox(A)$ is the set of homomorphisms $\rho_A\colon\Gamma\to \GL_n(A)$ such that $A^n$ is a filtered colimit of $\Gamma$-stable finitely generated $\O$-modules with continuous $\Gamma$-action. 
\end{lemma}
\begin{lemma}\label[lemma]{lemma_Xfinitelimits}
    Let $\XGammaGbox\colon \CAlgO\to \Ani$ be the functor of \Cref{def_XGammaMbox}. Then $\XGammaGbox$ commutes with finite limits, and for every $j\geq 0$, its restriction to $j$-truncated algebras commutes with filtered colimits. Moreover, if $A\in \CAlgO$ is an algebra such that $\pi_\ast(A)$ is a finitely generated $\O$-module, we have an equivalence
    \[
      \XGammaGbox(\lim_r A\slashL\varpi^r) \simeq \lim_r \XGammaGbox(A\slashL\varpi^r).
    \]
\end{lemma}
\begin{proof}
    The statement for filtered colimits is contained in \cite[Prop.~2.59]{zhu2025}. By definition, $\XGammaGbox$ is given by the composition
    \[
        \Hom_{\Fun(\Delta,\CAlgO)}(\O[\bfG^\bt],-) \circ \Csc(\Gamma,-) \circ \iota,
    \]
    where $\iota\colon \CAlgO\to \Fun(\Delta,\CAlgO)$ is the inclusion of constant cosimplicial objects. Being a right adjoint, $\iota$ preserves all limits. The functor $\Csc(\Gamma^\bullet,-)$ preserves finite limits and limits of the form $\lim A\slashL \varpi^r$ for $A$ as in the statement by \Cref{lemma_csc_fgA}. Finally, $\Hom_{\Fun(\Delta,\CAlgO)}(\O[\bfG^\bt],-)$ preserves all limits. This completes the proof.
\end{proof}
\subsection{Derived deformation functors}\label{ssec_ddf}
In this section, we define various derived deformation functors of Galois representations in terms of the stack $\XGammaGbox$ of \Cref{def_XGammaMbox}. Our main case of interest is when $\bfG = \GL_n$, but we state our results in broad generality.

\begin{definition}
    Let $\calX \colon \CAlgO \to \Ani$ be a functor and fix $\rhobar \in \calX(k)$. We define
    \begin{align*}
        \calX_{\rhobar}\colon\Artk &\to \Ani \\
        (B\to k) &\mapsto \calX(B) \times_{\calX(k)} \{\rhobar\},
    \end{align*}
    which we extend to $\calPro(\Artk)$ by setting $\calX_{\rhobar}(\{A_\alpha\}) \coloneqq \lim_{\alpha} \Xrhobar(A_\alpha)$. Similarly, if $\rho\in \calX(E)$, we obtain a functor $\calX_{\rho}\colon\ArtE\to\Ani$ given by $\calX(B\to E) = \calX(B) \times_{\calX(E)} \{\rho\}.$ 
\end{definition}
\begin{definition}\label[definition]{def_calDboxrhobar}
    Let $\XGammaGbox$ be the derived moduli stack of \Cref{def_XGammaMbox}, and suppose $\rhobar \in \XGammaGbox(k)$. The \emph{unrestricted framed derived deformation functor} of $\rhobar$ is the functor
    \[
        \Dboxrhobar \coloneqq (\XGammaGbox)_{\rhobar} \colon \Artk \to \Ani.
    \]
    Similarly, if $\rho\in \XGammaGbox(E)$, unrestricted framed derived deformation functor of $\rho$ is the functor
    \[
        \Dboxrho \coloneqq (\XGammaGbox)_{\rho} \colon \ArtE \to \Ani.
    \]
\end{definition}
\begin{remark}\cite[Eq.~(2.6)]{zhu2025}
    The derived deformation functors introduced above for a mod~$p$ representation coincide with the ones introduced by Galatius--Venkatesh in \cite[Def.~5.4]{galatius_venkatesh2018}.
\end{remark}
The functors in \Cref{def_calDboxrhobar} parametrise framed deformations. The unframed deformation functors are obtained by passing to a quotient with respect to conjugation. Let $\bfZ \subset \bfG$ denote the centre and define $\bfPG \coloneqq \bfG/\bfZ \colon \CAlgO\to\Ani$. We define functors $\PGhatk,\Ghatk \colon \Artk\to \Ani$ by
\[
    \PGhatk(B\to k) \coloneqq \bfPG(B) \times_{\bfPG(k)} \{\id\}, \quad \Ghatk(B\to k) \coloneqq \bfG(B) \times_{\bfG(k)} \{\id\}.
\]
Similarly, we define $\GhatE,\PGhatE\colon\ArtE\to\Ani$. Defined as above, the groups $\Ghatk$ and $\PGhatk$ act on $\Dboxrhobar$ by conjugation, i.e. there is a simplicial diagram
\[
    (\PGhatk^\bt \times \Dboxrhobar) = \left(\dots \ \PGhatk^2 \times \Dboxrhobar \mathrel{\substack{\textstyle\rightarrow\\[-0.6ex]
    \textstyle\rightarrow \\[-0.6ex]
    \textstyle\rightarrow}} \PGhatk \times \Dboxrhobar \rightrightarrows \Dboxrhobar\right) \in \Fun\big(\Delta^{\operatorname{op}},\Fun(\Artk,\Ani)\big).
\]
By taking a colimit, we obtain the functor of unframed deformations of $\rhobar$.
\begin{definition}\label[definition]{def_unframed}
    With notation as above, the \emph{(unframed) derived deformation functor of $\rho$} is defined as the geometric realisation
\[
    |\PGhatk^\bt \times \Dboxrhobar| \coloneqq \underset{\Delta^{\operatorname{op}}}{\colim}\ \big(\PGhatk^\bt \times \Dboxrhobar\big) \in \Fun(\Artk,\Ani).
\]
The functor $\Drho$ of derived deformations of a characteristic 0 representation $\rho$ is defined analogously.
\end{definition}
\begin{remark}
    \Cref{def_unframed} is modelled after classical Galois deformation theory. An alternative way of defining unframed deformation functors is passing to the quotient at the level of stacks by forming $\XGammaG \coloneqq |\XGammaGbox \times \bfPG^\bullet|$ as in \cite[Def.~2.20]{zhu2025}, then  defining $\calD_{\rhobar} \coloneqq (\XGammaG)_{\rhobar}$ for a point $\rhobar\in\XGammaG(k)$. This yields the same functor up to equivalence, which can be proved by considering the natural map $(\XGammaGbox)_{\rhobar} \to (\XGammaG)_{\rhobar}$ and comparing tangent complexes (see \Cref{section_tangentcomplexes}).
\end{remark}
\begin{definition}\label[definition]{def_deformationproblem}
    A (mod~$p$) \emph{derived deformation problem} is a tuple
    \[
        \bfS = (\rhobar,S,\{\Dv\}_{v\in S})
    \]
    where $S$ is a finite set of places of $F$, $\rhobar \in \XGammaFSGbox(k)$, and $\Dv\in\Fun(\Artk,\Ani)$ are equipped with actions by $\PGhatk$ and equivariant maps $\Dv\to \Dboxrhobarv$. Similarly, a \emph{characteristic 0 derived deformation problem} is a tuple $(\rho,S,\{\Dv\}_{v\in S})$ where $\rho\in \XGammaFSGbox(E)$ and $\Dv\to \Dboxrhov$.
\end{definition}
The following functor appears in \cite[Def.~3.22]{cai2021}.
\begin{definition}\label[definition]{def_localconditions}
    Let $\bfS = (\rhobar,S,\{\Dv\}_{v\in S})$ be a deformation problem. We define a functor
    \[
        \DboxS \coloneqq \Dboxrhobar \times_{\prod_{v\in S} \Dboxrhobarv} \prod_{v\in S} \Dv \colon \Artk \to \Ani
    \]
    which we call the \emph{functor of framed derived deformations of type $\bfS$}. Since the morphisms $\Dv\to \Dboxrhobarv$ are equivariant under conjugation, we may form the geometric realisation
    \[
        \DS \coloneqq |\PGhatk^\bullet\times\DboxS| \in \Fun(\Artk, \Ani).
    \]
    which we call the \emph{functor of (unframed) derived deformations of type $\bfS$}. We make the analogous definitions for a characteristic 0 derived deformation problem.
\end{definition}
If $\calD_v = \Dboxrhobarv$, we say $\bfS$ has trivial (or no) local condition at $v$. This terminology is justified -- at least when $\bfG = \GL_n$ -- since the local derived deformation functor $\Dboxrhobarv$ is represented by the classical framed deformation ring $\Rvbox$ by \cite[Thm.~1.1]{BIP2023} (for $v\in S_p$) and \cite[Thm.~2.5]{shotton2018} (for $v\in S\setminus S_p$) together with \cite[Lem.~7.5]{galatius_venkatesh2018}.
\begin{example}
    Although \Cref{def_localconditions} is quite general, in our application the local conditions will arise from classical deformation rings in the following familiar way. Let $\bfG=\GL_n/F$ and $\rhobar_v\colon \GFv \to \GL_n(k)$.  Any conjugation-invariant quotient $\Rvbox \surj R_{\rhobarv}^{D}$ defines a local condition $\Dv$ via
    \[
        \Dv \simeq \Map_{/k}(R_{\rhobarv}^{D},-) \to \Map_{/k}(\Rvbox,-) \simeq \Dboxrhobarv.
    \]
\end{example}
\begin{definition}
    Let $\bfS = (\rhobar,S,\{\Dv\}_{v\in S})$ be a deformation problem and suppose $T\subseteq S$ is a subset. We define the \emph{functor of derived $T$-framed liftings of type $\bfS$} as
    \[
        \calD^{\square,T}_{\bfS} = \DboxS \times \prod_{v\in T} \Ghatk \colon \Artk \to \Ani.
    \]
    There exists an action of $\Ghatk$ on $\calD^{\square,T}_{\bfS}$ given by
    \[
        \gamma \colon (\rho_A,\{\alpha_v\}_{v\in T}) \mapsto (\gamma \rho_A \gamma^{-1}, \{\gamma\alpha_v\}_{v\in T}).
    \]
    We define the \emph{functor of $T$-framed derived deformations of type $\bfS$} as the geometric realisation
    \[
        \calDTS \coloneqq |\calD^{\square,T}_{\bfS} \times \Ghatk^\bullet| \in \mathscr{Fun}(\Artk,\Ani).
    \]
\end{definition}
\begin{proposition}\label[proposition]{prop_Tframeddimensions}
    Let $\rhobar\colon \GFS\to\GL_n(k)$ and suppose $\bfS = (\rhobar,S,\{\Dv\}_{v\in S})$ is a deformation problem. Then, for any nonempty subset $T\subset S$, the natural morphism $\DS^{T}\to \DS$ is formally smooth of relative dimension $n^2|T|-1$.
\end{proposition}
\begin{proof}
    See e.g. \cite[Lem.~5.1.5]{a'campo2024}.
\end{proof}
We end this section by citing the facts we shall need about the local conditions featuring in the proof of our main theorem. Thus, for the moment we restrict to the case $\bfG = \GL_n$. At $v\in S\setminus S_p$, our local conditions will be trivial.
\begin{proposition}\label[proposition]{prop_dimension_vnmidp} \cite[Prop.~1.2.2]{allen2016} Let $v\nmid p$. Then
    $\dim R_{\rhobarv}^\square[1/\varpi] = n^2,$ and if $\rho_v\colon \GFv \to E$ is a representation such that the Weil--Deligne representation $\WD(\rho_v)$ is generic, the corresponding point of $\Spec R_{\rhobarv}^{\square}[1/\varpi]$ is smooth.
\end{proposition}
At the places dividing $p$, we will consider local conditions defined by a variant of Kisin's \cite{kisin2008} potentially semistable deformation rings introduced in \cite[§5.1]{hevesi2024}, defined by imposing conditions on Hodge--Tate weights and the action of inertia.
\begin{definition}
    Let $v\mid p$. A \textit{Weil--Deligne inertial type} of $F_v$ is an isomorphism class of pairs $\tau = (\rho_\tau, N_\tau)$ where $\rho_\tau\colon I_{F_v} \to \GL_n(E)$ is a representation with open kernel and $N_\tau \in M_n(E)$ is a nilpotent matrix such that $\tau$ is isomorphic to the restriction $(r,N)|_{I_{F_v}}$ of a Weil--Deligne representation of $\GFv$. 
\end{definition}
Hevesi defines a partial order $\preceq$ on Weil--Deligne inertial types, where $(\rho_\tau,N_\tau)\preceq (\rho_{\tau'},N_{\tau'})$ if $\rho_\tau \cong \rho_{\tau'}$ and, roughly, the nilpotence of $N_\tau$ is bounded from below by that of $N_{\tau'}$; see \cite[Def.~2.11]{hevesi2024} for the precise definition.
\begin{theorem}\label{theorem_semistabledeformationrings}
    Let $v\mid p$. Fix a dominant weight $\lambda \in (\Z^n)^{\Hom(F,E)}$ and a Weil--Deligne inertial type $\tau$ of $F_v$. Suppose $\rhobar_v\colon \GFv \to \GL_n(k)$. Then we have the following:
    \begin{enumerate}
        \item[(1)] \cite[Thm.~5.2]{hevesi2024} There exists a local condition $\calDvlambdatau \to \Dboxrhobarv$ represented by an $\O$-flat reduced quotient $\Rvbox \surj \Rvlambdatau$ such that a morphism $\Rvbox \to E$ factors through $\Rvlambdatau$ if and only if the corresponding representation is potentially semistable with labelled Hodge--Tate weights $(\lambda_{\iota,1}+n-1 > \dots > \lambda_{\iota,n})_{\iota \in \Hom(F,E)}$ and $\WD(\rho_v)^{\operatorname{ss}}|_{I_{F_v}} \preceq \tau_v$.
        \item[(2)] \cite[Thm.~3.3.8]{kisin2008} $
        \dim \Rvlambdatau[1/\varpi] \leq n^2 + \frac{1}{2}n(n-1)[F_v\colon \Q_p].  
        $
        \item[(3)] \cite[Prop.~1.3.12]{allen2016} If $\rho_v$ is a characteristic 0 lift of $\rhobar_v$ such that the Weil--Deligne representation $\WD(\rho_v)$ is generic, the corresponding point of $\Spec(\Rvlambdatau[1/\varpi])$ is smooth with tangent space equal to the local Bloch--Kato Selmer group $H^1_g(\GFv, \ad\rho)$. Localising and completing at this point yields $R^{\lambda_v}_{\rho_v}$, the deformation ring parametrising potentially semistable deformations of Hodge type $\lambda_v$.  
    \end{enumerate}
\end{theorem}
\subsection{Tangent complexes and representability}\label[section]{section_tangentcomplexes}
In this section, we recall the definition of the tangent complex and Lurie's representability theorem for formally cohesive functors, \cite[§6.2]{lurie_dag}.

\begin{definition}\label[definition]{def_fcf}\cite[Def.~6.2.1]{lurie_dag}
    A functor $\calD \colon \Artk \to \Ani$ (resp. $\calD \colon \ArtE \to \Ani$) is \emph{formally cohesive} if the following hold:
    \begin{itemize}
        \item[(1)] $\calD(k)$ is contractible (resp. $\calD(E)$ is contractible).
        \item[(2)] If $B' \to B$ and $B''\to B$ are surjective (on $\pi_0$) morphisms in $\Artk$ (resp. $\ArtE$), then \[\calD(B'\times_B B'')\simeq \calD(B')\times_{\calD(B)} \calD(B'').\] 
    \end{itemize}
\end{definition}
\begin{lemma}
    The derived deformation functors introduced in Section 3.5 are formally cohesive.
\end{lemma}
\begin{proof}
    \Cref{def_fcf}(1) follows from \Cref{lemma_classicalpoints}. To prove property (2), we let $\calX^\square = \XGammaGbox$ and define $\Drhobar^\square \coloneqq \Xboxrhobar$ for some $\rhobar\in \calX^\square(k)$. For any diagram $B' \to B \leftarrow B''$ in $\Artk$, we have
    \[
        \calD(B'\times_B B'') = \calXbox(B'\times_B B'') \times_{\calXbox(k)} \{\rhobar\} \simeq \Xboxrhobar(B') \times_{\Xboxrhobar(B)} \Xboxrhobar(B'') \simeq \calD(B') \times_{\calD(B)} \calD(B''),
    \]
    as needed. This equivalence is equivariant for the action of \[\PGhat(B'\times_B B'') = \PGhat(B')\times_{\PGhat(B)}\PGhat(B''),\] which implies the claim for the unframed version $|\PGhat^\bullet \times \calD^\square|$. The proof for a characteristic 0 deformation functor is identical.
\end{proof}
Given a formally cohesive functor $\calD\colon \Artk\to \Ani$, $r\geq 1$ and $\xi_r \in \calD(\O_r)$ where $\O_r = \O/\varpi^r$, one can define a functor on perfect connective $\O_r$-modules by
\[
    \opT \colon V \mapsto \calD(\O_r \oplus V).
\]
Here, perfect means $V$ has finitely generated homotopy. By \cite[Prop.~2.5.5, Lem.~6.2.2]{lurie_dag}, this functor has a unique filtered colimit-preserving and finite limit-preserving extension to the category $\AniMod_{\O_r}$ of all $\O_r$-modules, which is representable in the sense that there exists an $\O_r$-module $\opT^\bullet_{\O_r}(\calD) \in \AniMod_{\O_r}$ such that 
\[
    \opT(V) = \big(V \otimesL_{\O_r} \opT^\bullet_{\O_r}(\calD)\big)_0,
\]
the 0:th space of the spectrum $V\otimesL_{\O_r} \opT^\bullet_{\O_r}(\calD)$. The homotopy groups of $\opT^\bullet_{\O_r}(\calD)$ are given by
\[
    \opT^i_{\O_r}(\calD) \coloneqq \pi_i\big(\opT^\bullet_{\O_r}(\calD)\big) \cong \pi_0\big(\calD(\O_r \oplus \O_r[i]) \times_{\calD(\O_r)} \{\xi_r\}\big).
\]
The $\O_r$-module $\opT^\bullet_{\O_r}(\calD) \in \AniMod_{\O_r}$ is called the \emph{tangent complex of $\calD$ at $\xi_r$} (or at $\O_r$, if $\xi_r$ is inferrable from context). Similarly, for a formally cohesive functor $\calD\colon \ArtE \to \Ani$, we define its tangent complex $\opT^\bullet_E(\calD)$ at $E$ in terms of the functor mapping a perfect connective $V\in \AniMod_E$ to $\calD(E\oplus V)$. The construction of the tangent complex is functorial in $\calD$ and $\O_r$, and we define
\[
    \opT^\bullet_\O(\calD) \coloneqq \lim_r \opT^\bullet_{\O_r}(\calD).
\]
If $\calD$ is a formally cohesive functor represented by $\calR$, we write $\opT^\bullet_{?}(\calR) \coloneqq \opT^\bullet_{?}(\calD).$ The tangent complex of $\calR$ at $\xi_r\colon \calR \to \O_r$ is given by
    \[
        \opT^\bullet_{\O_r}(\calR) = \Map_{\AniMod_\calR}(\opL_{\calR/\O},\O_r) \cong \Map_{\AniMod_{\O_r}}(\opL_{\calR/\O} \otimesL_{\calR} \O_r, \O_r),
    \]
where $\opL_{\calR/\O}$ denotes the algebraic cotangent complex of $\calR$. If $A\to B \to \O$ are morphisms in $\CAlgO$, we denote by $\opT^\bullet_\O(B/A)$ the fibre of the induced map $\opT^\bullet_\O(B) \to \opT^\bullet_\O(A)$.

For any $\calR \in \CNLk$, the natural morphism $\calR \to \pi_0(\calR)$ induces an isomorphism
\[
    \opT^0_k(\pi_0(\calR)) \cong \opT^0_k(\calR)
\]
and an injection $\opT^1_k(\pi_0(\calR)) \cong \opT^1_k(\calR)$, and similarly with $\O_r$-coefficients. If $R=\pi_0(R) \in \CNLk$ is formally smooth, then $\opT^i_k(R)=0$ for $i\geq 1$. Similar statements hold for $\calR$ and $R$ in $\CNLE$.

\begin{theorem}\label{theorem_lurierepresentability}\cite[Cor.~6.2.14]{lurie_dag}
    Let $\calD\colon\Artk \to \Ani$ be a formally cohesive functor. The following are equivalent:
    \begin{itemize}
        \item[(1)] There exists an $\calR \in \CNLk$ (\Cref{definition_CNLk_CNLE}) and an equivalence of functors \[\calD \simeq \Map_{/k}(\calR,-).\]
        \item[(2)] The $k$-vector spaces $\opT^i_k(\calD)$ vanish for $i<0$ and are finite dimensional for $i\geq 0$.
    \end{itemize}
    The analogous statement holds for a formally cohesive functor $\calD\colon\ArtE\to \Ani$.
\end{theorem}
Given $\calR\in \CNLk$ with maximal ideal $\frakm\subset \pi_0(\calR)$, we have equivalences of functors
\[
    \Map_{/k}(\calR,-) \simeq \Map_{/k}(\lim_r \calR\slashL \frakm^r,-) \simeq \underset{r}{\colim}\Map_{/k}(\calR\slashL \frakm^r,-).
\]
The first equivalence holds because $\calR$ is derived $\frakm$-adically complete and the second holds by \Cref{prop_lurie_prohomspaces}. The quotient $\calR\slashL\frakm^r$ need not be truncated, but one can rephrase \Cref{theorem_lurierepresentability}(1) in terms of pro-representability using the pro-artinian completion $\{\tau_{\leq j}(\calR\slashL \frakm^r) \}_{j,r}\in \calPro(\Artk)$, as in \cite{galatius_venkatesh2018} (see also \cite[Rmk.~6.2.11]{lurie_dag}).
\begin{proposition}\label[proposition]{prop_Tifinitelygenerated}
    Let $\calR\in \CNLk$ and let $\calR \to \O_r = \O/\varpi^r$ be a morphism over $k$. Then, for every $i\geq 0$,
    $
        \opT^i_{\O_r}(\calR)
    $
    is a finitely generated $\O_r$-module.
\end{proposition}
\begin{proof}
    The case $r=1$ is \Cref{theorem_lurierepresentability}. By considering the short exact sequence
    \[
        0 \to \O_r \to \O_{r+1} \to k \to 0
    \]
    we deduce the analogous statement for $\O_r$-coefficients using induction.
\end{proof}
\begin{lemma}\label[lemma]{lemma_Otangentkomplex}
    Suppose $A \to B \to \O$ in $\CAlgO$ and that $B$ is almost of finite presentation over $A$. Then, if $\opT^i_\O(A)$ is a finitely generated $\O$-module for every $i\geq 0$, the same holds for $\opT^i_\O(B)$.
\end{lemma}
\begin{proof}
    By \cite[Prop.~3.2.14]{lurie_dag}, the relative cotangent complex $\opL_{B/A}$ is an almost perfect connective $B$-module, i.e. $\pi_i(\opL_{B/A})$ is a finitely generated $\pi_0(B)$-module for every $i\geq 0$. It follows that $\opL_{B/A}\otimesL_{B} \O$ is an almost perfect connective $\O$-module. Dualising the fibre sequence
    \[
        \opL_{A/\O} \otimesL_{A}\ \O \to \opL_{B/\O} \otimesL_{B}\ \O \to \opL_{B/A} \otimesL_{B}\ \O \stackrel{+1}{\to}
    \]
    we obtain a fibre sequence of tangent complexes with $\O$-coefficients
    \[
        \opT^\bullet_\O({B/A})  \to \opT^\bullet_\O(B) \to \opT^\bullet_\O(A) \stackrel{+1}{\to}.
    \]
    By assumption, $\opT^i_\O(A)$ is finitely generated over $\O$ for every $i$. The same now follows for $\opT^\bullet_\O(B)$.
\end{proof}
The tangent complexes of the deformation functors introduced in \Cref{ssec_ddf} are given by continuous cohomology in the usual sense.
\begin{proposition}\label[proposition]{prop_tangentcomplexesgaloiscohomology}\cite[Prop.~2.59]{zhu2025}\cite[Lem.~3.10,~3.12]{cai2021}
    Let $\XGammaGbox$ be the derived moduli stack of framed $\bfG$-valued $\Gamma$-representations, Fix $\rhobar \in \XGammaGbox(k)$ and let $\Dboxrhobar = (\XGammaGbox)_{\rhobar}$. Then the tangent complex of $\Dboxrhobar$ at $k$ is given by an object $\overline{C}_{\operatorname{sc}}^\bullet(\Gamma, \ad\rhobar)[1] \in \AniMod_k$ such that
    \[
        \opT^i_k(\Dboxrhobar) \cong \pi_{-i}\big(\overline{C}_{\operatorname{sc}}^\bullet(\Gamma, \ad\rhobar)[1]\big) \cong \begin{cases} Z^1(\Gamma, \ad\rhobar)&\text{ if }i=0,\\
            H^{2}(\Gamma, \ad\rhobar)&\text{ if }i=1,\\
            0&\text{ otherwise.}
        \end{cases} 
    \]
    The tangent complex of $\Drhobar$ at $k$ is given by an object $\Csc^\bullet(\Gamma, \ad\rhobar)[1] \in \AniMod_k$ such that
    \[
        \opT^i_k(\Drhobar) \cong \pi_{-i}\big(\Csc^\bullet(\Gamma, \ad\rhobar)[1]\big) \cong \begin{cases} 
            H^0(\Gamma, \ad^0\rhobar)&\text{ if }i=-1,\\
            H^1(\Gamma, \ad\rhobar)&\text{ if }i=0,\\
            H^{2}(\Gamma, \ad\rhobar)&\text{ if }i=1,\\
            0&\text{ otherwise.}
        \end{cases} 
    \]
    In particular, if $\rhobar$ is Schur, then
    $
        \opT^i_k(\Drhobar) \cong H^{i+1}(\Gamma, \ad\rhobar).
    $
    The analogous statements are true with $\xi \in \XGammaGbox(\O)$, $\xi_r \in \XGammaGbox(\O_r)$ or $\rho \in \XGammaGbox(E)$ in place of $\rhobar$.
\end{proposition}
\begin{corollary}
    Let $\Gamma=\Gal(F_S/F)$ or $\Gal(\overline{F}_v/F_v)$ and fix $\rhobar \in \XGammaGbox(k)$. The framed derived deformation functor $\Dboxrhobar = (\XGammaGbox)_{\rhobar}$ is representable and if $\rhobar$ is Schur, the derived deformation functor $\Drhobar$ is representable. The same statement is true for $\rho\in \XGammaGbox(E)$.
\end{corollary}
\begin{proof}
    This follows from \Cref{theorem_lurierepresentability,prop_tangentcomplexesgaloiscohomology}.
\end{proof}
We deduce the same statement for deformations of type $\bfS$.
\begin{proposition}
    Let $\bfS = (\rhobar,S,\{\Dv\}_{v\in S})$ be a deformation problem and let $T\subseteq S$ be a subset. Suppose $\rhobar$ is Schur. Then the functors $\DS$ and $\calDTS$ are representable.
\end{proposition}
\subsection{Localisation at a characteristic 0 point}\label{subsec_localisationatachar0pt}
In this section, we prove two results which generalise results of Kisin, relating the deformation ring of a characteristic 0 representation to that of its reduction mod~$p$ \cite[Lem.~2.3.3, Prop.~2.3.5]{kisin2009}.

In the following definition, which is a direct analogue of \cite[§2.3.1]{kisin2009}, we use the functor $\Int B \to (\CAlgO)_{/\O}$ of \Cref{lemma_IntB_to_C/O} to view $A\in \Int B$ as an object of $(\CAlgO)_{/\O}$.
\begin{definition}\label[definition]{def_Frhobarxi}
    Let $\calD \colon\Artk \to \Ani$ be a functor and $\xi \in \calD(\O) = \lim \calD(\O_r)$. We define
    \begin{align*}
        \calD_{(\xi)}\colon \ArtE &\to \Ani \\
        (B \to E) &\mapsto \colimAIntB \calD(A) \times_{\calD(\O)} \{\xi\}.
    \end{align*} 
\end{definition}
As before, we extend $\calD_{(\xi)}$ to $\calPro(\ArtE)$ by defining $\calD_{(\xi)}(\{A_\alpha\}) \coloneqq \lim_\alpha \calD_{(\xi)}(A_\alpha)$.
\begin{proposition}\label[proposition]{prop_FrhobarxiFrho}
    Let $\calXbox \colon \CAlgO \to \Ani$ be a functor which commutes with finite limits and limits of the form $A = \lim A\slashL \varpi^r$ when $\pi_\ast(A)$ is a finitely generated $\O$-module. Moreover, suppose the restriction of $\calXbox$ to $j$-truncated algebras commutes with filtered colimits. Fix $\rhobar\in\calXbox(k)$ and $\xi\in \calXbox_{\rhobar}(\O)$, and let $\rho$ be the image of $\xi$ in $\calXbox(E)$. Then there is an equivalence of functors
    \[
        (\calXbox_{\rhobar})_{(\xi)} \simeq \calXbox_\rho.
    \]
\end{proposition}
\begin{proof}
    Let $(B\to E)\in \ArtE$. By definition,
    \[
        \Xboxrhobarxi(B) = \underset{A\in \Int B}{\colim} (\Xboxrhobar(A) \times_{\Xrhobar(\O)}\{\xi\}).
    \]
    Now, $A \simeq \lim A\slashL \varpi^r$ and $A\slashL \varpi^r \in \Artk$ by \Cref{lemma_IntB}(3,4) Hence, if we let $\xi_r \in \calXbox(\O_r)$ denote the $\mod \varpi^r$ reduction of $\xi$, we have by definition that
    \[
        \Xboxrhobar(A) \times_{\Xboxrhobar(\O)} \{\xi\} \simeq \lim_r \Xboxrhobar(A\slashL \varpi^r)\times_{\Xboxrhobar(\O_r)} \{\xi_r\} \simeq \calXbox(A) \times_{\calXbox(\O)} \{\xi\}.
    \]
    The category $\Int B$ is filtered by \Cref{prop_IntBfiltered} and filtered colimits commute with finite limits. Hence, by the assumption on $\calXbox$ and \Cref{lemma_IntB}(1),
    \[
        \underset{A\in \Int B}{\colim} (\Xrhobar(A) \times_{\Xrhobar(\O)}\{\xi\}) \simeq \calXbox(B\times_E \O) \times_{\calXbox(\O)} \{\xi\} \simeq \calXbox(B) \times_{\calXbox(E)} \{\rho\} = \calXbox_\rho(B),
    \]
     since $A \in \Int B$ has non-zero homotopy in the same degrees as $B$, which is truncated, does.
\end{proof}
\begin{theorem}\label{theorem_charptochar0rings}
    Let $\calD \colon\Artk \to \Ani$ be a formally cohesive functor represented by $\calR\in \CNLk$, and $\xi \in \calD(\O).$ Then $\Dxi\colon \ArtE\to \Ani$ is represented by
    \[
        (\mathscr{R}\otimes_\O E)_\xi^\wedge \coloneqq \lim_s (\mathscr{R}\otimes_\O E)\slashL I_\xi^s,
    \]
    the derived completion of $\mathscr{R}\otimes_\O E$ at the maximal ideal $I_\xi = \ker(\pi_0(\calR)\otimes_\O E \to E)$ defined by $\xi$.
\end{theorem}
\begin{proof} Let $(B\to E)\in \ArtE$ and let $\xi_r\in \calD(\O_r)$ denote the $\mod \varpi^r$ reduction of $\xi$. By definition and \Cref{lemma_IntB}(3),
    \[
    \Dxi(B) = \underset{A\in \Int B}{\colim}(\calD(A)\times_{\calD(\O)} \{\xi\}) \simeq \underset{A\in \Int B}{\colim}\big(\lim_r \calD(A\slashL \varpi^r)\times_{\calD(\O\slashL \varpi^r)} \{\xi_r\}\big).
    \]
    Now, by representability of $\calD$, we have
    \[
        \lim_r \big(\calD(A\slashL \varpi^r)\times_{\calD(\O\slashL \varpi^r)} \{\xi_r\}\big) \simeq \lim_r \Map_{/k}(\calR,A\slashL \varpi^r) \times_{\Map_{/k}(\calR,\O\slashL \varpi^r)} \{\xi_r\},
    \]
    and from \Cref{lemma_IntB}(4) and commutativity of limits, we obtain
    \[
        \lim_r \Map_{/k}(\calR,A\slashL \varpi^r) \times_{\Map_{/k}(\calR,\O\slashL \varpi^r)} \{\xi_r\} \simeq \Map_{/\O}(\calR,A).
    \]
    For any $A\in \Int B$, the map $\pi_0(A) \to \pi_0(B)$ is injective and hence, for any $\calR \to A$ in $(\CAlgO)_{/\O}$, the induced map $\pi_0(\calR) \to \pi_0(A)$ factors through a quotient $\pi_0(\calR)/(I_\xi^\circ)^s$ for some $s\geq 1$, where $I_\xi^\circ \subset \pi_0(\calR)$ is the inverse image of $I_{\xi}=\ker(\pi_0(\calR) \otimes_\O E \to E)$. \Cref{prop_lurie_prohomspaces} implies
    \[
        \Map_{/\O}(\calR,A) \simeq \underset{s}{\colim} \Map_{/\O}\big(\calR\slashL (I_\xi^\circ)^s,A\big).
    \]
    For any $A\in \Int B$, $\pi_\ast(A)$ injects into $\pi_\ast(B)$, which is truncated. Hence, $A$ is $j$-truncated for $j$ sufficiently large depending only on $B$. Now, $\calR\slashL (I_\xi^\circ)^s$ is almost of finite presentation over $\O$, since $\pi_0(\calR\slashL (I_\xi^\circ)^s)\cong \pi_0(\calR)/(I_\xi^\circ)^s$ is a finitely generated $\O$-module and $\calR\slashL (I_\xi^\circ)^s$ is a derived tensor product of noetherian animated algebras (\Cref{prop_almostfinitelypresented}). Moreover, $\Int B$ is filtered (\Cref{prop_IntBfiltered}), and hence
    \[
        \Dxi(B) \simeq \underset{A\in \Int B}{\colim} \underset{s}{\colim} \Map_{/\O}(\calR\slashL (I_\xi^\circ)^s, A) \simeq \underset{s}{\colim} \Map_{/\O}(\calR\slashL (I_\xi^\circ)^s, \underset{A\in \Int B}{\colim} A).
    \]
    By \Cref{lemma_IntB}(1), $\colim A \simeq B\times_E \O$ and the space above is equivalent to
    \[
        \underset{s}{\colim}\Map_{/\O}\big((\calR\slashL (I_\xi^\circ)^s, B\times_E \O\big) 
        \simeq 
        \underset{s}{\colim}\Map_{/E}\big((\calRS\otimes_\O E)\slashL I_\xi^s, B\big).
    \]
    Thus, $\Dxi(B\to E) \simeq \Map_{/E}\big((\calR \otimes_\O E)_\xi^\wedge, B\big)$ as required.
\end{proof}
The following is a typical application of \Cref{prop_FrhobarxiFrho,theorem_charptochar0rings}.
\begin{corollary}
    Let $\XGammaGbox$ be the derived moduli stack of framed $\bfG$-valued $\Gamma$-representations, and let
    $
        \Dboxrhobar = (\XGammaGbox)_{\rhobar}
    $
    for some $\rhobar\in \XGammaGbox(k)$. Then, for any $\xi \in \Dboxrhobar(\O)$ with $\rho = \xi[1/p] \in \XGammaGbox(E)$,
    \[
        (\Dboxrhobar)_{(\xi)} \simeq \Dboxrho.
    \]
    The same statement holds for the functors $\Drhobar$ and  $\Drho$ of unframed deformations. In particular, if $\Drhobar$ and $\Drho$ are represented by $\calR_{\rhobar}$ and $\calRrho$, respectively. Then
    \[
        \calRrho \simeq \lim_s\ (\calR_{\rhobar} \otimes_\O E)_\xi^\wedge,
    \]
    the derived completion of $\calR_{\rhobar}$ at the ideal $I_\xi=\ker(\pi_0(\calR_{\rhobar}\otimes_\O E)\to E)$ defined by $\xi\colon \calR_{\rhobar}\to \O$.
\end{corollary}
\begin{proof}
    The statement for the framed functors follows from \Cref{lemma_Xfinitelimits,prop_FrhobarxiFrho}. To prove the analogous statement for the unframed functor, we let $(B\to E)\in \ArtE$. The proof of \Cref{prop_FrhobarxiFrho} comprises an equivalence
    \[
        \Dboxrhobarxi(B) = \underset{A\in \Int B}{\colim} \XGammaGbox(A) \times_{\XGammaGbox(\O)} \{\xi\} \simeq \XGammaGbox(B) \times_{\XGammaGbox(E)} \{\rho\},
    \]
    which is equivariant for the action of
    \[
        \underset{A\in \Int B}{\colim} \PGhat_k(A) \times_{\PGhat_k(\O)} \{\id\} \simeq \underset{A\in \Int B}{\colim} \bfPG(A) \times_{\bfPG(\O)} \{\id\} \simeq \bfPG(B) \times_{\bfPG(E)} \{\id\} = \PGhat_E(B),
    \]
    Hence, taking geometric realisations proves the equivalence in the unframed case. The final statement is obtained from \Cref{theorem_charptochar0rings}.
\end{proof}
\begin{proposition}\label[proposition]{prop_limcolim}
    Let $\calD \colon\Artk \to \Ani$ be a formally cohesive functor such that $\opT^i_k(\calD)$ is finite for every $i$, and fix $\xi \in \calD(\O)$. Then $\calD_{(\xi)}\colon\ArtE \to \Ani$ has tangent complex
    \[
        \opT^\bullet_E(\Dxi) \simeq \big(\lim_r \opT_{\O_r}^\bullet(\calD)\big)\otimes_\O E.
    \]
\end{proposition}
\begin{proof}
    It suffices to prove that, for every perfect connective $V\in \AniMod_\O$ (i.e. $\pi_\ast(V)$ is finitely generated over $\O$), we have an equivalence
    \[
        \calD(\O \oplus V) \times_{\calD(\O)} \{\xi\} \simeq \tau_{\geq 0}\big(\lim_r \opT^\bullet_{\O_r}(\calD) \otimesL_\O V\big) \simeq \lim_r \tau_{\geq 0}\big(\opT^\bullet_{\O_r}(\calD) \otimesL_\O V\big).
    \]
    Arguing as in \Cref{prop_Tifinitelygenerated}, one sees that $\opT^i_{\O_r}(\calD)$ is finitely generated for every $i$. Thus, if $V_r \coloneqq V\slashL \varpi^r$, and $\xi_r$ is the image of $\xi$ in $\calD(\O_r)$, then
    \begin{equation*}
        \begin{split}
            \tau_{\geq 0} \big(\lim_r \opT^\bullet_{\O_r}(\calD) \otimesL_\O V\big) 
            \simeq \lim_r \tau_{\geq 0}\big( \opT^\bullet_{\O_r}(\calD) \otimesL_{\O_r} V_r\big)
            &\simeq \lim_r \calD(\O_r \oplus V_r) \times_{\calD(\O_r)} \{\xi_r\}\\ &\simeq \calD(\O \oplus V) \times_{\calD(\O)} \{\xi\},
        \end{split}
    \end{equation*}
    where we have used that $\opT^\bullet_{\O_r}(\calD)\otimesL_{\O_r} V_r$ is almost perfect, so that the inverse limit has no higher derived functors.     
\end{proof}
\begin{proposition}\label[proposition]{prop_Tifinitelengthsmoothpoint}
    Let $R = \pi_0(R) \in \CNLk$ be static and suppose $\xi \colon R \to \O$ is a morphism such that $I_\xi = \ker(R\otimes_\O E \to E)$ is a smooth point. Then, for any $r\geq 1$ and $i\geq 1$, the $\O$-module $\opT^i_{\O_r}(R)$ has finite length bounded independently of $r$.
\end{proposition}
\begin{proof}
Let $i,r\geq 1$. Consider the short exact sequence $
    0 \to \O \to \O \to \O_r \to 0
$ given by multiplication by $\varpi^r$, and the induced exact sequence
\[
    \opT^i_\O(R) \to \opT^i_{\O_r}(R) \to \opT^{i+1}_\O(R).
\]
By \Cref{prop_Tifinitelygenerated,lemma_Otangentkomplex}, these are finitely generated $\O$-modules. Since $\xi$ defines a smooth point,  $\opT^i_\O(R)$ is torsion for $i\geq 1$ by \Cref{prop_limcolim}. It follows that
\[
    \length_\O\big(\opT^i_{\O_r}(R)\big) \leq \length_\O\big(\opT^i_\O(R)\big) + \length_\O\big(\opT^{i+1}_\O(R)\big)
\]
which is finite and independent of $r$. 
\end{proof}
The construction of \Cref{def_Frhobarxi} behaves well with respect to finite limits, as the following simple observation shows.
\begin{lemma}\label[lemma]{lemma_xihatpullbacks}
    Suppose given a pullback of formally cohesive functors
    \[
        \calD' \times_{\calD} \calD'' \in \Fun(\Artk,\Ani),
    \]
    and let $\tilde{\xi} \in (\calD' \times_{\calD} \calD'')(\O)$. Then
    \[
        (\calD' \times_{\calD} \calD'')_{(\tilde{\xi})} \simeq \calD'_{(\xi')} \times_{\calD_{(\xi)}} \calD''_{(\xi'')}
    \]
    where $\xi,\xi',\xi''$ denote the images of $\tilde{\xi}$ in $\calD(\O),\calD'(\O)$ and $\calD''(\O)$, respectively.
\end{lemma}
\begin{proof}
    This follows from the fact that filtered colimits commute with finite limits, and that $\Int B$ is filtered (\Cref{prop_IntBfiltered}). Indeed, for any $B\in \ArtE$,
    \[
        (\calD' \times_{\calD} \calD'')_{(\tilde{\xi})}(B) \simeq \underset{A \in \Int B}{\colim}  
        (\calD' \times_{\calD} \calD'')(A)\simeq \calD'_{(\xi')}(B) \times_{\calD_{(\xi)}(B)} \calD''_{(\xi'')}(B),
    \]
    as needed.
\end{proof}
\begin{proposition}\label[proposition]{prop_S(rho)}
    Let $\DS$ be the derived deformation functor of the deformation problem
    \[
        \bfS = (\rhobar, S, \{\Map_{/k}(R_{\rhobarv}^D,-)\}_{v\in S}).
    \]
    Suppose $\xi\in \DS(\O)$ and let $\rho = \xi[1/p]\colon \GFS \to \GL_n(E)$. Then $\DSxi\simeq\calD_{\bfS(\rho)}$, where $\bfS(\rho)$ is the characteristic 0 deformation problem
    \[
        \bfS(\rho)= (\rho,S,\{\Map_{/E}((R_{\rhobarv}^D\otimes_\O E)_{\xi_v}^\wedge,-)\}_{v\in S})
    \]
    and $(R_{\rhobarv}^D\otimes_\O E)_{\xi_v}^\wedge$ denotes the completion of $R_{\rhobarv}^D\otimes_\O E$ at the ideal $I_{\xi_v} = \ker(\pi_0(R_{\rhobarv}^D\otimes E)\to E)$ defined by $\xi_v = \xi|_{\GFv}$, as in \Cref{theorem_charptochar0rings}.
\end{proposition}
\begin{proof}
    Let $\bfS = (\rhobar,S,\{\calD_v\}_{v\in S})$ and $(B\to E)\in \ArtE$. 
    By definition, the functor of framed deformations of type $\bfS$ is given by the pullback
    $\DSbox = \Dboxrhobar \times_{\prod_{v\in S} \Dboxrhobarv} \prod_{v\in S} \Dv$. By \Cref{lemma_xihatpullbacks},
    \[
        \DSboxxi \simeq \bigg( \Dboxrhobar \times_{\prod_{v\in S} \Dboxrhobarv} \prod_{v\in S} \Dv \bigg)_{(\xi)} \simeq\  (\Dboxrhobar)_{(\xi)} \times_{\prod_{v\in S} (\Dboxrhobarv)_{(\xi_v)}} \prod_{v\in S} (\Dv)_{(\xi_v)}.
    \]
    Applying \Cref{prop_FrhobarxiFrho,theorem_charptochar0rings} yields an equivalence $\DSboxxi \simeq \calD^\square_{\bfS(\rho)}$ which is conjugation-equivariant. Taking geometric realisations completes the proof.
\end{proof}

\begin{corollary}\label[corollary]{cor_RSphat_equals_calRrhoss}
    With notation as in \Cref{theorem_semistabledeformationrings}, let $\DS$ be the deformation functor of
    \[
        \bfS = (\rhobar,S,\{\calDvlambdatau\}_{v \in S_p} \cup \{\Dboxrhobarv\}_{v\in S\setminus S_p}).
    \]
    Suppose $\rhobar\colon \GFS \to \GL_n(k) $ is absolutely irreducible, and let $\xi\in\DS(\O)$ be an $\O$-lattice in a characteristic 0 lift $\rho \colon \GFS \to \GL_n(E)$ such that $\WD(\rho|_{\GFv})$ is generic for every $v\in S$. Then 
    $
        (\DS)_{(\xi)} \simeq \DrhoSss,
    $
    where $\DrhoSss$ is the deformation functor of the characteristic 0 deformation problem
    \[
        (\rho,S,\{\calD_{\rho_v}^{\lambda_v}\}_{v\in S_p} \cup \{\calD_{\rho_v}^\square\}_{v\in S^p}).
    \]
    Moreover, $\opT^i_E(\DrhoSss)=0$ for $i\neq 0,1$, and there are canonical isomorphisms
    \[
        \opT^0_E\big(\DrhoSss\big) \cong H^1_g(\GFS, \ad\rho),\quad \opT^1_E\big(\DrhoSss\big) \cong H^1_g(\GFS, \ad\rho(1))^\vee.
    \]
\end{corollary}
\begin{proof}
    By \Cref{prop_S(rho)}, we have an equivalence $\DSxi \simeq \calD_{\bfS(\rho)}$, where $\bfS(\rho)$ has trivial local conditions at $v\in S\setminus S_p$, and at $v\in S_p$, the local conditions are defined by
    \[
        (\Rvlambdatau\otimes_\O E)_{\xi}^\wedge \simeq \varprojlim_s\ (\Rvlambdatau\otimes_\O E)/I_\xi^s \simeq \Rvlv,
    \]
    where we have used \Cref{theorem_semistabledeformationrings}(3) and the assumption on $\WD(\rho_v)$, and $\varprojlim$ denotes the $1$-categorical inverse limit in the category of rings. Note that the local conditions are defined by smooth algebras  by \Cref{prop_dimension_vnmidp} and \Cref{theorem_semistabledeformationrings}(3). To describe the tangent complex, we use the argument in \cite[Thm.~9.2, App.~B]{galatius_venkatesh2018} which proves the corresponding statement for mod~$p$ deformation problems defined by formally smooth local conditions. The same argument, applied to $\DrhoSss$, identifies $\opT^\bt_E(\Drhoss)$ with the complex constructed in \Cref{cor_localconditions}; the result follows.
\end{proof}
Having determined the tangent complex $\opT^\bullet_E(\calRrhoss)$, we are able to identify the graded ring $\pi_\ast(\calRrhoss)$ under the assumption that the adjoint Bloch--Kato Selmer group vanishes.
\begin{corollary}\label[corollary]{cor_extalgebra}
    With assumptions and notation as in \Cref{cor_RSphat_equals_calRrhoss}, let $\calRrhoss$ represent $\DrhoSss$ and suppose $H^1_g(\GFS, \ad\rho)=0$. Then there is a natural isomorphism
    \[
        \pi_\ast(\calRrhoss) \cong \bigwedge H^1_g(\GFS, \ad\rho(1)).
    \]
\end{corollary}
\begin{proof}
    By \Cref{cor_RSphat_equals_calRrhoss} and the assumed vanishing of $H^1_g(\GFS, \ad\rho)$, the tangent complex $\opT^\bt_E(\calRrhoss)$ is concentrated in degree $1$. It follows from \cite[Cor.~2.1.6]{arinkin_gaitsgory2015} that $\pi_\ast(\calRrhoss)$ is an exterior algebra on $\pi_1(\calRrhoss)$ which is of dimension $l_0 = \dim_E \opT^1_E(\calRrhoss)$ over $E$ (\Cref{prop_h1gdimension}). Thus, it suffices to define a natural surjection
    \[
        \pi_1(\calRrhoss)\surj H^1_g(\GFS,\ad\rho(1)).
    \]
    By \Cref{cor_RSphat_equals_calRrhoss}, we have an isomorphism
    \[
        H^1_g(\GFS, \ad\rho(1))^\vee \cong \opT^1_E(\calRrhoss) = \pi_0\big(\Map_{/E}(\calRrhoss,E\oplus E[1])\big).
    \]
    Any class $[f\colon \calRrhoss \to E\oplus E[1]]$ induces a morphism $\pi_1(\calRrhoss) \to E$. This defines a natural map
    \[
        \varphi\colon \opT^1_E(\calRrhoss) \to \pi_1(\calRrhoss)^\vee.
    \]
    We claim that $\varphi$ is injective. Indeed, suppose $\varphi([f])=0$, i.e. $f\colon\calRrhoss \to E\oplus E[1]$ induces the trivial map on $\pi_1$. \cite[Prop.~3.3.3]{lurie_dag} implies that the space of such maps is equivalent to the subspace
    \[
        \Map_{/E}\big(\pi_0(\calRrhoss),E\oplus E[1]\big) \hookrightarrow \Map_{/E}\big(\calRrhoss,E\oplus E[1]\big).
    \]
    By assumption, $\pi_0(\calRrhoss)\cong E$ and thus, this space is contractible. It follows that $f$ factors through the augmentation $\calRrhoss\to E$, and that $\varphi$ is injective. Passing to duals, we obtain a canonical bijection $\varphi^\vee \colon \pi_1(\calRrhoss) \cong H^1_g(\GFS, \ad\rho(1))$, as required.   
\end{proof}
\section{Patching and Venkatesh's conjecture}\label{section_mainresults}
The remainder of this article is devoted to proving the following result.
\begin{theorem}\label{theorem_mainresult}
    Let $\Pi$ be a regular algebraic cuspidal automorphic representation of $\GL_n(\A_F)$, cohomological with respect to an algebraic representation $\Vlam$ where $\lambda\in(\Z^n)^{\Hom(F,E)}$. Let $\tau = (\tau_v)_{v\in S_p}=(\operatorname{LL}(\Pi_v))_{v\in S_p}$ denote the inertial type of $\Pi,$ and let $\rho = \rho_{\Pi,\iota} \colon \GFS \to \GL_n(E)$ denote the Galois representation attached to $(\Pi,\iota)$ where $\iota\colon\overline{\Q}_p \simto \C$ and $E/\Q_p$ is a finite extension containing all embeddings $F \hookrightarrow \overline{\Q}_p$. Suppose the following:
    \begin{enumerate}
        \item[(1)] $p \nmid 2n$.
        \item[(2)] $\Pi$ is tempered at almost all finite places.
        \item[(3)] $\zeta_p \notin F$ and $\rhobar|_{\Gamma_{F(\zeta_{p^\infty})}}$ has enormous image \cite[Def.~6.2.28]{10authorpaper}.
        \item[(4)] $\rhobar$ is absolutely irreducible and decomposed generic \cite[Def.~4.3.1]{10authorpaper}.
        \item[(5)] $F$ contains an imaginary quadratic subfield in which $p$ splits.
        \item[(6)] For every pair $\vbar\neq\vbar'$ of places of $F^+$ above $p$, we have (with $S^+_p$ the places of $F^+$ above $p$)
         \[
            \sum_{\wbar \in S^+_p \setminus \{\vbar,\vbar'\}} [F^+_{\wbar}\colon \Q_p] > \frac{1}{2}[F^+\colon \Q].
        \]
        \item[(7)] For every $v\in S$, the Weil-Deligne representation $\WD(\rho|_{\Gamma_{F_v}})$ is generic.
    \end{enumerate}
    Let $\calRrhoss\in\CNLE$ denote the derived deformation ring representing the characteristic 0 deformation problem
    $
        {(\rho,S,\{\calD_{\rho_v}^{\lambda_v}\}_{v\in S_p}\cup \{\calD_{\rho_v}^\square\}_{v\in S\setminus S_p})}.
    $
    Then there is a canonical isomorphism
    \[
        \pi_\ast(\calRrhoss) \cong \bigwedge H^1_g(\GFS, \ad\rho(1)),
    \]
    and there is an action of $\calRrhoss$ on $\RG(K,\Vlt)_\p$ inducing a free action
    \[
        \pi_\ast(\calRrhoss) \circlearrowright H_\ast(X_K,\Vlt)_\Pi.
    \]
\end{theorem}
\begin{remark}
    We briefly describe the roles played by the assumptions in \Cref{theorem_mainresult} and what is currently known.
    \begin{itemize}
        \item \emph{Taylor--Wiles method:} We require (1)-(4) to prove the existence of Taylor--Wiles primes in \Cref{prop_taylorwilesprimes} below. We do not expect (1) to be essential. The Ramanujan conjecture predicts that (2) always holds; the case $n=2$ is known by \cite{bcgnt2025} and the polarisable case holds by \cite{caraiani2012,caraiani2014}. One might expect to be able to replace the condition in (3) on $\rhobar|_{\Gamma_{F(\zeta_{p^\infty})}}$ with a `big image' assumption on the characteristic 0 representation, and to remove the `decomposed generic' part of (4).
        \item \emph{Local-global compatibility:}  Assumptions (5) and (6) are necessary to apply Hevesi's local-global compatibility result \cite[Thm.~1.4]{hevesi2024} which implies that $\rho_v$ has the expected Hodge--Tate weights; this compatibility is expected to hold in general --  when stronger results are available in the literature, we expect to be able to relax assumptions (5) and (6).
        \item \emph{Vanishing of the Bloch--Kato Selmer group / `$R=T$':} All of the assumptions (1)-(7) are required to apply the abstract patching lemma of A'Campo \cite[Lem.~2.6.4]{a'campo2024}, though our assumptions are somewhat weaker than those in \emph{loc. cit.}, owing to \cite{hevesi2024}. Assumption (7) would follow from a strong form of local-global compatibility at the ramified places. In the conjugate self-dual case, it follows from the main results of \cite{caraiani2014,caraiani2012}. 
    \end{itemize}
\end{remark}
\begin{remark}\label[remark]{rmk_noncanonical}
    A priori, the action of \Cref{theorem_mainresult} depends on non-canonical choices made in the Taylor--Wiles argument below. With stronger local assumptions on $\rhobar$ at $v\mid p$, this ambiguity is ruled out in \cite[§15]{galatius_venkatesh2018} by making a comparison between the mod~$p$ derived deformation ring $\calRS$ and the derived Hecke algebra introduced in \cite{venkatesh2019}.
\end{remark}

\textbf{From now on}, we fix $\Pi$ and $\rho$ as in \Cref{theorem_mainresult}, and an $\O$-lattice $\xi \colon \GFS \to \GL_n(\O)$ of $\rho$ with reductions $\xi_r\colon \GFS \to \GL_n(\O_r)$. The main result of \cite{hevesi2024} then implies the following.
\begin{theorem}\cite{hevesi2024}
    With assumptions as in \Cref{theorem_mainresult}, \[\bfS = (\rhobar, S, \{\calDvlambdatau\}_{v\in S_p} \cup \{\Dboxrhobarv\}_{v\in S\setminus S_p})\]
is a well-defined deformation problem. 
\end{theorem}
\subsection{Taylor--Wiles primes}
The basis of our patching argument later on is a delicate choice of Taylor--Wiles places, detailed in this subsection.
\begin{definition}
    Let $\bfS=(\rhobar,S,\{\calD_v\}_{v\in S})$ be a deformation problem. A Taylor--Wiles datum of level $N$ for $\bfS$ is a tuple $(Q_N,\{(\alpha_{v,1},\dots,\alpha_{v,n})\}_{v\in Q_N})$ where:
    \begin{enumerate}
        \item[(1)] $Q_N$ is a finite set of places of $F$, disjoint from $S$, such that, for every $v\in Q_N$, $v\equiv 1 \mod p^N$.
        \item[(2)] For every $v\in Q_N$, $(\alpha_{v,1},\dots,\alpha_{v,n})\in k^n$ are the eigenvalues of $\rhobar(\Frob_v)$, assumed to be pairwise distinct and $k$-rational.
    \end{enumerate}
\end{definition}
We generally omit $(\alpha_{v,i})$ and write simply $Q_N$ for Taylor--Wiles data. \Cref{theorem_mainresult} contains the assumption (3) that $\rhobar|_{\Gamma_{F(\zeta_{p^\infty})}}$ has enormous image. In \cite[Def.~2.23]{newton_thorne2023}, a notion of `enormous image' for characteristic 0 representations is defined, which is implied by the `enormous image' property of its reduction $\rhobar$, and features in the proof of \cite[Lem.~2.30]{newton_thorne2023}.
\begin{proposition}\label[proposition]{prop_taylorwilesprimes} Suppose the following:
    \begin{enumerate}
        \item[\textbullet] $p \nmid 2n$.
        \item[\textbullet] $\zeta_p \notin F$ and $\rhobar|_{\Gamma_{F(\zeta_{p^\infty})}}$ has enormous image.
        \item[\textbullet] $\Pi$ is tempered at almost all finite places. 
    \end{enumerate}
    Then there exist integers $q,l\geq 0$ such that, for any $N\geq 1$, there exists a Taylor--Wiles datum $Q_N$ for $\bfS'=(\rhobar,S,\{\Dboxrhobarv\}_{v\in S})$ such that
    \begin{enumerate}
        \item[(1)] $|Q_N|=q$.
        \item[(2)] For every $v\in Q_N$, $v\equiv 1 \mod p^N$ and the rational prime below $v$ splits in $F_0$.
        \item[(3)] There is a local $\O$-algebra surjection $\RSSloc[[X_1,\dots,X_g]]\surj \RSSQN$.
        \item[(4)] Let $\mathscr{L}_S=\{\mathscr{L}_{v,r}\}_{v\in S}$ be a Selmer system for $\ad\xi_r$. Then, for every $r\leq N$, with $SQ_N \coloneqq S \cup Q_N$ we have
        \[
            \operatorname{length}_\O H^1_{\mathscr{L}_{SQ_N}^\bot}\big(\GFSQN, \ad \xi_r(1)\big) \leq l.
        \]
        for $\mathscr{L}_{SQ_N}=\mathscr{L}_S \cup \{H^1(\GFv, \ad\xi_r)\}_{v\in Q_N}.$
    \end{enumerate}
\end{proposition}
\begin{proof} By \cite[Prop.~5.3.2]{a'campo2024}, for $q$ large enough we can find, for every $N\geq 1$, a Taylor--Wiles datum $(Q_N,\{(\alpha_{v,1},\dots,\alpha_{v,n})\}_{v\in Q_N})$ such that (1)-(3) hold. To achieve (4), we will enlarge the sets $Q_N$ by adding suitable Taylor--Wiles primes, following the proof of \cite[Lem.~2.30]{newton_thorne2023}. Thus, in (1) we replace if necessary $q$ by a larger number. The splitting condition in (2) defines a set of Dirichlet density $1$, so we may restrict our attention to places satisfying the condition. Enlarging $Q_N$ does not affect (3), as the statement is equivalent to the vanishing of
    \[
        H^1_{\mathscr{L}_{SQ_N}^\bot}(\GFSQN, \ad\xi_r(1)) = \ker\bigg(H^1(\GFS, \ad\rhobar(1))\to \prod_{v\in Q_N}H^1(\GFv, \ad\xi_r(1))\bigg)
    \]
    for $r=1$, and adding more factors to the product over $v$ preserves the vanishing of this kernel.

    Let $\ad_{E/\O} = \ad\xi \otimes_\O E/\O$ and similarly for $\ad_\O, \ad_E, \ad_{\O_r}$ so that $\ad_{E/\O} = \colim \ad_{\O_r}$. Fix
    \[
        s \geq \operatorname{corank}_\O H^1(\GFS, \ad_{E/\O}(1)) \coloneqq \rank_\O \Hom_\O\left(H^1(\GFS, \ad_{E/\O}(1)),E/\O\right).
    \]
    We will prove the existence of Frobenius elements $\sigma_1,\dots,\sigma_s \in G_{F(\zeta_{p^\infty})}$ such that $\rho(\sigma_i)$ has distinct eigenvalues in $E$ and the kernel of the map
        \[
            H^1\big(\GFS, \ad_{E/\O}(1)\big) \to \bigoplus_{i=1}^s H^1\big(\langle \sigma_i \rangle, \ad_{E/\O}(1)\big) \cong \bigoplus_{i=1}^s \frac{\ad_{E/\O}(1)}{(\sigma-1)\ad_{E/\O}(1)}
        \]
        is of finite length as an $\O$-module. Given such $\sigma_i\in G_{F(\zeta_{p^\infty})}$, for every $r\leq N$, consider the commutative diagram with exact columns
\[\begin{tikzcd}
	0 & 0 \\
	{H^0\big(\GFS, \ad_{E/\O}(1)\big)/\varpi^r} & {\bigoplus_{i=1}^{s} H^0\big(\langle\sigma_i\rangle, \ad_{E/\O}(1)\big)/\varpi^r} \\
	{H^1\big(\GFS, \ad_{\O_r}(1)\big)} & {\bigoplus_{i=1}^{s} H^1\big(\langle\sigma_i\rangle, \ad_{\O_r}(1)\big)} \\
	{H^1\big(\GFS, \ad_{E/\O}(1)\big)[\varpi^r]} & {\bigoplus_{i=1}^{s} H^1\big(\langle\sigma_i\rangle, \ad_{E/\O}(1)\big)[\varpi^r]} \\
	0 & {0.}
	\arrow[from=1-1, to=2-1]
	\arrow[from=1-2, to=2-2]
	\arrow[from=2-1, to=2-2]
	\arrow[from=2-1, to=3-1]
	\arrow[from=2-2, to=3-2]
	\arrow[from=3-1, to=3-2]
	\arrow[from=3-1, to=4-1]
	\arrow[from=3-2, to=4-2]
	\arrow[from=4-1, to=4-2]
	\arrow[from=4-1, to=5-1]
	\arrow[from=4-2, to=5-2]
\end{tikzcd}\]
Note that $H^0(\GFS, \ad_{E/\O}(1))$ has finite length by irreducibility of $\rho$. It follows from the snake lemma that the kernel of
\[
    H^1(\GFS, \ad_{\O_r}(1)) \to \bigoplus_{i=1}^s H^1(\langle \sigma_i\rangle, \ad_{\O_r}(1))
\]
has length bounded independently of $r$ and $N$. By the Chebotarev density theorem, we can then find places $v_1,\dots v_s$ of $F$ with Frobenius elements $\Frob_{v_1},\dots,\Frob_{v_s}$ approximating $\sigma_1,\dots,\sigma_s$, and thus we obtain statement (4) of the proposition. Now, we have an isomorphism
\[
    H^1(\GFS, \ad_{E/\O}(1)) \cong (E/\O)^{s_0} \oplus M 
\]
where $M$ is a finite length $\O$-module and $s_0\leq s$ is the corank defined above. Thus, the problem reduces to finding, for every non-zero
\[
    f \colon E/\O \to H^1(\GFS, \ad_{E/\O}(1)),
\]
an element $\sigma\in G_{F(\zeta_{p^\infty})}$ such that $\rho(\sigma)$ has $n$ distinct eigenvalues and the composition
\[
    E/\O \to H^1(\GFS, \ad_{E/\O}(1)) \to H^1(\langle\sigma\rangle, \ad_{E/\O}(1)) \cong \frac{\ad_{E/\O}(1)}{(\sigma-1)\ad_{E/\O}(1)}
\]
is non-zero. Let $F_\infty = F(\zeta_{p^\infty})$ and $L_\infty = \overline{F}^{\ker\ad_{E}(1)}$. Then $F_\infty \subseteq L_\infty$ since $E(1)\subset \ad_E(1)$. The short exact sequence
\[
    0 \to \ad_\O \to \ad_E \to \ad_{E/\O} \to 0
\]
induces an exact sequence (with $\Gamma(L_\infty/F) = \Gal(L_\infty/F)$)
\[
    H^1\big(\Gamma(L_\infty/F), \ad_E(1)\big) \to H^1\big(\Gamma(L_\infty/F), \ad_{E/\O}(1)\big) \to H^2\big(\Gamma(L_\infty/F), \ad_\O(1)\big).
\]
The first term vanishes by \cite[Lem.~6.2]{kisin2004} and the assumption that $\Pi$ is tempered at all finite places. The third term is a finitely generated $\O$-module. It follows that the middle term is annihilated by a power of $p$. Thus, the exactness of the inflation-restriction sequence
\[
    0 \to H^1\big(\Gamma(L_\infty/F), \ad_{E/\O}(1)\big) \to H^1\big(\GFS, \ad_{E/\O}(1)\big) \to H^1\big(\Gamma(F_S/L_\infty), \ad_{E/\O}(1)\big)^{\Gamma(L_\infty/F)}
\]
implies that the composition
\begin{equation*}
    \begin{split}
        E/\O \stackrel{f}{\to} 
        H^1\big(\GFS, \ad_{E/\O}(1)\big) 
        &\to H^1\big(\Gamma(F_S/L_\infty), \ad_{E/\O}(1)\big)^{\Gamma(L_\infty/F)}\\
        &\cong \Hom_{\Gamma(F_S/F)}\big(\Gamma(L_{\infty,S_{L_\infty}}/L_\infty), \ad_{E/\O}(1)\big)
    \end{split}
\end{equation*}
is non-zero. Let $W \subset \ad_{E/\O}(1)$ be the $\O$-submodule spanned by $\{f(x)(\sigma) \mid x\in E/\O,\sigma\in G_{L_\infty}\}$. Then $W$ is non-zero and divisible. By \cite[Lem.~2.22]{newton_thorne2023} there exists $\sigma \in G_{F_\infty}$ such that $\rho(\sigma)$ has $n$ distinct eigenvalues in $E$ and $W \not\subset (\sigma-1)\ad_{E/\O}(1)$. Thus, for some $m \geq 0$ and $\tau \in G_{L_\infty}$ we have
\[
    f(1/\varpi^m)(\tau) \notin (\sigma-1)(\ad_{E/\O}(1)).
\]
If $f(1/\varpi^m)(\sigma) \notin (\sigma-1)(\ad_{E/\O}(1))$ then we are done. If not, then $\tau\sigma$ has the sought property, i.e.
\[
    E/\O \to H^1(F_S, \ad_{E/\O}(1)) \to \frac{\ad_{E/\O}(1)}{(\tau\sigma-1)\ad_{E/\O}(1)}
\]
is non-zero, since $\tau$ acts trivially on $\ad_E(1)$.
\end{proof}

\subsection{Patching}\label[section]{ssec_patching}
In this section, we set up the patching argument, basically following \cite{a'campo2024} and \cite{galatius_venkatesh2018} but with some adjustments. Define deformation problems
\begin{align*}
    \bfS &\coloneqq \big(\rhobar, S, \{\calDvlambdatau\}_{v\in S_p} \cup \{\Dboxrhobarv\}_{v\in S\setminus S_p}\big) \\
    \bfS' &\coloneqq \big(\rhobar, S, \{\Dboxrhobarv\}_{v\in S}\big)
\end{align*}
We let $\calRS,\calRSpr \in \CNLk$ denote the corresponding derived deformation rings.

For a Taylor--Wiles-datum $Q$, we define an augmented deformation problem
\[
    \bfS_{Q} = \big(\rhobar, S\cup Q, \{\calDvlambdatau\}_{v\in S_p} \cup \{\Dboxrhobarv\}_{v\in S\setminus S_p} \cup \{\Dboxrhobarv\}_{v\in Q}\big).
\]
and similarly for $\bfS'_Q$. Applying \Cref{prop_taylorwilesprimes}, we obtain $q,l \geq 0$ and for every $N\geq 1$ a Taylor--Wiles datum $Q_N$ satisfying the conclusions of the proposition. We let
\[
    \calRN \coloneqq \calR_{\bfS_{Q_N}}^{S}\ \text{ and }\ \calRNpr \coloneqq \calR_{\bfS'_{Q_N}}^S
\]
be the $S$-framed type-$\bfS_{Q_N}$ and type-$\bfS'_{Q_N}$ derived deformation rings. We also let 
\[
    \RNSSQN \coloneqq \pi_0(\calR_{\bfS_{Q_N}}^{S})\ \text{ and }\ \RNSSQNpr \coloneqq \pi_0(\calR_{\bfS'_{Q_N}}^S)
\]
denote the underlying classical deformation rings. These objects and the ones introduced below depend on the specific choice of $Q_N$ and not just $N$, but in our applications, we will fix a choice of $Q_N$ for every $N\geq 1$, motivating our simplified notation.

In \cite{galatius_venkatesh2018}, the authors consider deformations of $\rhobar_v$ for $v\in Q_N$ valued in the diagonal torus $\bfT \subset  \bfG$. The local $\bfT$-valued deformation functors at Taylor--Wiles places are proved in \cite[Lem.~8.3]{galatius_venkatesh2018} to be equivalent to their $\bfG$-valued counterparts, and we will not make a distinction.

We let $\calT \coloneqq \O[[T_1,\dots,T_{n^2|S|-1}]]$ be the power series ring in the `$S$-frame variables' (\Cref{prop_Tframeddimensions}) and define, for a Taylor--Wiles datum $Q_N$,
\begin{align*}
    \calS_v &\coloneqq \text{representing object of }\Dboxrhobarv, \ v\in Q_N, \\
    \calS_N &\coloneqq  \big(\widehat{\bigotimes}_{v\in Q_N} \calS_v \big)\widehat{\otimes}_{\O} \calT \in \CNLk.\\
    \calS_v^{\operatorname{ur}} &\coloneqq \text{representing object of }(\calX^\square_{\GFv/\IFv,\GL_n})_{\rhobarv}, \ v\in Q_N, \\
    \calSNur &\coloneqq \widehat{\bigotimes}_{v\in Q_N} \calS_v^{\operatorname{ur}}.
\end{align*}
\begin{lemma}\label[lemma]{lemma_RNSNSNur}\cite[Eq.~11.4]{galatius_venkatesh2018}
    For every $N\geq 1$, we have natural maps $\calSN \to \calRNpr \to \calRN$ and $\RN'$ and $\RN$ are local $\SN$-algebras. Moreover, there are equivalences
    \[
        \calRS \simeq \calRN \otimesL_{\calSN} \calSNur \ \text{ and }\ \calRSpr \simeq \calRNpr \otimesL_{\calSN} \calSNur.
    \]
\end{lemma}
The rings $\calSN$ and $\calSNur$ are in fact static, and the latter is equivalent to a power series ring over $\O$ (see \cite[Lem.~8.6]{galatius_venkatesh2018}). It follows that the descent formula afforded by \Cref{lemma_RNSNSNur} can be simplified by, informally speaking, cancelling these power series variables.
\begin{lemma}\label[lemma]{lemma_Srings}
    Let $Q_N$ be a Taylor--Wiles datum and define
    \[
    \SNo \coloneqq \widehat{\bigotimes}_{v\in Q_N} \O[[\bfT(k(v))_p)]],
    \]
    where $\bfT\subset \bfG$ is the diagonal torus. Then we have the following:
    \begin{itemize}
    \item[(1)] There is a natural morphism $\SNo \to \calSN$ which is formally smooth. That is, $\calSN\simeq \pi_0(\calSN)$ is equivalent to a power series ring over $\SNo$.
    \item[(2)] There is an equivalence $\calSNur \simeq \calSN\otimesL_{\SNo} \O$. Consequently, there are equivalences
    \[
        \calRS \simeq \calRN \otimesL_{\SNo} \O \ \text{ and }\ \calRSpr \simeq \calRNpr \otimesL_{\SNo} \O.
    \]
    \item[(3)] Let $r< s(N)$. Then, for the augmentation $\SNo \to \O_r$,
    \[
        \opT^1_{\O_r}(\SNo) \cong \prod_{v\in Q_N} H^{2}(\GFv, \ad\xi_r).
    \]
    \end{itemize}
\end{lemma}
\begin{proof}
    (1) and (2) are contained in \cite[Lem.~8.6]{galatius_venkatesh2018}. (3) follows from \Cref{prop_tangentcomplexesgaloiscohomology} and (1).
\end{proof}
As is typical in the Taylor--Wiles method, we will pass to a certain limit over $N$ and thus obtain objects of a simple form. We define
\begin{align*}
    \Soinf &\coloneqq \calT[[Y_{1},\dots,Y_{qn}]], \\
    \mathfrak{a}_{N} &= \big\langle (1+Y_i)^{p^{N}}-1,(1+T_i)^{p^{N}}-1\big\rangle \subset \Soinf, \\
    \mathfrak{c}_{N} &= \big\langle \varpi^N, \mathfrak{a}_{N} \rangle \subset \Soinf,\\
    \SNobar &\coloneqq \Soinf/\mathfrak{c}_N, 
\end{align*}
so that $\Soinf\cong \varprojlim_N \SNobar$ as profinite topological rings. For every $N\geq 1$, we fix an isomorphism
$
    \SNobar \cong \SNo/\mathfrak{c}_N.
$
and a surjection (whose existence is guaranteed by \Cref{prop_taylorwilesprimes})
\[
    \Rinf \coloneqq \RSSloc[[X_1,\dots,X_{nq-n^2[F^+\colon\Q]}]] \surj \RN = \RSSQN.
\]
Note that $\frakm_{\RN} = \frakm_{\Rinf}R_N$ by surjectivity; we will write simply $\frakm$ for the maximal ideal of $\RN$.
\begin{theorem}\label{theorem_acampo}\cite{a'campo2024}
    With the notation above and under the assumptions of \Cref{theorem_mainresult}, we have the following:
    \begin{enumerate}
        \item[(i)]
        There exists an increasing sequence $\N^\ast = \{N_s\}_{s\geq 2} \subseteq \N_{\geq 2}$ such that, with
        \begin{align*}
            \RNbar &\coloneqq \RN/(\fraka_{N}R_{N} + \frakm^{s(N)}), \\
             s(N) &\coloneqq \#\{n \in \N^\ast\colon n\leq N\}+1,
        \end{align*}
        there are isomorphisms of $R_\infty$-algebras
        $
            \RNbar/\frakm^{s(M)} \cong \RMbar
        $
        for $M\leq N$ in $\N^\ast$. Moreover, there is an isomorphism
        \[
            \Rinfphat \simto \big(\varprojlim_{N\in\mathbb{N}^\ast} \RNbar\big)_\p^\wedge, 
        \]
        where $\p=\ker(\Rinf \to \O)$ corresponds to $\xi$ and the transition maps are given by
        \[
            \RNbar \surj \RNbar / \frakm^{s(M)} \cong \RMbar.
        \]
        \item[(ii)] There exists an $R_\infty$-module $M_\infty$ such that $\Minfphat$ is free over $\Rinfphat$ and
        \[
            \Minfphat \otimesL_{\Sinfahat} E \simeq  \RHom_E\big(\RG(X_K,\Vlt)_\p, E\big).
        \]
        Consequently, the graded $E$-algebra
        \[
            \pi_\ast\big(\Rinfphat \otimesL_{\Sinfahat} E\big) \cong \Tor_\ast^{\Sinfahat}\big(\Rinfphat,E\big)
        \]
        acts freely on the graded $E$-vector space
        \[
            \pi_\ast\big(\Minfphat \otimesL_{\Sinfahat} E\big) \cong H_\ast(X_K,\Vlt)_\Pi.
        \]
        \item[(iii)] The geometric adjoint Bloch--Kato Selmer group of $\rho$ vanishes, i.e.
        \[
            H^1_g(\GFS, \ad\rho) = 0.
        \]
    \end{enumerate}
\end{theorem}
\begin{proof}
    All three essentially follow from the proof of \cite[Thm.~5.3.5]{a'campo2024}. Since our local conditions at $v\mid p$ differ slightly, we spell out the proof which is a straightforward application of the abstract patching lemma \cite[Lem.~2.6.4]{a'campo2024}. For $N\geq 1$, one defines (see \cite[Prop.~5.3.3]{a'campo2024} for details):
    \begin{itemize}
        \item An open compact subgroup $K_0(Q_N) = \prod_{v} K_{0,v}(Q_N) \subset K$ where $K_{0,v}(Q_N)=K_v$ for $v\notin Q_N$ and $K_{0,v}(Q_N)=\operatorname{Iw}_v$ is the matrices which are upper triangular mod~$\varpi_v$ for $v\in Q_N$.
        \item An open compact subgroup $K_1(Q_N) = \prod_{v} K_{1,v}(Q_N) \subset K$ where $K_{1,v}(Q_N) = K_v$ for $v\notin Q_N$ and $K_{1,v}(Q_N)=\operatorname{Iw}_{1,v} \subset \operatorname{Iw}_v$ is the matrices which are upper triangular and unipotent mod~$\varpi_v$ for $v\in Q_N$.
        \item For any object $M$ with a $\TS$-action, an enlarged (faithful) Hecke algebra $\TSQNQN(M)$ obtained by adjoining Hecke operators at the places $v\in Q_N$.
        \item Maximal ideals $\frakm_1,\frakm_0$ such that, if $\frakm\subset \TS(\RG(X_K,\Vlt))$ is the ideal corresponding to $\Pi$, we have morphisms
        \[
        \TSQNQN\big(\RG(X_{K_1(Q_N)},\Vlt)\big)_{\frakm_1} \to \TSQN\big(\RG(X_{K_0(Q_N)},\Vlt)\big)_{\frakm_0} \simto \TS\big(\RG(X_{K},\Vlt)\big)_{\frakm}.
        \]
        where the first is given by taking $K_0(Q_N)/K_1(Q_N)$-invariants, and corresponding quasi-isomorphisms of complexes
        \[
            \RG\big(K_0(Q_N)/K_1(Q_N),\RG(X_{K_1(Q_N)},\Vlt)_{\frakm_1}\big) \simto \RG(X_{K_0(Q_N)},\Vlt)_{\frakm_0} \simto \RG(X_K,\Vlt)_{\frakm}.
        \]
    \end{itemize}
    In the numbering of \cite[Lem.~2.6.4]{a'campo2024}, we have 7 statements to verify.
    \begin{itemize}
        \item[(1)] $R_\infty$ is a complete noetherian local $\O$-algebra.
        \item[(2)] $R_N/\fraka_N$ is a local $\Soinf/\fraka_N$-algebra and a quotient of $R_\infty$.
        \item[(3)] Define $C_N \in \AniMod^{\operatorname{cn}}_\O$ for $N\geq 0$ by
        \begin{align*}
            C_0 &\coloneqq \RHom_\O\left(R\Gamma(X_K,\Vlt)_\frakm,\O\right), \\
            C_N &\coloneqq \RHom_\O\left(R\Gamma(X_{K_1(Q_N)},\Vlt)_{\frakm_1},\calT\right)/\fraka_N,   
        \end{align*}
        where the right-hand sides denote minimal representatives in the derived category of $\O$ and $\Soinf/\fraka_N$, respectively. Then $C_N/\fraka \simeq C_0$ by the fourth bullet point above.
        \item[(4)] By \Cref{prop_spreadingout}, $C_0[1/\varpi]$ has homology concentrated in degrees $[q_0,q_0+l_0]$ and is not exact. Moreover, $\dim (\Soinf)_{\fraka} = n^2|S|-1+qn$, and from \Cref{prop_dimension_vnmidp} and \Cref{theorem_semistabledeformationrings}(2) we deduce \[
            \dim R_\infty[1/\varpi] \leq n^2|S|+qn-n[F^+\colon \Q] = \dim (\Soinf)_{\fraka} - l_0.
        \]
        \item[(5)] For $N\geq 1$, we let $\frakm_1,\frakm_0$ denote the ideals in the fourth bullet point above, and define
        \[
            T_N \coloneqq \bbT^{S\cup Q_N}_{Q_N}(C_N)_{\frakm_1}\cdot \calT \subset \End_{D(\Soinf/\fraka_N)}(C_N),
        \]
        so that the image of $T_N$ in  $\End_{D(S_0)}(C_0)$ equals 
        \[
            \bbT^{S\cup Q_N}(C_0) \subset \bbT^{S}(C_0) \eqqcolon T_0.
        \]
        \item[(6)] Using \Cref{prop_Tframeddimensions} and \cite[Thm.~2.24]{hevesi2024}, we obtain for every $N$ an $\Soinf/\fraka_N$-algebra surjection
        $
            R_N \to T_N/I_N
        $
        (after taking a tensor product with $\calT$),
        where $I_N^\delta = 0$ for some $\delta$ independent of $N$, obtaining a commutative square
\[\begin{tikzcd}
	{R_N} & {T_N/I_N} \\
	{R_0} & {T_0/I_0}
	\arrow[two heads, from=1-1, to=1-2]
	\arrow[from=1-1, to=2-1]
	\arrow[from=1-2, to=2-2]
	\arrow[two heads, from=2-1, to=2-2]
\end{tikzcd}\]
        \item[(7)] Since the Weil--Deligne representations $\WD(\rho_v)$ are generic for every $v\in S$, the representation $\rho$ defines a smooth point of $R_\infty$ by \Cref{prop_dimension_vnmidp,theorem_semistabledeformationrings}(3).
    \end{itemize}
    Now, \cite[Lem.~2.3.3, Lem.~2.6.4]{a'campo2024} provide a subsequence $\N^\ast \subseteq \N$ indexed by $s\colon \N^\ast\to \N_{\geq 2}$ and morphisms $\RNbar \to \RMbar$ for every $M\leq N$ in $\N^\ast$ as in the statement, and an isomorphism as in (i): in the notation of \textit{loc. cit.}, 
    \[
        \Rinfphat \cong \mathbf{P}_{\Rinf}\big((\RN)_{N\geq 1}\big)_\p^\wedge \cong (\varprojlim_N \RNbar)_\p^\wedge,
    \]
    where $\mathbf{P}_{\Rinf}$ is the patching functor. Statements (ii) and (iii) follow from \cite[Thm.~5.3.5]{a'campo2024} and \cite[Lem.~2.6.4]{a'campo2024}, respectively.
\end{proof}
In \Cref{theorem_acampo}, the subset $\N^\ast \subseteq \N$ and the maps $\RNbar \to \RMbar$ depend on a choice of non-principal ultrafilter on $\N$. As such, the isomorphism provided by (1) is non-canonical. By construction, we have $s(N)\leq N$, so that we have points $\xi_r \colon \RNbar \to \O_r$ for $r<s(N)$.

\textbf{From now on}, we fix a choice of ultrafilter and consider the objects afforded by \Cref{theorem_acampo} to be fixed. In the sequel, limits and colimits over the variables $N$ and $M$ are taken over $\N^\ast$.

In light of part (ii) of \Cref{theorem_acampo}, we see that to define an action as in \Cref{theorem_mainresult}, it suffices to prove an equivalence
\[
    (\calRS \otimes_\O E)_\xi^\wedge \simeq \Rinfphat \otimesL_{\Sinfahat} E,
\]
and to identify the right-hand side with the exterior algebra on the Bloch--Kato Selmer group.

\subsection{The approximating deformation functor}\label{section_GVextended}
Our proof of \Cref{theorem_mainresult} will rely on constructing rings which approximate $\calRS$, at least locally around the point corresponding to $\rho$. In this section, we construct these approximating rings and prove a few statements about their tangent complexes.

For every $M\leq N$ in $\N^\ast$, \Cref{theorem_acampo}(1) provides a commutative diagram
\[\begin{tikzcd}
	\RNbar & \SNobar & \ON \\
	\RMbar & \SMobar & \OM.
	\arrow[from=1-1, to=2-1]
	\arrow[from=1-2, to=1-1]
	\arrow[from=1-2, to=1-3]
	\arrow[from=1-2, to=2-2]
	\arrow[from=1-3, to=2-3]
	\arrow[from=2-2, to=2-1]
	\arrow[from=2-2, to=2-3]
\end{tikzcd}\]
where $\RNbar$ and $\SNobar$ are the artinian quotients introduced in \Cref{ssec_patching}. We let 
\[
  \calC_N \coloneqq \RNbar \otimesL_{\SNbar} \SNurbar    
\]
and denote by $e_{N,M}\colon\calC_N \to \calC_M$ the map induced by the diagram. Note that there is a morphism
\[
    \calRS \stackrel{4.2.2}{\simeq} \calRN \otimesL_{\SNo} \O \to \RNbar \otimesL_{\SNobar} \ON.
\]
Roughly speaking, we will show that, as $N\to \infty$, the right-hand converges to a ring which approximates $\calRS$ near $\rho$. In order for this to be a viable strategy, the artinian quotients $\RNbar$ and $\SNobar$ ought to be somewhat large. The following proposition makes this precise.
\begin{proposition}\label[proposition]{prop_RNSNtangentspaces}
    Let $r<s(N)$. Then the surjections $\RN \surj \RNbar$ and $\SNo \surj \SNobar = \SNo/\mathfrak{c}_N$ induce isomorphisms on tangent spaces at the point corresponding to $\xi_r$; that is,
    \[
        \opT^0_{\O_r}(\RNbar) \cong \opT^0_{\O_r}(\RN) \text{ and } \opT^0_{\O_r}(\SNobar) \cong \opT^0_{\O_r}(\SNo).
    \]
\end{proposition}
\begin{proof}
    The maps in question are the quotients
    \[
        \RNbar = \RN/(\fraka_N,\frakm^{s(N)}),\quad \SNobar = \SNo/\mathfrak{c}_N 
    \]
    The natural injective map on tangent spaces is given by precomposition of derivations:
    \[
        \operatorname{Der}_{\O}(\RNbar,\O_r) \hookrightarrow \operatorname{Der}_{\O}(\RN,\O_r).
    \] 
    For the other statement, it suffices to show that any derivation $d\colon \RN \to \O_r$ annihilates $\fraka_NR_N$ and $\frakm^{s(N)}$. We have homomorphisms
    \[
        \Soinf \to \RN \to \RNbar \to \O_r
    \]
    where the composite is augmentation map $\Soinf\to\O_r$. Reduced mod~$\varpi^r$, the generators $(1+Y_i)^{p^N}-1$ of $\fraka_N\subset \Soinf$ have only quadratic and higher terms, hence the containment
    $
        \fraka_N \subseteq \ker(\Soinf \to \O_r)^2
    $
    and consequently
    \[
        \fraka_NR_N \subseteq \ker(\RN \to \O_r)^2.
    \]
    The Leibniz rule now forces any derivation $d\colon \RN \to \O_r$ annihilates $\fraka_N R_N$. Moreover, the fixed map $\xi_r \colon \RN \to \O_r$ factors through $\RN/\frakm^{r}$, so that any derivation $d\colon \RN \to \O_r$ factors through $\RN/\frakm^{r+1}$ by the Leibniz rule. Since $r+1\leq N$, this completes the proof for $\RN\surj \RNbar$. The argument for $\SNo\surj \SNobar$ is similar. 
\end{proof}

\begin{lemma}\label[lemma]{lemma_tgtcomplexcalculation}
    With notation as above, let
    \[
        \FN = \Map_{/k}(\calC_N,-) \colon \Artk \to \Ani
    \]
    and $\Finfty = \colim_N \FN$. For our fixed representation $\xi\in \Finfty(\O)$, we have the following:
    \begin{enumerate}
    \item [(1)] For every $i$ and every $r\geq 1$, $\opT^i_{\O_r}(\Finfty)$ is a finitely generated $\O_r$-module.
    \end{enumerate}
    \begin{enumerate}
    \item [(2)] The functor $(\Finfty)_{(\xi)}$ (see \Cref{def_Frhobarxi}) is represented by $\Rinfphat \otimesL_{\Sinfahat} E \in \CNLE$.
    \item [(3)] The tangent complex of $(\Finfty)_{(\xi)}$ has cohomology concentrated in degrees $[0,1]$ and satisfies
    \[
        \dim_E \opT^0_E\big((\Finfty)_{(\xi)}\big) - \dim_E \opT^1_E\big((\Finfty)_{(\xi)}\big) =   -l_0.
    \]
    \end{enumerate}
\end{lemma}
\begin{remark}\label[remark]{rmk_Cinfty}
    Due to the lack of information about the rings $\RNbar$, it is unclear to the author whether $\Finfty$ is necessarily represented by $\calC_\infty \coloneqq \lim_N\ \calC_N$, i.e. whether $\{\calC_N\}_{N\geq 1}\in\calPro(\Artk)$ is equivalent to the pro-artinian completion $\{\tau_{\leq j}(\calC_\infty\slashL \frakm_{\pi_0(\calC_\infty)}^r)\}_{j,r} \in \calPro(\Artk)$.
\end{remark}
\begin{proof}
    The fact that filtered colimits commute with finite limits implies the equivalence of functors
    \[
        \Finfty \simeq \big(\underset{N}{\colim} \Map_{/k}(\RNbar,-)\big) \times_{\colim_N \Map_{/k}(\SNobar,-)} \big(\colim_N \Map_{/k}(\ON,-)\big). 
    \]
    Now, by  \cite[\href{https://stacks.math.columbia.edu/tag/06SF}{Tag 06SF}]{stacks-project} and the fact that $\Rinf$ surjects onto $\RNbar$ for every $N\geq 1$, the $\frakm$-adic topology on the inverse limit $\varprojlim \RNbar$ can be defined using the system of ideals $\ker(\varprojlim \RNbar \to \RNbar)$. Similarly, $\Soinf\cong \varprojlim \SNobar$ is an isomorphism of topological rings. It then follows from \cite[Prop.~2.1.4]{halpern-leistner_preygel2023} that there are equivalences of functors
    \begin{align*}
        \underset{N}{\colim} \Map_{/k}(\RNbar, -) &\simeq \Map_{/k}(\varprojlim \RNbar, -), \\
        \underset{N}{\colim} \Map_{/k}(\SNobar, -) &\simeq \Map_{/k}(\Soinf, -),
    \end{align*}
    and similarly for $\O$. Thus, we have proved
    \[
        \Finfty \simeq \Map_{/k}(\varprojlim \RNbar,-) \times_{\Map_{/k}(\Soinf,-)}\Map_{/k}(\O,-).
    \]
    It follows that there is a fibre sequence of tangent complexes
    \[
        \opT^0_{\O_r}(\Finfty) \to \opT^\bt_{\O_r}(\varprojlim \RNbar) \to \opT^\bt_{\O_r}(\Soinf) \stackrel{+1}{\to}.
    \]
    Since the rings $\varprojlim \RNbar$ and $\Soinf$ are noetherian, \Cref{prop_Tifinitelygenerated} and the fibre sequence above implies that $\opT^i_{\O_r}(\Finfty)$ is a finitely generated $\O_r$-module. 
    
    To prove (2), we use the isomorphism $\Rinfphat \cong (\varprojlim \RNbar)_\p^\wedge$ of  \Cref{theorem_acampo}(1). Applying \Cref{theorem_charptochar0rings,lemma_xihatpullbacks,prop_limcolim}, we obtain an equivalence
    \[
        (\Finfty)_{(\xi)} \simeq \Map_{/E}(\Rinfphat,-) \times_{\Map_{/E}(\Sinfahat,-)} \Map_{/E}(E,-) \simeq \Map_{/E}(\Rinfphat \otimesL_{\Sinfahat} E,-).
    \]
    Thus, we obtain a fibre sequence of tangent complexes over $E$:
    \[
        \opT^\bt_E\big((\Finfty)_{(\xi)}\big) \to \opT^\bt_E(\Rinfphat) \to \opT^\bt_E(\Sinfahat) \stackrel{+1}{\to},
    \]
    which induces an exact sequence
    \[
        0 \to \opT^0_E\big((\Finfty)_{(\xi)}\big) \to \opT^0_E(\Rinfphat) \to \opT^0_E(\Sinfahat) \to \opT^1_E\big((\Finfty)_{(\xi)}\big) \to 0,
    \]
    where we have used that $\Rinfphat$ is formally smooth by \Cref{theorem_semistabledeformationrings}(3) and \Cref{prop_dimension_vnmidp}, under the genericity assumption on $\WD(\rho_v)$ for every $v\in S$. The higher terms of the same exact sequence yield isomorphisms $\opT^i_E\big((\Finfty)_{(\xi)}\big) \cong 0$ for $i\neq 0,1$, which completes the proof.
\end{proof}

\begin{proposition}\label[proposition]{prop_T1RN}
    Let $r< s(N)$. With assumptions and notation as above, the kernel
    \[
        \ker\big( \opT^1_{\O_r}(\calRN) \to \opT^1_{\O_r}(\SNo)\big)
    \]
is a finite length $\O$-module with length bounded independently of $r$ and $N$.
\end{proposition}
\begin{proof}
    To begin, we prove the same statement for $\calRNpr$, the ($S$-framed) derived deformation ring with trivial local conditions at all $v\in S$. That is, we let $r< s(N)$ and prove that the kernel
    \[
        \ker\big( \opT^1_{\O_r}(\calRNpr) \to \opT^1_{\O_r}(\SNo)\big)
    \]
is a finite length $\O$-module with length bounded independently of $r$ and $N$. By \Cref{prop_tangentcomplexesgaloiscohomology}, \Cref{lemma_Srings} and the discussion in \Cref{ssec_localconditions}, the map in question is the product of restrictions
\[
    \opT^1_{\O_r}(\calRNpr) \cong H^{2}(\GFSQN, \ad\xi_r) \to \prod_{v\in Q_N} H^{2}(\GFv, \ad\xi_r) \cong
        \opT^1_{\O_r}(\SNo)  
\]
    We factor $\opT^1_{\O_r}(\calRNpr) \to \opT^1_{\O_r}(\SNo)$ as $\beta\circ \alpha$, where
    \[
        H^2(\GFSQN, \ad\xi_r) \stackrel{\alpha}{\to} \prod_{v\in SQ_N} H^2(\GFv, \ad\xi_r) \stackrel{\beta}{\to} \prod_{v\in Q_N} H^2(\GFv, \ad\xi_r)
    \]
    to obtain an exact sequence
    \[
        0 \to \ker(\alpha) \to \ker(\beta \circ \alpha) \to \ker(\beta).
    \]
    It then suffices to show that $\ker(\alpha)$ and $\ker(\beta)$ have finite length bounded independently of $r$ and $N$. To see this, consider the Selmer system \[\{\mathscr{L}_{v,r}\}_{v\in SQ_N} = \{H^1(\GFv, \ad\xi_r)\}_{v\in SQ_N}\cup \{H^1_{\operatorname{ur}}(\GFv, \ad\xi_r)\}_{v\notin SQ_N}.\]
    It follows from Poitou--Tate duality that we have an exact sequence 
    \[
        H^1(\GFSQN, \ad\xi_r(1)) \to H^2(\GFSQN, \ad\xi_r) \stackrel{\alpha}{\to} \prod_{v\in SQ_N} H^2(\GFv, \ad\xi_r).
    \]
    In particular, we have 
    \[
        \operatorname{length}_\O  \ker(\alpha) \leq \operatorname{length}_\O H^1(\GFSQN, \ad\xi_r(1)),
    \]
    which is bounded from above independently of $r,N$ by \Cref{prop_taylorwilesprimes}(4). Now,
    \[
        \ker(\beta) = \prod_{v\in S} H^2(\GFv, \ad\xi_r) \cong \opT^1_{\O_r}(\RSprSloc),
    \]
    and since the characteristic 0 point of $\Spec \RSprSloc$ corresponding to $\xi$ is smooth by \Cref{prop_dimension_vnmidp} and the genericity assumption at all $v\in S$, $\ker(\beta)$ has finite length bounded independently of $r$ (and $N$) by \Cref{prop_Tifinitelengthsmoothpoint}. This completes the proof of the initial claim. Consider now the kernel
    \[
        \ker\big( \opT^1_{\O_r}(\calRN) \to \opT^1_{\O_r}(\SNo)\big)
    \]
    in the proposition. The equivalence
$
    \calRN \simeq \calRNpr \otimesL_{\RSprSloc} \RSSloc
$
induces a commutative diagram
\[\begin{tikzcd}
	0 & Q_{r,N} & {\opT^1_{\O_r}(\calR_N)} & {\opT^1_{\O_r}(\calR_N^\prime) \oplus \opT^1_{\O_r}(\RSSloc)} \\
	&& {\opT^1_{\O_r}(\SNo)} & {\opT^1_{\O_r}(\calR_N^\prime)}
	\arrow[from=1-1, to=1-2]
	\arrow[from=1-2, to=1-3]
	\arrow[from=1-3, to=1-4]
	\arrow["{\gamma}", from=1-3, to=2-3]
	\arrow[two heads, from=1-4, to=2-4]
	\arrow[from=2-4, to=2-3]
\end{tikzcd}\]
where 
\[
    Q_{r,N} \coloneqq \coker\big(\opT^0_{\O_r}(\calR_N^\prime) \oplus \opT^0_{\O_r}(\RSSloc) \to \opT^0_{\O_r}(\RSprSloc)\big).
\]
The top row is exact and the right vertical map is the projection onto the first coordinate. Thus, we see that
$
    Q_{r,N} \hookrightarrow \ker(\gamma),
$
so that $\gamma$ factors as
\[
    \opT^1_{\O_r}(\calRN) \surj \opT^1_{\O_r}(\calRN)/Q_{r,N} \stackrel{\tilde{\gamma}}{\to} \opT^1_{\O_r}(\SNo),
\]
and there is a short exact sequence
\[
    0 \to Q_{r,N} \to \ker(\gamma) \to \ker(\tilde{\gamma}) \to 0.
\]
By \Cref{prop_taylorwilesprimes}(4), $Q_{r,N}$ is annihilated by a power of $p$ which is independent of $r$ and $N$. Thus, it suffices to establish the same property for $\ker(\tilde{\gamma})$. By the first diagram, we have
\begin{align*}
    \ker(\tilde{\gamma}) &\subseteq \ker\big(\opT^1_{\O_r}(\calR_N^\prime) \oplus \opT^1_{\O_r}(\RSSloc)\surj \opT^1_{\O_r}(\calR_N^\prime) \to \opT^1_{\O_r}(\SNo)  \big) \\
    &\cong \opT^1_{\O_r}(\RSSloc) \oplus \ker\big(\opT^1_{\O_r}(\calR_N^\prime)\to \opT^1_{\O_r}(\SNo)\big).
\end{align*}
We have established above that $\ker(\opT^1_{\O_r}(\calR_N^\prime)\to \opT^1_{\O_r}(\SNo))$ is annihilated by a power of $p$ that is independent of $r$ and $N$. The same holds for $\opT^1_{\O_r}(\RSSloc)$ by \Cref{prop_Tifinitelengthsmoothpoint,prop_dimension_vnmidp,theorem_semistabledeformationrings} under the genericity assumption on $\operatorname{WD}(\rho_v)$. This completes the proof.
\end{proof}
\subsection{Proof of the main theorem}
Recall that we have defined, for every $N\geq 1$,
\[
    \calC_N = \RNbar \otimesL_{\SNobar} \ON \in (\CAlgO)_{/k}.
\]
Our fixed choice of $\xi\colon \GFS \to \GL_n(\O)$ lifting $\rhobar$ determines morphisms
\[
    \xi_r\colon \calC_N \to \O_r,\ r< s(N),
\]
and \Cref{theorem_acampo} provides compatible $e_{N,M}\colon \calC_N\to \calC_M$ lifting $\xi_r$ for $M,N \in \N^\ast$.
\begin{theorem}\label{theorem_RtoCinfty}
    With $\calC_N $ as above, let 
    \[
        \Finfty = \underset{N}{\colim} \Map_{/k}(\calC_N,-) \colon \Artk \to \Ani
    \]
    be the functor considered in \Cref{section_GVextended}. There exists a morphism
    $
      \Finfty \to \DS  
    $
    inducing an equivalence 
    \[
        (\Finfty)_{(\xi)} \simeq (\DS)_{(\xi)} \in \Fun(\ArtE,\Ani).
    \]
\end{theorem}
\begin{remark}
    Recall that the inverse system $(\calC_N)_{N\in\N^\ast}$ is not canonically defined. The proof of \Cref{theorem_RtoCinfty} is also non-constructive and, a priori, the provided equivalence depends on our fixed choice of non-principal ultrafilter on $\N$.
\end{remark}
The proof of \Cref{theorem_RtoCinfty} is subdivided into separate results. The morphism in the theorem is induced by a morphism of rings
\[
    \calRS \to \calC_\infty = \lim_N \ \calC_N.
\]
We begin by showing that a homotopy class of maps $\calRS \to \calC_\infty$ is simply a compatible sequence of homotopy classes of maps $\calRS\to \calC_N$. We use the notation
\[
        [\calRS,-] \coloneqq \pi_0\big(\Map_{/k}(\calRS,-)\big).
\]
The morphisms $e_{N,M} \colon \calC_N \to \calC_M$ induce maps $[\calRS,\calC_N] \to [\calRS,\calC_M]$.
\begin{lemma}\label[lemma]{lemma_maps_to_Cinfty}
    Let $\calR \in \CNLk$ and suppose $(\calC_N)_{N\geq 1} \in \calPro(\CNLk)$ is a projective system such that $\tau_{\leq j}(\calC_N) \in \Artk$ for every $j\geq 0$, with limit $\calC_\infty = \lim \calC_N$. Then the natural maps
    \[
        [\calR,\calC_\infty] \to {\varprojlim}_N\ [\calR,\calC_N] \to {\varprojlim}_{j,N}\ [\calR,\tau_{\leq j}(\calC_N)]
    \]
    are bijections (of sets), where $\varprojlim$ denotes the 1-categorical inverse limit in $\Set$.
\end{lemma}
\begin{proof}
    Let $\calC_N^j = \tau_{\leq j}(\calC_N)$. By \cite[Lem.~A.9]{galatius_venkatesh2018}, there exists a short exact sequence of pointed sets (for any basepoint)
    \[
        0 \to {\varprojlim}^1_{j,N}\pi_{1}\big(\Map_{/k}(\calR,\calC_N^j)\big) \to [\calR,\calC_\infty] \to {\varprojlim}_{j,N}[\calR,\calC_N^j] \to 0.
    \]
    Since $\calC_N^j\in\Artk$, \Cref{lemma_artinian_mapping_spaces_finite} implies that $\Map_{/k}(\calR,\calC_N^j)$ has finite homotopy groups. Hence, by the Mittag--Leffler criterion, the first term in the sequence vanishes, and we obtain
    \[
        [\calR,\calC_\infty] \cong \varprojlim_{j,N}\ [\calR,\calC_N^j].
    \]
    It follows that $[\calR,\calC_\infty] \to {\varprojlim}_{N}[\calR,\calC_N]$ is injective. It is automatically surjective by the existence of the short exact sequence
    \[
    0 \to {\varprojlim}^1_{N}\pi_{1}\big(\Map_{/k}(\calR,\calC_N)\big) \to [\calR,\calC_\infty] \to {\varprojlim}_{N}[\calR,\calC_N] \to 0,
    \]
    also deduced from \cite[Lem.~A.9]{galatius_venkatesh2018}.
\end{proof}
\newpage
\begin{proposition}\label{prop_compatible_maps_Cinfty}
    For an integer $a\geq 0$ and $N\in \N^\ast$, let $X_N(a) \subseteq [\calRS,\calC_N]$ denote the subset consisting of classes $[f]$ such that, for every $M\leq N$ in $\N^\ast$ and $r< s(M)$,
    \begin{itemize}
        \item[(i)] $f \circ e_{N,M} \colon \calRS \to \calC_M$ lifts $\xi_r$,
        \item[(ii)] $\opT^0_{\O_r}(e_{N,M} \circ f)$ is an isomorphism,
        \item[(iii)] $p^a \coker(\opT^1_{\O_r}(e_{N,M} \circ f)) = 0$,
    \end{itemize}
    where 
\[
    \opT^i_{\O_r}(e_{N,M}\circ f) \colon \opT^i_{\O_r}(\calC_M) \to \opT^i_{\O_r}(\calR)
\]
is the induced homomorphism of $\O_r$-modules. Then we have the following:

$(X_N(a))_{N \in \N^\ast}$ is an inverse system of non-empty profinite sets with continuous transition maps induced by $e_{N,M}$, and there exists an $a\in \Z$ such that $X_N(a)$ is nonempty; in particular,
\[
    \emptyset \neq {\varprojlim}_N X_N(a) \subset [\calRS,\calC_\infty].
\]
\end{proposition}
\begin{proof}
    We subdivide the proof into four steps. In Steps 1-3, we show that $(X_N(a))$ form an inverse system of profinite sets with continuous transition maps, independently of $a\in \Z$. Finally, in Step 4, we show that there exists an $a$ such that $X_N(a) \neq \emptyset$ for all $N$; the conclusion then follows from an application of Tychonoff's theorem.

\textbf{Step 1: $(X_N(a))$ is an inverse system.}
Let $[f]\in X_N(a)$. Then, for any $L\leq M \leq N$ in $\N^\ast$, we have
\[
    [e_{M,L} \circ e_{N,M} \circ f]=[e_{N,L} \circ f]
\]
and hence $[e_{N,M} \circ f] \in X_M(a)$ and the transition maps are compatible.

\textbf{Step 2: $X_N(a)$ is profinite.} Let $\calC_N^j = \tau_{\leq j}(\calC_N)$. By \Cref{lemma_artinian_mapping_spaces_finite}, 
\[
    [\calRS,\calC_N]=\varprojlim_j\ [\calRS,\calC_N^j]
\]
is profinite. By functoriality, $e_{N,M}$ defines a morphism between truncations, also denoted $e_{N,M}$ below, and to prove that $X_N(a)$ is profinite, we show that it is the limit of an inverse system
\[
    X_N^j(a) \subseteq [\calRS,\calC_N^j],\ j\geq 0.
\]
In other words, we may represent any $[f] \in [\calRS,\calC_N]$ by a compatible sequence of classes
\[
    [f_j\colon \calRS \to \calC_N^j] \in [\calRS,\calC_N^j],\ j\geq 0,
\]
 and the claim is that $[f]\in X_N(a)$ if and only if, for every $j\geq 0$ and $r< s(M)\leq s(N)$:
\begin{itemize}
    \item[(i')] The composition $e_{N,M}\circ f_j\colon \calRS \to \calC_N^j \to \calC_M^j$ lifts $\xi_r \colon \calRS \to \O_r$.
    \item[(ii')] $\opT^0_{\O_r}(e_{N,M}\circ f_j)$ is an isomorphism.
    \item[(iii')]
    $
        p^a\coker\big(\opT^1_{\O_r}(e_{N,M}\circ f_j)\big)=0.
    $
\end{itemize}
We compare these conditions to (i), (ii) and (iii) in the statement of the theorem. Condition (i') is equivalent to (i) since $\xi_r\colon \calRS \to \O_r$ factors as 
\[
    \calRS \to \pi_0(\calR)\to\pi_0(\calC_M)\to \O_r,    
\]
and condition (ii') is equivalent to (ii) since $\opT^0_k(\tau_{\leq j}(\calC_M))\cong \opT^0_k(\pi_0(\calC_M))$ for every $j\geq 0$. To prove that (iii') is equivalent to (iii), it suffices to show that
\[
    \opT^1_{\O_r}(\calC_M^j) \cong \opT^1_{\O_r}(\calC_M^{j+1})
\]
for $j$ large enough, independent of $M$. We have a cofibre sequence of cotangent complexes
\[
    \calC_M^j \otimesL_{\calC_M^{j+1}} \opL_{\calC_M^{j+1}/\O} \to \opL_{\calC_M^{j}/\O} \to \opL_{\calC_M^{j}/\calC_M^{j+1}} \stackrel{-1}{\to}.
\]

By \cite[Prop.~3.2.16]{lurie_dag} (see also \cite[Cor.~7.4.3.6]{lurie_ha}), $\opL_{\calC_M^{j}/\calC_M^{j+1}}$ has trivial homology in degrees below $j+2$. Applying $\RHom(- \otimesL_{\calC_M^{j}} \O_r, \O_r)$ and passing to cohomology, we find that 
\[
    \opT^1_{\O_r}(\calC_M^j) \cong \opT^1_{\O_r}(\calC_M^{j+1})\ \text{ for }j\geq 2,  
\]
as needed. Thus, we define finite subsets $X_N^j(a) \subseteq [\calRS,\calC_N^j]$ using (i'), (ii') and (iii') and conclude that $X_N(a)\cong\varprojlim_j X_N^j(a)$ is profinite.

\textbf{Step 3: $X_N(a)\to X_M(a)$ is continuous.} Let $M \leq N$ in $\N^\ast$. For any $j$, we have commutative diagrams (the right one a consequence of the left one)
\[\begin{tikzcd}
	{\calC_N} & {\calC_M} & {[\calRS,\calC_N]} & {[\calRS,\calC_M]} \\
	{\calC_N^j} & {\calC_M^j} & {[\calRS,\calC_N^j]} & {[\calRS,\calC_M^j]}
	\arrow[from=1-1, to=1-2]
	\arrow[from=1-1, to=2-1]
	\arrow[from=1-2, to=2-2]
	\arrow[from=1-3, to=1-4]
	\arrow[from=1-3, to=2-3]
	\arrow[from=1-4, to=2-4]
	\arrow[from=2-1, to=2-2]
	\arrow[from=2-3, to=2-4]
\end{tikzcd}\]
It follows that $[\calRS,\calC_N]\to[\calRS,\calC_M]$ is continuous for the profinite topologies, and hence the restriction $X_N(a)\to X_M(a)$ is also continuous.

\textbf{Step 4: $X_N(a)$ is nonempty.}
Let $N\geq 1$ be arbitrary. We prove that there exists an $a\in \Z$, independent of $N$, such that $[f]\in X_N(a)$, where $f=f_N$ is the composite morphism
\[
    f \colon \calRS \simeq \calRN \otimesL_{\SNo} \O \to \RNbar \otimesL_{\SNobar} \O_N = \calC_N,
\]
where the equivalence follows from \Cref{lemma_Srings}. Note that we do not prove any compatibility between $f_N$'s for varying $N$.

For any $r< s(M)\leq s(N)$, the map $e_{N,M} \circ f$ induces a morphism of fibre sequences; passing to cohomology we obtain a commutative diagram with exact rows
\[\begin{tikzcd}
	0 & {\opT^0_{\O_r}(\calC_M)} & {\opT^0_{\O_r}(\RMbar)} & {\opT^0_{\O_r}(\SMobar)} \\
	0 & {\opT^0_{\O_r}(\calRS)} & {\opT^0_{\O_r}(\calR_N)} & {\opT^0_{\O_r}(\SNo).}
	\arrow[from=1-1, to=1-2]
	\arrow[from=1-2, to=1-3]
	\arrow[from=1-2, to=2-2]
	\arrow[from=1-3, to=1-4]
	\arrow[from=1-3, to=2-3]
	\arrow[from=1-4, to=2-4]
	\arrow[from=2-1, to=2-2]
	\arrow[from=2-2, to=2-3]
	\arrow[from=2-3, to=2-4]
\end{tikzcd}\]
The two rightmost vertical maps are isomorphisms by \Cref{prop_RNSNtangentspaces} and the natural isomorphism $\opT^0_{\O_r}(A)\cong\opT^0_{\O_r}(\pi_0(A))$, valid for any $A\in(\CNLk)_{/\O_r}$. It follows from the five lemma that the leftmost vertical map is also an isomorphism. This establishes property (ii) in the definition of $X_N(a)$, irrespective of $a$. To prove property (iii), we consider the continuation of the diagram above:
\[\begin{tikzcd}
	{\opT^0_{\O_r}(\SMobar)} & {\opT^1_{\O_r}(\calC_M)} & {\opT^1_{\O_r}(\RMbar) \oplus \opT^1_{\O_r}(\O_M)} & {\opT^1_{\O_r}(\SMobar)} \\
	{\opT^0_{\O_r}(\SNo)} & {\opT^1_{\O_r}(\calRS)} & {\opT^1_{\O_r}(\calR_N)} & {\opT^1_{\O_r}(\SNo)}.
	\arrow[from=1-1, to=1-2]
	\arrow["\cong", from=1-1, to=2-1]
	\arrow[from=1-2, to=1-3]
	\arrow["\alpha", from=1-2, to=2-2]
	\arrow[from=1-3, to=1-4]
	\arrow[from=1-3, to=2-3]
	\arrow[from=1-4, to=2-4]
	\arrow[from=2-1, to=2-2]
	\arrow[from=2-2, to=2-3]
	\arrow[from=2-3, to=2-4]
\end{tikzcd}\]
The claim to be proven is that $p^a \coker(\alpha)=0$ for an $a\geq 0$ which is independent of $r,M$ and $N$. By the diagram, we have a surjection
\[
    \ker\big(\opT^1_{\O_r}(\calR_N) \to \opT^1_{\O_r}(\SNo)\big) \cong \coker\big(\opT^0_{\O_r}(\SNo) \to \opT^1_{\O_r}(\calRS)\big)\surj \coker(\alpha).
\]
We recall that $\ker\big(\opT^1_{\O_r}(\calR_N) \to \opT^1_{\O_r}(\SNo)\big)$ has finite length bounded independently of $r$ and $N$ by \Cref{prop_T1RN}. Thus, $\coker(\alpha)$ is annihilated $p^a$ for $a$ sufficiently large and independent of $N$. Hence $[f]\in X_N(a)$ and, in particular, $X_N(a)\neq \emptyset$ for every $N$.
\end{proof}
Having established the existence of the morphism featuring in \Cref{theorem_RtoCinfty}, we turn to the task of proving that it induces an equivalence after passing to characteristic 0. 
\begin{proof}[Proof of \Cref{theorem_RtoCinfty}]
By \Cref{lemma_maps_to_Cinfty,prop_compatible_maps_Cinfty}, for $a\in \Z$ sufficiently large, there exists $[g] \in \varprojlim_N X_N(a)$ represented by a compatible sequence $(g_N)$ of morphisms 
\[
    g_N \colon\calRS \to \calC_N
\]
such that $\opT^0_{\O_r}(g_N)\colon \opT^0_{\O_r}(\calC_N)\to \opT^0_{\O_r}(\calRS)$ is an isomorphism and the cokernel
\[
    Q_{r,N}(g) \coloneqq  \coker\big(\opT^1_{\O_r}(\calC_N)\to \opT^1_{\O_r}(\calRS)\big)
\]
is $p^a$-torsion. We will prove that the corresponding morphism 
\[
    (\calF_\infty)_{(\xi)} \to (\DS)_{(\xi)}
\]
induced by $g$ is an equivalence, by showing that it is an isomorphism on $\opT^0_E$ and an epimorphism on $\opT^1_E$, then counting dimensions. The first statement follows from the fact that $g_N$ induces isomorphisms on $\opT^0$. To prove the second statement, we note that $p^a$-torsion modules are preserved under filtered colimits and inverse limits, so that
\[
    p^a Q_{r,N}(g) = 0 \implies \big(\varprojlim_r \underset{N}{\colim\ } Q_{r,N}(g)\big) \otimes_\O E = 0,
\]
where $\varprojlim$ denotes the 1-categorical inverse limit of modules. By \Cref{lemma_tgtcomplexcalculation}(2), $\colim_N \opT^1_{\O_r}(\calC_N)$ is a finite abelian group. It follows from \Cref{prop_limcolim} that we have an isomorphism
\begin{equation*}
    \begin{split}
        \big(\varprojlim_r \underset{N}{\colim\ } Q_{r,N}(g)\big) \otimes_\O E
            &\cong \coker\left(\big(\lim_r \underset{N}{\colim} \opT^1_{\O_r}(\calC_N)\big)\otimes_\O E \to \big(\lim_r \opT^1_{\O_r}(\calRS)\big)\otimes_\O E\right)\\
            &\cong \coker\bigg( \opT^1_E\big((\Finfty)_{(\xi)}\big) \to \opT^1_E\big(\RxExihat\big) \bigg)  
    \end{split}
\end{equation*}
since colimits commute and $\varprojlim$ is exact on inverse systems of finite abelian groups. Thus, this cokernel vanishes and we have proved that the map
\[
    \opT^\ast_E\big((\Finfty)_{(\xi)}\big) \to \opT^\ast_E\big(\RxExihat\big)
\]
induced by $g$ is an isomorphism in degree $0$ and an epimorphism in degree $1$. Both are graded vector spaces concentrated in degrees $[0,1]$ and have the same dimensions by \Cref{prop_S(rho),lemma_tgtcomplexcalculation}. It follows that $(\Finfty)_{(\xi)} \simeq (\DS)_{(\xi)}$ in $\Fun(\ArtE,\Ani)$.
\end{proof}

We are now in a position to piece together our results to deduce our main theorem.
\begin{proof}[Proof of \Cref{theorem_mainresult}]
\Cref{lemma_tgtcomplexcalculation}(2) provides an equivalence
\[
    \Finfty \simeq \Map_{/E}(\Rinfphat \otimesL_{\Sinfahat} E,-),
\]
and thus \Cref{cor_RSphat_equals_calRrhoss,theorem_RtoCinfty} imply that
\[
    (\calRS\otimes_\O E)_\xi^\wedge \simeq \calRrhoss \simeq \Rinfphat \otimesL_{\Sinfahat} E.
\]
By \Cref{theorem_acampo}, there exists a free $\Rinfphat$-module $\Minfphat$ such that 
\[
\Minfphat \otimesL_{\Sinfahat} E \simeq \RHom_E\big(\RG(K,\Vlt)_\p, E\big) .
\]
Thus, $\pi_\ast(\calRrhoss)$ acts freely on $H_\ast(X_K,\Vlt)_\Pi$. By \Cref{cor_extalgebra}, we obtain a free action
\[
    \pi_\ast\big(\calRrhoss\big) \cong \bigwedge H^1_g(\GFS, \ad\rho(1)) \circlearrowright H_\ast(X_K,\Vlt)_\Pi.
\]

\end{proof}

%



\bibliographystyle{alpha}
\bibliography{bib}

\end{document}